\documentclass[final]{siamltex}
\usepackage{amsmath,amsfonts,amssymb}
\usepackage{latexsym}
\usepackage{graphicx,overpic}
\usepackage{subfig}
\usepackage[table,xcdraw]{xcolor}
\usepackage{pgfplots}
\pgfplotsset{compat=newest}
\usepackage{color}
\usepackage{booktabs}
\usepackage{float}
\usepackage{multirow}

\newcommand{\R}{\mathbb R}

\newcommand{\bj}{\mathbf j}
\newcommand{\blf}{\mathbf f}
\newcommand{\bn}{\mathbf n}

\newcommand{\bx}{\mathbf x}

\newcommand{\T}{\mathcal T}

\newcommand{\divG}{{\mathop{\,\rm div}}_{\Gamma}}
\newcommand{\gradG}{\nabla_{\Gamma}}
\newcommand{\nablaG}{\nabla_{\Gamma}}
\newcommand{\laplG}{\Delta_{\Gamma}}

\newcommand{\OGamma}{\Omega^\Gamma_h}

\renewcommand{\div}{\textrm{div}\ \!}

\newcommand{\tr}{{\rm tr}}
\DeclareGraphicsExtensions{.pdf,.eps,.ps,.eps.gz,.ps.gz,.eps.Y}

\def\el {\nonumber }

\begin{document}
\title{A computational study of lateral phase separation in biological membranes}
\author{
Vladimir Yushutin\thanks{Department of Mathematics, University of Houston, Houston, Texas 77204 (yushutin@math.uh.edu).}
\and Annalisa Quaini\thanks{Department of Mathematics, University of Houston, Houston, Texas 77204 (quaini@math.uh.edu); Partially supported by NSF through grant  DMS-1620384.}
\and  Sheereen Majd\thanks{Department of Biomedical Engineering, University of Houston, Houston, Texas 77204 (smajd9@uh.edu); Partially supported by NSF through grant DMR-1753328.}
\and
Maxim Olshanskii\thanks{Department of Mathematics, University of Houston, Houston, Texas 77204 (molshan@math.uh.edu); Partially supported by NSF through grant  DMS-1717516.}
}
\maketitle
%
\begin{abstract}Conservative and non-conservative phase-field models are considered for the numerical simulation of
lateral phase  separation and coarsening in biological membranes. An unfitted finite element method is devised for these models to allow for a flexible treatment of complex shapes in the absence of an explicit surface parametrization.
For a set of biologically relevant shapes and parameter values, the paper compares the dynamic coarsening produced by conservative and non-conservative numerical models, its dependence on certain geometric characteristics and convergence to the final equilibrium.
\end{abstract}

\section{Introduction}\label{sec:intro}

Over the past 20 years, there has been a growing interest in studying phase separation in cell membranes,
mostly motivated by its important role in a variety of cellular processes.
In fact, lipid-driven separation of immiscible liquid phases is likely a factor in the formation of rafts in cell membranes \cite{simons1997functional,VEATCH20033074}. 
Lipid rafts in eukaryotic cells have been associated with important biological processes, such as
endocytosis, adhesion, signaling, and protein transport; see, e.g., \cite{Edidin2003}.
Furthermore, lipid phase separation has recently been utilized to enhance the delivery performance
of targeted lipid vesicles \cite{Bandekar_et_al2013,KARVE20104409}.
The formation of reversible phase-separated (nano)patterns (resembling lipid-rafts) on the
vesicle surface is expected to increase target selectivity, cell uptake and overall efficacy
\cite{Bandekar_et_al2013}.

Phase separation and pattern formation in lipid bilayers has been studied theoretically, experimentally,
and numerically.
For theoretical investigations of the equilibrium configuration we refer to
\cite{Andelman_et_al1992,Seifert1993,Kawakatsu_et_al1993,Kawakatsu_et_al1994,Harden_MacKintosh1994} and references therein.
Experimental studies have been mostly conducted to visualize
pattern formation on giant vesicles; see, e.g., \cite{Baumgart_et_al2003,VEATCH20033074}.
In particular, in \cite{VEATCH20033074} various dynamics have been observed with
fluorescence microscopy, including domain ripening, spinodal decomposition, and viscous fingering.
However, experimental investigations proved to be challenging due to the frail nature of
giant vesicles. Computational studies help observe dynamics and gain insights that
may not be obtained experimentally. Multicomponent vesicle have been investigated
with different numerical approaches: molecular dynamics \cite{Marrink_Mark2003,Vries_et_al2004},
dissipative particle dynamics \cite{Laradji_Sunil2004,B901866B},
and continuum based methods \cite{Wang2008,Lowengrub2009,sohn2010dynamics,nitschke_voigt_wensch_2012,Li_et_al2012,Funkhouser_et_al2014}.
In this paper, we choose a continuum based method based on phase-field description.

The phase-field method has emerged as a powerful computational approach
to modeling and predicting phase separation in materials and fluids.
It describes the system using a set of conserved and/or non-conserved field variables
that are continuous across the interfacial regions separating phases.
The main reasons for the success of the phase-field methodology are two: it replaces sharp interfaces
by thin transition regions called diffuse interfaces, making front-tracking unnecessary; and it is based on
rigorous mathematics and thermodynamics.
We consider a diffuse-interface description of phase separation developed
by Cahn and Hilliard in \cite{Cahn_Hilliard1958,CAHN1961}
and by Allen and Cahn in \cite{ALLEN1979}.
The classical Cahn--Hilliard theory for phase transformation in closed systems is
characterized by a conserved order parameter (concentration), while the Allen--Cahn equation
describes the evolution of non-conserved order parameters during phase transformation.
The Allen--Cahn and Cahn--Hilliard equations have been widely used
in many complex moving interface problems in materials science and fluid
dynamics through a phase-field approach; see, e.g.,
\cite{Lowengrub1998,Anderson_McFadden1998,boettinger2002phase,Chen2002,LIU2003}.

While there is still some controversy on the basis of lipid raft formation, function and even existence  \cite{munro2003lipid,edidin2003state}, both conservative (Cahn--Hilliard)   \cite{mcwhirter2004coupling,mercker2012multiscale,funkhouser2010dynamics,Funkhouser_et_al2014,garcke2016coupled,magi2017modelling} and non-conservative (Allen--Cahn , a.k.a. Landau--Ginzburg) \cite{ayton2005coupling,elliott2010modeling,Wang2008,ElliottStinner,elliott2013computation}  phase field models are considered in the literature  for phase separation in lipid bilayers.
By conservative it is meant that the model respects the conservation law for the phase concentration of species for any membrane surface element.
In relation to this, we note that in experimental setting, a number of molecules are known to preferentially partition into one of lipid phases on phase-separated vesicles.
Examples include membrane proteins caveolin-3, peripheral myelin protein 22 \cite{Baumgart3165,Schlebach2016}
and membrane dyes \cite{Baumgart3165,BAUMGART2007}. In such settings,
the Cahn--Hilliard equation provides the correct model. Similarly, the conservative model
seems to be suited to describe membrane separation in bacteria as well as mammalian cells.
The use of a non-conservative model (Allen--Cahn equation) may be justified
when phase separation induces ``high'' curvature to the membrane that leads to vesicle budding
or basically formation and separation of individual vesicles from the original parent vesicle
\cite{Baumgart_et_al2003,HURLEY2010875}. We also note that a phase-field model with global volume constraints was proposed in   \cite{du2006simulating,ALAND201432} for the modeling of an impermeable closed membrane embedded in an incompressible fluid. Finally, for additional models to simulate mobile lipid rafts we refer to
\cite{PhysRevLett.105.148102,Camley_Brown,C2CP41274H,Adkins:2017}.

In this paper, we apply and critically compare both surface Cahn--Hilliard and surface Allen--Cahn equations for  the numerical  simulation of lateral  phase separation in biological membranes.
Some simplifying assumptions are made in the present  study. First, we consider rigid shapes thus ignoring the coupling of lateral sorting to radial deformations caused by the minimization of elastic or bending energy (cf., e.g.,  \cite{Wang2008,ElliottStinner,funkhouser2010dynamics,elliott2013computation,Funkhouser_et_al2014}). Although most bio-membranes are compliant, the assumption is reasonable for some manufactured  vesicles  designed for intracellular drug delivery~\cite{Laouini_et_al2012}.
Second, we do not account here for possible viscous dissipation and fluidity of cell membrane (see, e.g., \cite{magi2017modelling} for Cahn--Hilliard--Navier--Stokes surface models). Finally, we neglect the effect of an external fluid; see, e.g., \cite{sohn2010dynamics} for modelling of a multicomponent vesicle in an external
viscous flow.

The Allen--Cahn and Cahn--Hilliard equations are challenging to solve numerically due to non-linearity,
and stiffness. The Cahn--Hilliard has the added difficulty of a
fourth order derivative in space. Although there exists an extensive literature on numerical methods for 
these models in planar and volumetric  domains (see, e.g., recent publications~\cite{guillen2014second,tierra2015numerical,liu2015stabilized, cai2017error, hou2017discrete} and references therein), there are not so many papers where the equations are treated on surfaces.
Solving equations numerically on general surfaces poses additional difficulties that are related to the discretization of tangential differential operators and the approximate recovery of complex  shapes.
Several authors have opted for a finite difference method.
For example,  the closest point finite difference  method  was applied
to solve the surface Allen--Cahn equation in \cite{kim2017finite} and  the
Cahn--Hilliard equation on torus  in \cite{gera2017cahn} and in \cite{jeong2015microphase} on more general domains.
A finite difference method for a diffuse volumetric representation of the surface Cahn--Hilliard equation was introduced in \cite{greer2006fourth}.
A finite difference method in a reference domain was used in \cite{Funkhouser_et_al2014}
to model phase separation with radial deformation on sphere like domains.
However, a finite element method (FEM) is often considered to be the most flexible numerical approach to handle complex geometries.

Finite element solutions to the surface Allen--Cahn equation can be found in \cite{dziuk2007surface,elliott2010modeling,elliott2013computation}.
The convergence of a FEM for the surface Cahn--Hilliard equation
was studied in \cite{du2011finite}, where numerical examples
for a sphere and saddle surface are provided; the Cahn--Hilliard equation on more general surfaces was treated by a finite element method more recently in \cite{garcke2016coupled,li2017unconditionally}.
Surface FEMs were extended to more general systems
in several papers: in \cite{nitschke_voigt_wensch_2012} solutions to the  surface
Cahn--Hilliard--Navier--Stokes equation were computed  on a sphere and torus;
the authors of \cite{barrett2017finite} studied a finite element method for the bulk Navier--Stokes equations coupled to the surface Cahn--Hilliard model; and
in \cite{elliott2015evolving} the Cahn--Hilliard equation was 
solved numerically by a finite element method on evolving surfaces.
All of the above references use a sharp surface representation and
a discretization mesh \emph{fitted} to the computational surface.

In the present paper, we study for the first time a \emph{geometrically unfitted} finite element method
for the simulation of lateral phase separation on surfaces.
Our approach builds on earlier work on a unfitted FEM for elliptic PDEs posed on surfaces~\cite{ORG09}
called TraceFEM. Unlike some other geometrically unfitted methods for surface PDEs,
TraceFEM employs sharp surface representation. The surface can be defined implicitly and no knowledge of the surface parametrization is required. This approach allows for flexible numerical treatment of complex shapes,
like the ones found in cell biology.
After validating  the accuracy of  the numerical method with benchmark problems,
we apply it to simulate phase transition on a series of surfaces of increasing
geometric complexity using both the surface Allen--Cahn and Cahn--Hilliard models. The
surfaces are chosen to resemble shapes of biological membranes known to
exhibit phase separation.
We compare the dynamic coarsening produced by conservative and
non-conservative numerical models, its interplay with the curvature and convergence to the final equilibrium.
An advantage of our approach is that it can handle ``small'' (of the order of 1\% of the shape characteristic length, which is consistent with available experimental observations)
interface thicknesses between phases without stability issues, as long as the time step is properly set.
Finally, we note that TraceFEM can be naturally combined with a level-set surface representation
and works  well for  PDEs posed on evolving surfaces, including cases with topological changes (cf., e.g., \cite{lehrenfeld2018stabilized}).
We will not make use of this rather unique property in this paper. However,
such property will become very convenient for the numerical simulation of certain phenomena
like uptake of drug carriers that use reversible phase-separated (nano)patterns on the vesicle surface
\cite{Bandekar_et_al2013,KARVE20104409}.

The remainder of the paper is organized as follows. In Sec.~\ref{s_cont}
we state the surface Allen--Cahn and the Cahn--Hilliard equations and the respective variational
formulations.
The application of TraceFEM to both models is presented in Sec.~\ref{sec:num_method}.
After a validation of our implementation of TraceFEM, in Sec.~\ref{sec:num_res}
we study phase separation modeled by both the Allen--Cahn and the Cahn--Hilliard equations
on different surfaces of increasing complexity.
Concluding remarks are provided in Sec.~\ref{sec:concl}.

\section{Mathematical model}\label{s_cont}
In order to state the surface Allen--Cahn and the Cahn--Hilliard equations, we first need to introduce some notation.

Let $\Gamma$ be a closed sufficiently smooth surface in $\mathbb{R}^3$, with
the outward pointing unit normal $\bn$. The surface represents a 
bio-membrane, e.g., a closed bilayer composed of multiple lipids, embedded in a bulk fluid.
In the present study, we do not account for  bulk phenomena and concentrate on the lateral behavior of the system.
For any sufficiently smooth function $g$ in a neighborhood of $\Gamma$ its tangential gradient is
defined as $\nablaG g=\nabla g-(\bn\cdot\nabla g)\bn$. The tangential (surface) gradient $\nablaG g$ then depends only on values of $g$ restricted to $\Gamma$ and $\bn\cdot\nablaG g=0$ holds. For a vector field $\blf$ on $\Gamma$ we define $\nablaG \blf$ componentwise.
We need the surface divergence operator for $\blf$ and
the Laplace-Beltrami operator for $g$ :
\[
\divG \blf  := \tr (\gradG \blf),  \quad \laplG g := \divG (\gradG g),
 \]
where $\tr(\cdot)$ is the trace of a matrix.

Further $L^2(\Gamma)$ is the Lebesgue space of square-integrable functions on $\Gamma$ and $H^1(\Gamma)$ is the Sobolev space of all functions $g\in L^2(\Gamma)$ such that $\nabla_\Gamma g\in L^2(\Gamma)^3$.

\subsection{Allen--Cahn equation}\label{sec:AC}

Let $\eta \in [0, 1]$ be an order parameter on $\Gamma$, i.e. a measure of the degree of order across the boundaries in a phase transition system, with $\eta = 0$ indicating complete lack of order and
 $\eta = 1$ indicating full order.
Such parameter normally ranges between zero in one phase 
and nonzero in the other.
For example, if solid-liquid phase transition happens on $\Gamma$, then the order parameter is the difference of the surface material densities.
If the total specific free energy $f$ is not at a minimum with respect to a local variation in $\eta$,
Allen and Cahn postulated in \cite{ALLEN1979} that there is an immediate change in $\eta$ given
by:
\begin{align}\label{eq:evol_eta}
\frac{\partial \eta}{\partial t} = - \alpha \frac{\delta f}{\delta \eta} \quad \text{on}~\Gamma \times (0, T],
\end{align}
where $\alpha$ is a positive kinetic coefficient and $(0, T]$ is the time interval of interest.
The total specific (i.e., per unit surface)  free energy $f$
is a function of the order parameter:
\begin{align}\label{eq:total_free_e_eta}
f(\eta) = f_0(\eta) + \frac{1}{2} \epsilon^2 | \gradG \eta |^2.
\end{align}
In \eqref{eq:total_free_e_eta}, $\epsilon$ is the gradient energy coefficient and $f_0$ is
the specific free energy of a homogeneous phase, which
is a function of $\eta$ with a characteristic double-well
form.
The kinetic equation \eqref{eq:evol_eta} reflects the fact that the order parameter $\eta$ is not
a conserved quantity.
The 
functional derivative of $f$
with respect to $\eta$ 
 is given by:
\begin{align}\label{eq:var_der_f}
\frac{\delta f}{\delta \eta} = f'_0(\eta) - \epsilon^2  \laplG \eta.
\end{align}
Plugging \eqref{eq:var_der_f} into eq.~\eqref{eq:evol_eta}, we obtain:
\begin{equation}\label{eq:AC_eq}
\frac{\partial \eta}{\partial t} + \alpha  f'_0(\eta) - \alpha \epsilon^2 \laplG \eta = 0 \quad \text{on}~\Gamma\times (0, T].
\end{equation}
We remark that $\alpha \epsilon^2$ in eq.~\eqref{eq:AC_eq} has dimensions of a diffusion coefficient, i.e.
m$^2$/s. Obviously, eq.~\eqref{eq:AC_eq} needs to be supplemented with initial condition $\eta = \eta_0$
on $\Gamma\times \{0\}$, for a given $\eta_0$.
A classical choice for $f_0$ is the Ginzburg--Landau double-well potential
\begin{align}
f_0(\eta) = \frac{1}{4} (\eta -1)^2\eta^2, \el
\end{align}
which makes eq.~\eqref{eq:AC_eq} non-linear.

For the numerical method, we need a weak (variational) formulation of the surface Allen--Cahn equation. To devise it, one multiplies \eqref{eq:AC_eq} by $v\in H^1(\Gamma)$, integrates over $\Gamma$ and employs the integration by parts identity.
For a closed smooth surface $\Gamma$, the integration by parts identity reads:
 \begin{equation}\label{int_parts}
 \int_{\Gamma} v\div_\Gamma\blf \, ds= -\int_{\Gamma} \blf\cdot\nabla_\Gamma v\, ds+\int_{\Gamma} \kappa v\blf\cdot\bn \, ds,\quad \text{for}~\blf\in H^1(\Gamma)^3,\, v\in H^1(\Gamma),
 \end{equation}
here $\kappa$ is the sum of principle curvatures. Identity \eqref{int_parts} is applied to the diffusion term
in \eqref{eq:AC_eq}, i.e. $\blf=\nabla_\Gamma \eta$, which makes the curvature term vanish.
  This leads to the week formulation: Find $\eta \in H^1(\Gamma) $ 
such that
\begin{equation}\label{eq:AC_weak}
\int_\Gamma \left(\frac{\partial \eta}{\partial t} + \alpha f_0'(\eta)\right) \,v \, ds + \alpha \epsilon^2\int_\Gamma  \gradG \eta \, \gradG v \, ds = 0,\quad \forall\, v \in H^1(\Gamma).
\end{equation}

\subsection{Cahn--Hilliard equation}\label{sec:CH}

On $\Gamma$ we consider a heterogeneous mixture of two species with mass
concentrations $c_i = m_i/m$, $i = 1, 2$, where $m_i$ are the masses of the components
and $m$ is the total mass.  Since $m = m_1 + m_2$,
we have $c_1 + c_2 = 1$. Let $c_1$ be the representative concentration $c$, i.e. $c = c_1$.
Unlike the order parameter $\eta$ in Section~\ref{sec:AC}, concentration $c \in [0, 1]$
is a conserved quantity.
Moreover let $\rho$ be the constant total density of the system $\rho = m/S$, where
$S$ is the surface area of $\Gamma$.
Phase separation in this two component system
can be modelled by the Cahn--Hilliard equation \cite{Cahn_Hilliard1958,CAHN1961}.

In order to describe the evolution of the concentration profile $c(\bx, t)$, we consider
the conservation law:
\begin{align}\label{eq:evol_c}
\rho \frac{\partial c}{\partial t} +  \divG \bj = 0 \quad \text{on}~\Gamma  \times (0, T],
\end{align}
where $\bj$ is a diffusion flux. The flux $\bj$ is defined according to (empirical) Fick's law:
\begin{align}\label{eq:flux}
\bj = - M \gradG \mu \quad \text{on}~\Gamma, \quad \mu = \frac{\delta f}{\delta c},
\end{align}
where $M$ is the so-called mobility coefficient (see \cite{Landau_Lifshitz_1958}) and
$\mu$ is the chemical potential, which is defined as the functional derivative of the total
specific free energy $f$ with respect to the concentration $c$. Thus, we introduce
the total specific free energy:
\begin{align}\label{eq:total_free_e}
f(c) = f_0(c) + \frac{1}{2} \epsilon^2 | \gradG c |^2.
\end{align}
Just like in eq.~\eqref{eq:total_free_e_eta}, $f_0(c)$ is the free energy per unit surface, while the second term represents
the interfacial free energy based on the concentration gradient.
We recall that in order
to have phase separation, $f_0$ must be a non-convex function of $c$.
A fundamental fact of the chemical thermodynamics 
is that even when phase separation has occurred, there
is a limited miscibility between the components. In model
\eqref{eq:evol_c}--\eqref{eq:total_free_e},
the interface between the two components is a layer of size ${\epsilon}$
where thermodynamically
unstable mixtures are stabilized by a gradient term in the energy.
Further details concerning the thermodynamics of
partially miscible mixtures can be found, for example, in \cite{Landau_Lifshitz_1958}.

By combining eq.~\eqref{eq:evol_c}, \eqref{eq:flux}, and \eqref{eq:total_free_e}, we obtain
the surface Cahn--Hilliard equation:
\begin{align}\label{eq:CH}
\rho \frac{\partial c}{\partial t} -  \divG \left(M \gradG \left(f_0' - \epsilon^2 \laplG c\right)\right) = 0 \quad \text{on}~\Gamma.
\end{align}
Eq.~\eqref{eq:CH} is a fourth-order equation, so casting it in a weak form
would result in the presence of second-order spatial derivatives. From the numerical point of view,
it is beneficial to avoid higher order spatial derivatives. Hence, it is common to rewrite eq.~\eqref{eq:CH} in mixed form, i.e. as
two coupled second-order equations:
\begin{align}
&\rho \frac{\partial c}{\partial t} -  \divG \left(M \gradG \mu \right)  = 0 \quad \text{on}~\Gamma \label{eq:sys_CH1}, \\
&\mu = f_0' - \epsilon^2 \laplG c \quad \text{on}~\Gamma. \label{eq:sys_CH2}
\end{align}

System \eqref{eq:sys_CH1}--\eqref{eq:sys_CH2} needs to be supplemented with the
definitions of mobility $M$ and free energy per unit surface $f_0$. A possible choice for $M$
is given by
\begin{align}\label{eq:M}
M = M(c) = c(1 - c).
\end{align}
This mobility is referred to as a degenerate mobility, since it is not strictly positive.
We note that in many of the existing analytic studies, as well as numerical simulations, mobility is assumed to be constant.
At the same time, concentration dependent mobility was already considered by Cahn
\cite{CAHN1961} and \eqref{eq:M} is also a popular choice for numerical studies.
Although it is known that the dependence between the mobility and the concentration difference
produces important changes during the coarsening process, only a few authors consider more complex
mobility functions; see, e.g., \cite{PhysRevE.60.3564}. In the absence of studies on the
appropriate mobility function for lateral phase separation in biological membranes, here we choose to use \eqref{eq:M}.
Again, a common choice for $f_0$ is given by
\begin{align}\label{eq:f0}
f_0(c) = \frac{\xi}{4} c^2(1 - c)^2,
\end{align}
where $\xi$ defines the barrier height, i.e. the local maximum at $c = 1/2$ \cite{Emmerich2011}. We set  $\xi = 1$ for the rest of the paper.
With mobility as in \eqref{eq:M} and
specific free energy as in \eqref{eq:f0}, problem \eqref{eq:sys_CH1}--\eqref{eq:sys_CH2} is a coupled system of nonlinear PDEs posed on $\Gamma$.

Following a procedure similar to the one used to get eq.~\eqref{eq:AC_weak}, we obtain the weak (variational) formulation of the surface Cahn--Hilliard problem \eqref{eq:sys_CH1}--\eqref{eq:sys_CH2}: Find $(c,\mu) \in H^1(\Gamma) \times H^1(\Gamma)$ 
such that
\begin{align}
&\int_\Gamma \rho \frac{\partial c}{\partial t} \,v \, ds + \int_\Gamma M \gradG \mu \, \gradG v \, ds = 0, \label{eq:sys_CH1_weak} \\
&\int_\Gamma  \mu \,q \, ds - \int_\Gamma f_0'(c) \,q \, ds - \int_\Gamma \epsilon^2 \gradG c \, \gradG q \, ds = 0, \label{eq:sys_CH2_weak}
\end{align}
for all $ (v,q) \in H^1(\Gamma) \times H^1(\Gamma)$.

\section{Numerical method}\label{sec:num_method}

Biological membranes exhibit asymmetric complex sha\-pes, which may affect the phase separation on the bilayer in an intricate way. To model the process numerically, we apply the trace finite element method~\cite{olshanskii2017trace}. This method allows to solve for a scalar quantity or a vector field on the surface $\Gamma$ and does not require parametrization or triangulation of $\Gamma$. To discretize surface equations, TraceFEM relies on a tessellation of a bulk computational domain $\Omega$ ($\Gamma\subset\Omega$ holds) into shape regular tetrahedra untangled to the position of $\Gamma$.

Biological membranes also exhibit shape transitions and shape instabilities.
Thus, a realistic model of phase separation on biological membranes requires
to solve either the Allen--Cahn or Cahn-Hilliard equation
on evolving shapes. This is our longer term goal, which dictates the choice of the numerical approach.
In fact, as already mentioned in Sec.~\ref{sec:intro}, evolving surfaces with possible topology changes
can be relatively easily handled in the present numerical framework \cite{lehrenfeld2018stabilized}.

Let the 3D bulk domain $\Omega$ be such that $\Gamma\subset\Omega$. Surface
$\Gamma$ is defined implicitly as the zero level set of a sufficiently smooth (at least Lipschitz continuous) function $\phi$, i.e. $\Gamma=\{\bx\in\Omega\,:\, \phi(\bx)=0\}$, such that $|\nabla\phi|\ge c_0>0$ in a 3D neighborhood $U(\Gamma)$ of the surface. The vector field $\bn=\nabla\phi/|\nabla\phi|$ is normal on $\Gamma$ and  defines  quasi-normal directions in $U(\Gamma)$.
Let $\T_h$ be  the collection of all tetrahedra, such that $\overline{\Omega}=\cup_{T\in\T_h}\overline{T}$.
The subset of tetrahedra that have a \emph{nonzero intersection} with $\Gamma$ is denoted by $\T_h^\Gamma$.
The grid is refined towards $\Gamma$, however the tetrahedra from $\T_h^\Gamma$ form a quasi-uniform tessellation with the characteristic tetrahedra  size  $h$.
The domain formed by all tetrahedra in $\T_h^\Gamma$ is denoted by $\OGamma$.

\begin{figure}[h]
\centering
\includegraphics[width=.4\textwidth]{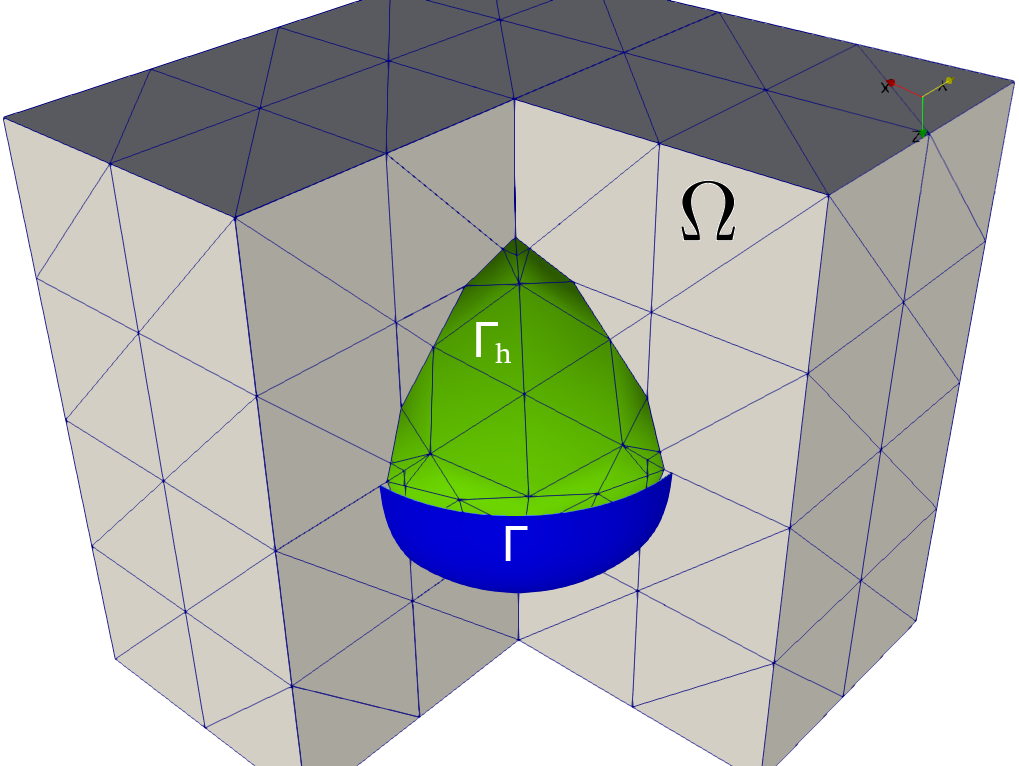}
  \caption{{Example of bulk domain $\Omega$ with a given triangulation and sphere surface $\Gamma$ with the
  corresponding ``discrete'' surface $\Gamma_h$.}}
  \label{fig:Gamma_h}
\end{figure}

On $\T_h^\Gamma$ we use a standard finite element space of continuous functions that are piecewise-polynomials
of degree $1$.
Higher order finite elements are possible (see, e.g., \cite{grande2018analysis}) but will not
be addressed in this paper.
This bulk (volumetric) finite element space is denoted by $V_h$,
\[
V_h=\{v\in C(\OGamma)\,:\, v\in P_1(T)~\text{for any}~T\in\T_h^{\Gamma}\}.
\]
For the purpose of numerical integration, we approximate $\Gamma$ with a ``discrete'' surface $\Gamma_h$, which is defined as the zero level set of a $P_1$ Lagrangian interpolant for $\phi$ on the one time refined mesh:
\[
\Gamma_h=\{\bx\in\Omega\,:\, \phi_h(\bx)=0,\quad \phi_h:=I_{h/2}(\phi(\bx))\in V_{h/2}\}, \quad \text{with}~\bn_h=
\frac{\nabla I^2_{h/2}(\phi)}{|\nabla I^2_{h/2}(\phi)|},
\]
where $I^2_{h/2}(\phi)$ is a $P_2$ nodal interpolant of the level set function. See Fig.~\ref{fig:Gamma_h}.

For time-stepping we apply a semi-implicit stabilized schemes from~\cite{Shen_Yang2010}.
At time instance $t^k=k\Delta t$, with time step $\Delta t=\frac{T}{N}$,  $\eta^k$ denotes the approximation of the order parameter $\eta(t^k, \bx)$;  similar notation is used for other quantities of interest. Further, we need second order approximation of first and second time derivatives \cite{Shen_Yang2010}:
 \begin{equation}\label{BDF2}
\left[\eta\right]_t^{k} =\frac{3\eta^{k}-4\eta^{k-1}+\eta^{k-2}}{2\Delta t},\quad
\left[\eta\right]_{tt}^{k} =\frac{\eta^{k+1}-2\eta^{k}+\eta^{k-1}}{|\Delta t|^2},
\end{equation}
and linear extrapolation for $f_0'$ at time $t^k$: $\widetilde f_0'(\eta)^k=2f_0'(\eta^{k-1})- f_0'(\eta^{k-2})$. Same notations for differences and extrapolation will be used with other variables.

The finite element discretizations are based on the weak
formulations of  the surface Allen--Cahn and Cahn--Hilliard problems
\eqref{eq:AC_weak} and \eqref{eq:sys_CH1_weak}--\eqref{eq:sys_CH2_weak}, respectively.
For the Allen--Cahn equation, the semi-implicit stabilized TraceFEM reads: Given $\eta^{k-1}_h,\eta^{k-2}_h\in V_h$, find $\eta_h^k \in V_h$ solving
\begin{equation}\label{eq:AC_FE}
\begin{split}
\int_{\Gamma_h}\left\{\left[\eta_h\right]_t^{k} + \beta_s|\Delta t|^2\left[\eta_h\right]_{tt}^{k-1}\right. &\left. + \alpha  \widetilde f'_0(\eta_h)^{k}\right\}v_h\,ds + \alpha\epsilon^2 \int_{\Gamma_h} \nabla_{\Gamma_h}\eta_h\cdot\nabla_{\Gamma_h}v_h\,ds\\
&+\epsilon^2 h \int_{\OGamma} (\bn_h\cdot\nabla \eta_h^k) (\bn_h\cdot\nabla v_h) \, dx = 0
\end{split}
\end{equation}
for all  $v_h \in V_h$, $k=2,3,\dots,N$. For $k=1$ an obvious first order modification is used.
The last term in \eqref{eq:AC_FE} is included to stabilize the {resulting algebraic systems}~\cite{grande2018analysis}. The term is consistent up to geometric errors related to the approximation of $\Gamma$ by $\Gamma_h$ and $\bn$ by $\bn_h$ in the following sense: any smooth  solution $\eta$ of equations~\eqref{eq:AC_FE} can be always extended off the surface along (quasi)-normal directions so that $\bn\cdot\nabla \eta=0$ in $\OGamma$.  Another stabilization term $|\Delta t|^2\left[\eta_h\right]_{tt}^{k-1}$ scaled with user defined parameter $\beta_s$ is included to relax the stability restriction for the time step $\Delta t$. It introduces the consistency error of second order in time. We set $\beta_s =1$. 

Likewise, the semi-implicit stabilized TraceFEM for the Cahn--Hilliard equations
{(see Sec.~3.1.2 in \cite{Shen_Yang2010})}
reads: Given $c^{k-1}_h,c^{k-2}_h\in V_h$ and $\mu^{k-1}_h,\mu^{k-2}_h\in V_h$,  find $c_h^k, \mu^k \in V_h$ solving
\begin{equation}\label{eq:CH_FE}
\begin{split}
\int_{\Gamma_h} \rho \left[c_h\right]_t^{k}v_h\,ds    +  \int_{\Gamma_h} M(\widetilde c^k) \nabla_{\Gamma_h}\mu_h^k\cdot\nabla_{\Gamma_h}v_h\,ds +  h \int_{\OGamma} (\bn_h\cdot\nabla \mu_h^k) (\bn_h\cdot\nabla v_h) \, dx\\
+\int_{\Gamma_h}\left\{  \mu_h^k - \beta_s|\Delta t|^2\left[{c}_h\right]_{tt}^{k-1} - \widetilde  f'_0( c_h)^{k}\right\}q_h\,ds
- \int_{\Gamma_h} \epsilon^2 \gradG c_h^k \, \gradG q_h \, ds \\ -\epsilon^2 h \int_{\OGamma} ( \bn_h\cdot\nabla c_h^k) (\bn_h\cdot\nabla q_h) \, dx = 0
\end{split}
\end{equation}
for all  $v_h\in V_h$ and $q_h \in V_h$, $k=2,3,\dots,N$. Again, for $k=1$ an obvious first order modification is used and we set $\beta_s=1$.

The Allen--Cahn and Cahn--Hilliard equations define gradient flows of the energy functional
$
E(u)=\int_{\Gamma}f(u)\,ds
$
 in $L^2(\Gamma)$ and $H^{-1}(\Gamma)$ (a dual space to $H^1(\Gamma)$), respectively. More precisely, the following energy minimization properties hold:
\begin{align}
  \frac{d}{dt}E(\eta) &= -\int_{\Gamma}\alpha|\epsilon^2\Delta_\Gamma\eta+f'(\eta)|^2\,ds<0\quad\text{Allen--Cahn} \label{E_AC} \\
  \frac{d}{dt}E(c) &= -\int_{\Gamma}|M(c)\nabla_\Gamma(\epsilon^2\Delta_\Gamma c+f'(c))|^2\,ds<0\quad\text{Cahn--Hilliard} \label{E_CH}
\end{align}

It is important for a physically consistent numerical method to adhere to fundamental properties
\eqref{E_AC} and \eqref{E_CH}.
According to \cite{Shen_Yang2010}, the semi-discrete counterpart of~\eqref{eq:AC_FE} in a planar domain  is stable with (slightly) modified potential $\hat{f}_0$ such that $\max_{x\in\R}|\hat{f}_0^{''}|\le L$  and under condition $\Delta t \le \frac23L^{-1}$.
For the Cahn-Hillard problem, the stability  of a semi-discrete plain counterpart of~\eqref{eq:CH_FE}  was shown  for
\begin{equation}\label{eq:CH_cond}
\Delta t \le \epsilon^2L^{-2}
\end{equation}
and $M(c)=1$.
A straightforward extension of arguments from \cite{Shen_Yang2010} proves that the finite element solution to \eqref{eq:AC_FE} is stable under the same condition and it holds:
\begin{equation*}
E_h(\eta^{k+1}_h)\le E_h(\eta^k_h)\quad\text{for all}~k=1,2,\dots,
\end{equation*}
with the numerical  energy functional,
\begin{equation*}
E_h(\eta^k)= \int_{\Gamma_h} \alpha (f(\eta^k)+\beta_s|\Delta t|^2|\left[\eta\right]_{\hat t}^{k}|^2)\,ds+ \epsilon^2 h \int_{\OGamma} \left|\bn_h\cdot\nabla \eta^k\right|^2 \, dx,
\end{equation*}
here $\left[\eta\right]_{\hat t}^{k}:=(\eta^k-\eta^{k-1})/\Delta t$ is the first order finite difference derivative.
Extension of the stability analysis for the fully discrete method  \eqref{eq:CH_FE} with a concentration dependent mobility coefficient is less straightforward and we shall address it elsewhere.
We remark that stability conditions are independent of $\beta_s$. It is, however, noted in \cite{Shen_Yang2010} that in practice the restriction for time step is much less severe if stabilization parameter $\beta_s$ is not too small.
This also agrees with our numerical experience; see Sec.~\ref{sec:beta_s}.

\section{Numerical experiments}\label{sec:num_res}

After validating the accuracy of the numerical method outlined in
Sec.~\ref{sec:num_method},
a series of numerical tests is presented to study phase separation modeled
by the Allen--Cahn and Cahn--Hilliard equations on surfaces of biological interest
and to demonstrate the flexibility of our approach. We use
free energy per unit surface \eqref{eq:f0} with $\xi = 1$  for both equations
for the sake of comparing the evolution of phase separations.

We start by comparing the numerical results produced by the two models
on a sphere in Sec.~\ref{sec:sphere}. We find that the Cahn--Hilliard model
successfully reproduces the spinodal decomposition experimentally
observed in giant vesicles in \cite{VEATCH20033074}. Next, in Sec.~\ref{sec:spindle}
we compare both models on the surface of a spindle
with the aim of getting a preliminary insight into the formation of microdomains
in bacteria \cite{Bramkamp_Lopez2015}. Finally, we present in Sec.~\ref{sec:scell} the results on
a more complex surface that represents an idealized cell. For both the sphere and the idealized cell,
we let the simulations run until sufficiently close to the steady state. All implementation are done in the FE package DROPS~\cite{DROPS}.

\subsection{Validation of the numerical method} \label{sec:NumVal}
Before presenting the results on phase separation,
we proceed with checking the spatial accuracy of the finite element method described in Sec.~\ref{sec:num_method}
with a benchmark test. The aim is to validate our implementation of the method.
For this purpose, we consider the following synthetic solution to
the Allen--Cahn equation on the unit sphere, centered at the origin:
\begin{align*}
\eta^*={\frac{1}{2}\left(1-0.8e^{-0.4t}\right)\left(Y(x_1,x_2)+ 1\right), \quad Y(x_1,x_2)=x_1 x_2}, \quad t\in[0,5].
\end{align*}
Here and in the following, $\bx = (x_1, x_2, x_3)^T$ denotes a point in $\mathbb{R}^3$.
Thus,  $\eta^*$ is the exact solution to the non-homogeneous equation $
\eta_t + \alpha  f'_0(\eta) - \alpha \epsilon^2 \laplG \eta = g$,
with free energy  \eqref{eq:f0} and non-zero right-hand side $g$, which is easy to compute since $Y$ is a real spherical harmonic function. We set $\alpha = 1$ and $\epsilon = 0.1$.

To apply the method, we characterize the surface $\Gamma$ as the zero level set of function $\phi(\bx) = \|\bx\|_2 -1$,
and embed $\Gamma$ in an outer cubic domain $\Omega=[-5/3,5/3]^3$.
The initial triangulation $\T_{h_\ell}$ of $\Omega$ consists of 8 sub-cubes,
where each of the sub-cubes is further subdivided into 6 tetrahedra.
Further, the mesh is refined towards the surface, and $\ell\in\Bbb{N}$ denotes the level of refinement, with the associated mesh size $h_\ell= \frac{10/3}{2^{\ell+1}}$. The refined tetrahedra cut by the surface form the computational mesh $\T^\Gamma_{h_\ell}$.
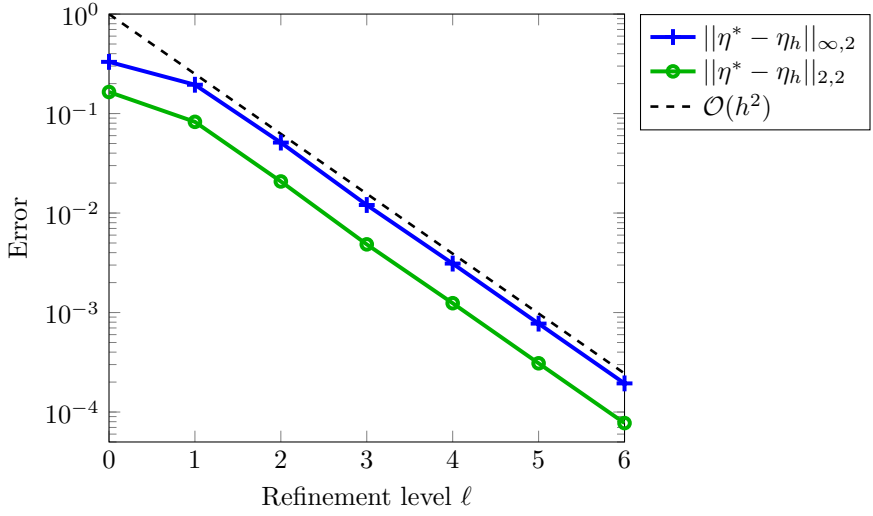
\begin{figure}
	\begin{center}
		\begin{tikzpicture}
		\def\vara{1}
		\def\varb{1}
		\begin{semilogyaxis}[ xlabel={Refinement level $\ell$},xmax=6,xmin=0, ylabel={Error}, ymin=5E-5
		, ymax=1, legend style={ cells={anchor=west}, legend pos=outer north east} ]

		\addplot+[blue,mark=+, solid,mark size=3pt,mark options={solid},line width=1.5pt] table[x=level, y=C_chi] {test11_h-t.dat};
		
		\addplot+[black!30!green,mark=o, solid,mark size=2pt,mark options={solid},line width=1.5pt] table[x=level, y=t_L2_chi] {test11_h-t.dat};
		
		%
		
		\addplot[dashed,line width=1pt] coordinates { 
			(0,\varb) (1,\varb*0.5*0.5) (2,\varb*0.25*0.25) (3,\varb*0.125*0.125) (4,\varb*0.0625*0.0625) (5,\varb*0.03125*0.03125) (6,\varb*0.03125*0.03125*0.5*0.5)
		};
		\legend{
			$|| \eta^* - \eta_h ||_{\infty,2}$, 
			$|| \eta^* - \eta_h ||_{2,2}$, 
			$\mathcal{O}(h^{2})$}
		\end{semilogyaxis}
		\end{tikzpicture}
		\caption{{Discrete $L_\infty(0,T; L^2(\Gamma_h))$ norm \eqref{eq:i2norm} and discrete $L_2(0,T; L^2(\Gamma_h))$
			norm \eqref{eq:22norm} of the error for order parameter $\eta$
			plotted against the refinement level $\ell$ 
			with $\Delta t=2^{-(1+\ell)}$.}
		}
		\label{test7_convergence}
	\end{center}
\end{figure}

We report in Fig.~\ref{test7_convergence} the discrete $L_\infty(0,T; L^2(\Gamma_h))$ norm:
\begin{align}\label{eq:i2norm}
|| \eta^* - \eta_h ||_{\infty, 2}= \max_k\|\eta^*(t_k) - \eta_h(t_k)\|_{L^2(\Gamma_h)}
\end{align}
and the discrete  $L_2(0,T; L^2(\Gamma_h))$ norm:
\begin{align}\label{eq:22norm}
|| \eta^* - \eta_h ||_{2,2}= \left(\frac{1}{T}\sum_k\, \Delta{}t\|\eta^*(t_k) -\eta_h(t_k)\|^2_{L^2(\Gamma_h)}\right)^{1/2}
\end{align}
of the error for order parameter $\eta$ plotted against the refinement level $\ell$. The time step was refined together with the mesh size according to $\Delta t=2^{-(1+\ell)}$. Norms in \eqref{eq:i2norm}--\eqref{eq:22norm} naturally appear in the error analysis for Allen--Cahn equation; see \cite{Shen_Yang2010}. The second order convergence observed in  Fig.~\ref{test7_convergence} is optimal for $P^1$ finite elements and consistent with the second order time-stepping method in \eqref{eq:AC_FE}.
We note that all the norms reported in Fig.~\ref{test7_convergence} are
computed on the approximate surface $\Gamma_h$, where $\eta^*$ was defined through its normal extension from $\Gamma$.

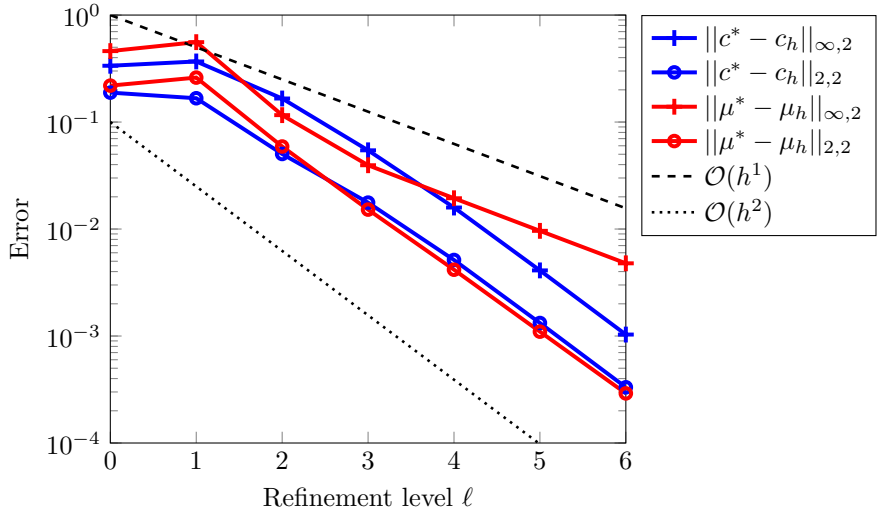
\begin{figure}
	\begin{center}
		\begin{tikzpicture}
		\def\vara{1}
		\def\varb{0.1}
		\begin{semilogyaxis}[ xlabel={Refinement level $\ell$},xmax=6,xmin=0, ylabel={Error}, ymin=1E-4
		, ymax=1, legend style={ cells={anchor=west}, legend pos=outer north east} ]
		
		\addplot+[blue,mark=+, solid,mark size=3pt,mark options={solid},line width=1.5pt] table[x=level, y=C_chi] {test10_h-t.dat};
		
		\addplot+[blue,mark=o, solid,mark size=2pt,mark options={solid},line width=1.5pt] table[x=level, y=t_L2_chi] {test10_h-t.dat};
		
		\addplot+[red,mark=+, solid,mark size=3pt,mark options={solid},line width=1.5pt] table[x=level, y=C_omega] {test10_h-t.dat};
		
		\addplot+[red,mark=o, solid,mark size=2pt,mark options={solid},line width=1.5pt] table[x=level, y=t_L2_omega] {test10_h-t.dat};
		
		\addplot[dashed,line width=1pt] coordinates { 
			(0,\vara) (1,\vara*0.5) (2,\vara*0.25) (3,\vara*0.125) (4,\vara*0.0625) (5,\vara*0.03125) (6,\vara*0.03125/2)
		};
		\addplot[dotted,line width=1pt] coordinates { 
			(0,\varb) (1,\varb*0.5*0.5) (2,\varb*0.25*0.25) (3,\varb*0.125*0.125) (4,\varb*0.0625*0.0625) (5,\varb*0.03125*0.03125) (6,\varb*0.03125*0.03125*0.5*0.5)
		};
		\legend{
			$|| c^* - c_h ||_{\infty,2}$,
			$|| c^* - c_h ||_{2,2}$,
			$|| \mu^* - \mu_h ||_{\infty,2}$,
			$|| \mu^* - \mu_h ||_{2,2}$,
			$\mathcal{O}(h^1)$ ,
			$\mathcal{O}(h^{2})$}
		\end{semilogyaxis}
		\end{tikzpicture}
		\caption{{Discrete $L_\infty(0,T; L^2(\Gamma_h))$ norm \eqref{eq:i2norm} and discrete $L_2(0,T; L^2(\Gamma_h))$ norm \eqref{eq:22norm} of the error for concentration $c$ (blue lines) and chemical potential $\mu$ (red lines) plotted against the refinement level $\ell$ with time step $\Delta t=2^{-(1+\ell)}$.}
		}
		\label{test6_convergence}
	\end{center}
\end{figure}

To assess the numerical solver for the Cahn--Hilliard model, we similarly consider the following exact solution to  non-homogeneous surface Cahn--Hilliard equations with free energy per unit surface \eqref{eq:f0}:
\begin{align*}
c^*={\frac{1}{2}\left(1-0.8e^{-0.4t}\right)\left(Y(x_1,x_2)+ 1\right), \quad Y(x_1,x_2)=x_1 x_2}, \quad t\in[0,5],
\end{align*}
i.e.,  same exact solution we used for order parameter in the Allen--Cahn equation.
The exact chemical potential $\mu^*$ can be readily computed from eq.~\eqref{eq:sys_CH2}.
We set $\rho = 1$, mobility $M$ as in \eqref{eq:M}, and $\epsilon = 0.1$.

We report in Fig.~\ref{test6_convergence} the discrete $L_\infty(0,T; L^2(\Gamma_h))$ norm \eqref{eq:i2norm}
and discrete  $L_2(0,T; L^2(\Gamma_h))$ norm \eqref{eq:22norm} of the error for concentration $c$ (blue lines)
and chemical potential $\mu$ (red lines) plotted against the refinement level $\ell$,
with time step $\Delta t=2^{-(1+\ell)}$.
All the norms reported in Fig.~\ref{test6_convergence} are
computed on the approximate surface $\Gamma_h$, where $c^*$ and $\mu^*$
were defined through their normal extensions from $\Gamma$.
As for the Allen--Cahn equation, we see optimal second order of convergence in the
discrete $L_2(0,T; L^2(\Gamma_h))$ norm for both variables. In the discrete $L_\infty(0,T; L^2(\Gamma_h))$ norm, the concentration converges to the exact solution with the second order, while chemical potential converges with first order. {Note that the definition of the chemical potential involves the derivatives of the concentration and hence the loss of convergence order can be anticipated.}
Both convergence trends are consistent with the numerical analysis  found in the literature for the planar Cahn--Hilliard equations; see, e.g., \cite[Th.3.1]{elliott1989second}, \cite[Th.3.1]{Shen_Yang2010}, \cite[Th.3]{feng2004error}. {
The optimal convergence of concentration for the finite element approximations to surface Cahn--Hilliard was shown
also in \cite{du2011finite,elliott2015evolving}, while for the discrete chemical potential
only boundedness was shown.}

\begin{figure}
  \centering
    \subfloat[$\epsilon = 0.02$, $\ell = 5$]{\includegraphics[width=.3\textwidth]{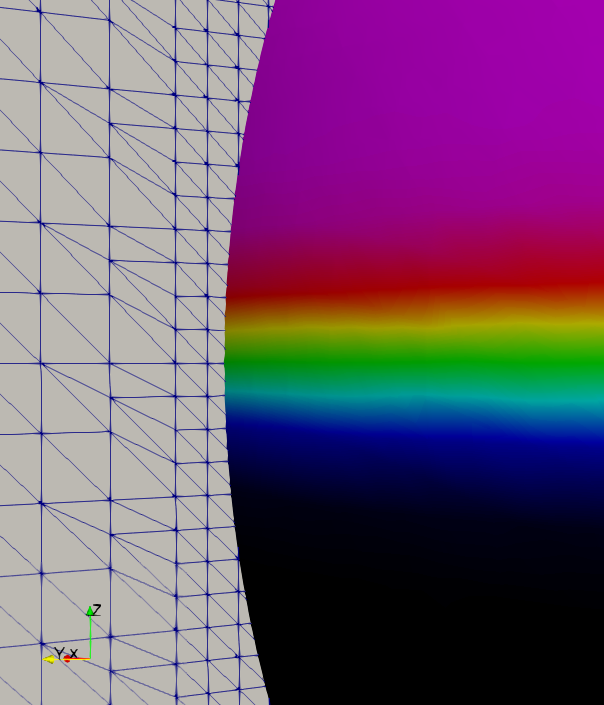}}~
  \subfloat[$\epsilon = 0.02$, $\ell = 6$]{\includegraphics[width=.3\textwidth]{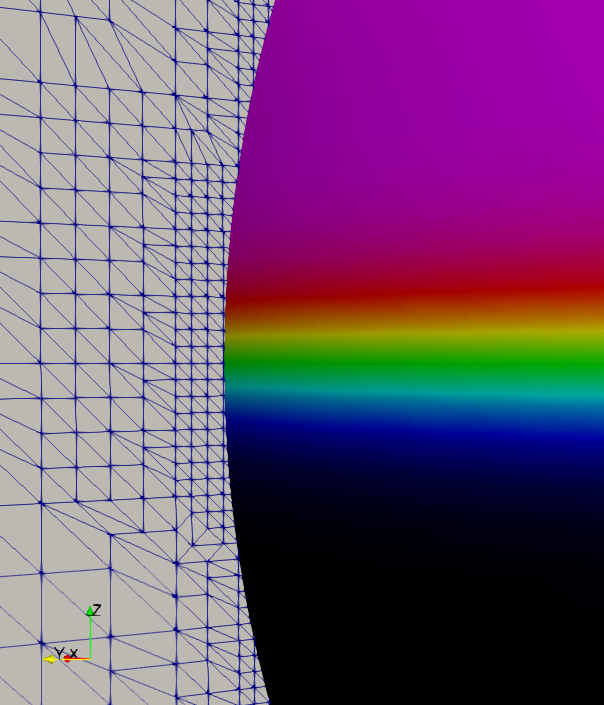}}~
  \subfloat[$\epsilon = 0.01$, $\ell = 6$]{\includegraphics[width=.3\textwidth]{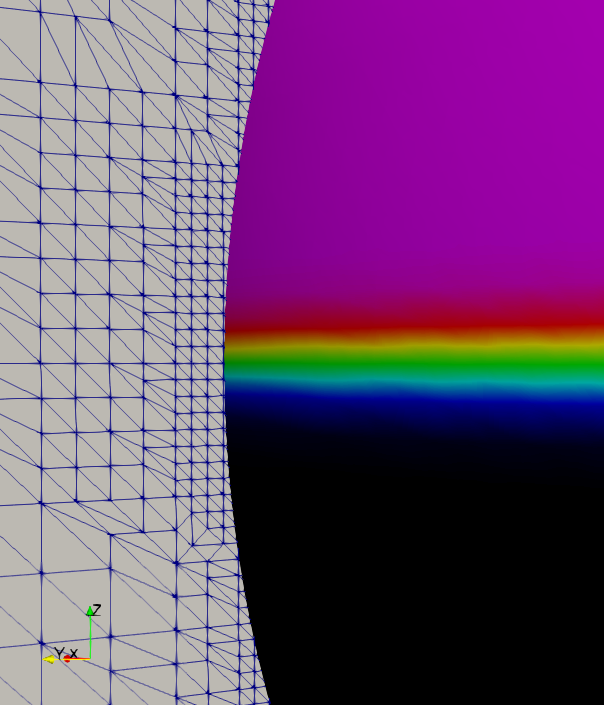}}
  \hfill
  \caption{Magnified view of the interface thickness and bulk mesh near $\Gamma$ for (a) $\epsilon = 0.02$ and mesh $\ell = 5$,
  (b) $\epsilon = 0.02$ and mesh $\ell = 6$, and (c) $\epsilon = 0.01$ and mesh $\ell = 6$.}
  \label{fig:eps_ell}
\end{figure}

In the convergence tests above, we used $\epsilon = 0.1$. The quantity $\epsilon$
is a crucial modeling parameter, since it defines the thickness of a layer where
phase transition takes place. The ability of method to resolve this interfacial phenomena
is critical for the physical fidelity of the simulations. The case of smaller $\epsilon$
is numerically challenging as seen, for example, from stability restriction \eqref{eq:CH_cond} and the blow-up of known error estimates for $\epsilon\to0$ \cite{feng2004error,Shen_Yang2010}.
To ensure that the interfacial region is properly resolved
(i.e.~enough elements fit in its thickness), it is common to apply mesh adaptivity techniques; see, e.g., \cite{yue2006phase,bavnas2008adaptive,nitschke_voigt_wensch_2012}. We will not address
mesh adaptivity in this paper. Therefore, we need to make sure that the finer meshes we consider are suited for thin interfacial regions such as the ones for $\epsilon = 0.02, 0.01$. Notice that these are realistic values of $\epsilon$.
In fact, if we consider a typical giant vesicle with an average diameter of 30 $\mu$m,
on which phase separation can be visualized using fluorescence microscopy
\cite{VEATCH20033074} with a resolution of about 300 nm, the thickness of transition region between
the phases is approximately 1\% of vesicle diameter.

For this purpose, we consider the eq.~\eqref{eq:AC_eq} with free energy per unit surface \eqref{eq:f0}
and initial condition is $\eta_0(\bx)=x_3+1/2$. We let the simulation run until changes in the interface thickness
cannot be visually appreciated. Fig.~\ref{fig:eps_ell} show a magnified view of the interface thickness
at $t=30$ for $\epsilon = 0.02$ and meshes $\ell = 5,6$, and at $t=60$ for for $\epsilon = 0.01$ and mesh $\ell = 6$,
together with the bulk mesh near the surface. From Fig.~\ref{fig:eps_ell} (a) and (b), we see that the thickness of the computed  interface does not vary significantly from mesh $\ell = 5$ to mesh
$\ell = 6$, indicating that both meshes are sufficiently refined for $\epsilon = 0.02$.
Fig.~\ref{fig:eps_ell} (b) and (c) show that the thickness of the interface computed with mesh
$\ell = 6$ gets halved when $\epsilon$ goes from 0.02 to 0.01.
Thus, we use  mesh $\ell = 6$ to perform all numerical simulation reported further.

\subsection{Phase separation on a sphere} \label{sec:sphere}
The surface of the sphere is appealing not only for its simplicity, which makes it ideal
for benchmarking, but also for its relevance in practical applications. In fact,
lipid vesicles used as drug carriers have a spherical shape \cite{Laouini_et_al2012}.
Therefore, we first focus on the process of phase separation on the surface of the sphere.

We consider both Allen--Cahn and Cahn--Hilliard equations,
i.e. eq.~\eqref{eq:AC_eq} and \eqref{eq:CH}, with free energy per unit surface \eqref{eq:f0} and $\epsilon = 0.01$. The initial condition is still $\eta_0(\bx)=\text{rand}(\bx)$. Notice that this initial state corresponds to having a single thermodynamic phase, i.e. a ``homogeneous'' mix
of the components. This miscible chemical mixture is unstable and proceed to separate
into two distinct phases by diffusion. This process is characterized by two time scales:
an initial fast stage followed by a slower process.
In fact, the minimization of the chemical energy results in very fast development
in the early stage of the phase separation.
Later on, in the coarsening process, the dissipation of the interfacial energy
is orders of magnitude slower.

\subsubsection{Allen--Cahn model}
\begin{figure}
\begin{center}
\begin{overpic}[width=.20\textwidth,viewport=170 20 630 470, clip,grid=false]{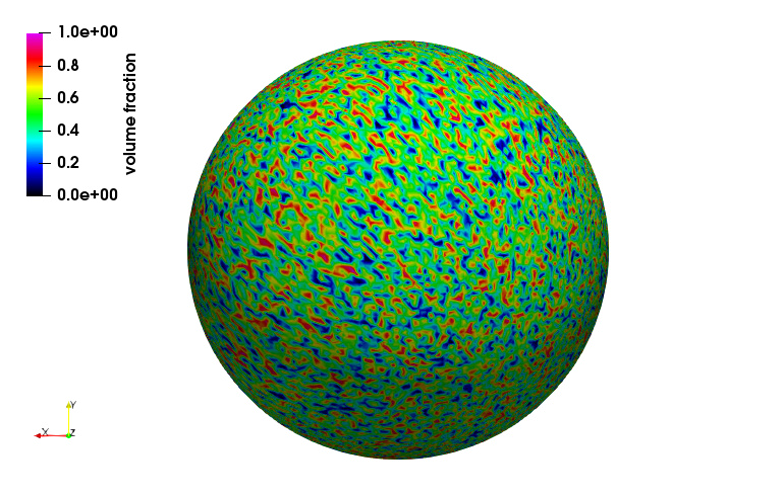}
        \put(38,100){\small{$t = 0$}}
\end{overpic}
\begin{overpic}[width=.20\textwidth, viewport=170 20 630 470, clip,grid=false]{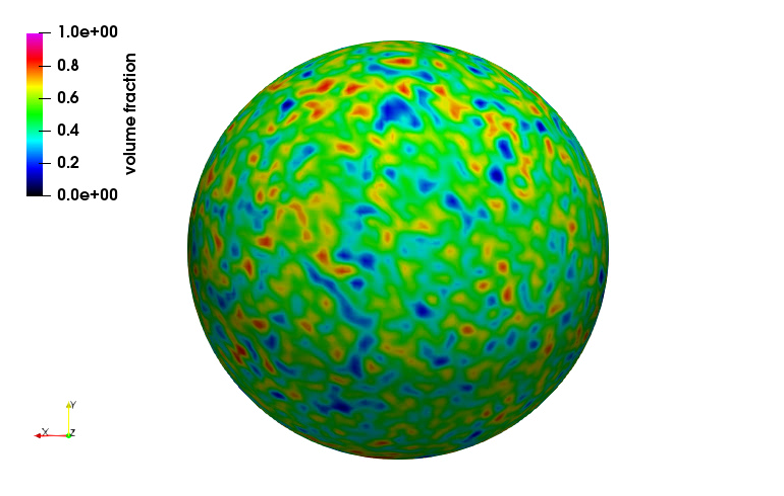}
         \put(38,100){\small{$t = 5$}}
\end{overpic}
\begin{overpic}[width=.20\textwidth, viewport=170 20 630 470, clip,grid=false]{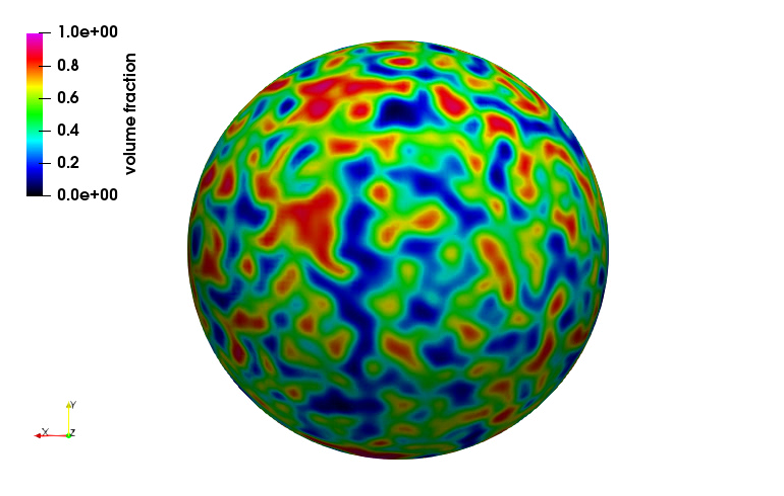}
         \put(36,100){\small{$t = 10$}}
\end{overpic}
\begin{overpic}[width=.20\textwidth, viewport=170 20 630 470, clip,grid=false]{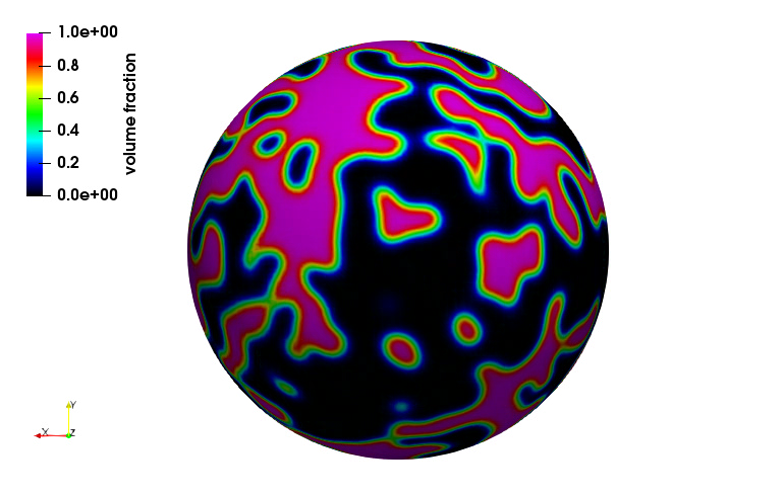}
         \put(32,100){\small{$t = 25$}}
\end{overpic}
\vskip .3cm
\begin{overpic}[width=.20\textwidth, viewport=170 20 630 470, clip,grid=false]{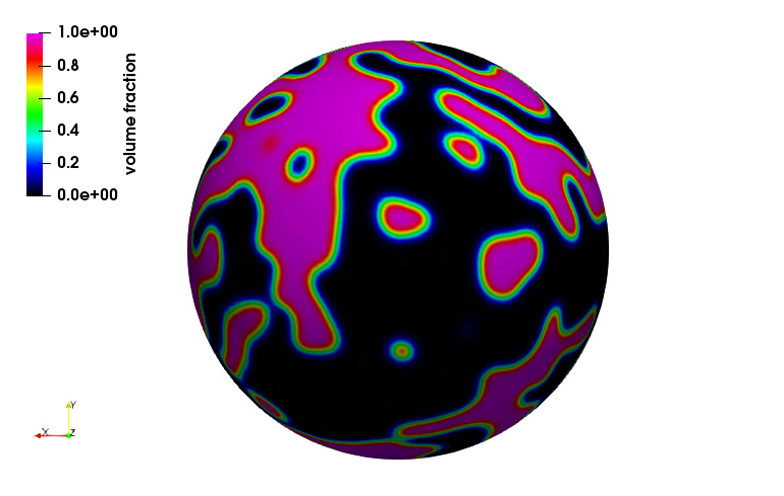}
         \put(38,100){\small{$t = 50$}}
\end{overpic}
\begin{overpic}[width=.20\textwidth, viewport=170 20 630 470, clip,grid=false]{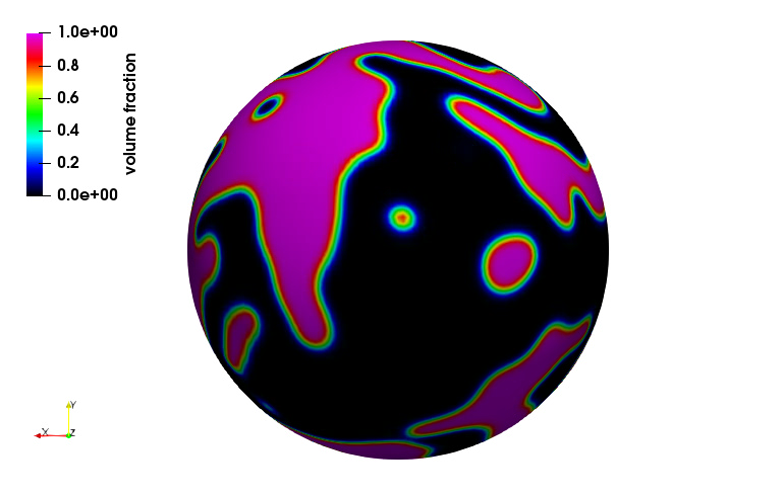}
         \put(35,100){\small{$t = 75$}}
\end{overpic}
\begin{overpic}[width=.20\textwidth, viewport=170 20 630 470, clip,grid=false]{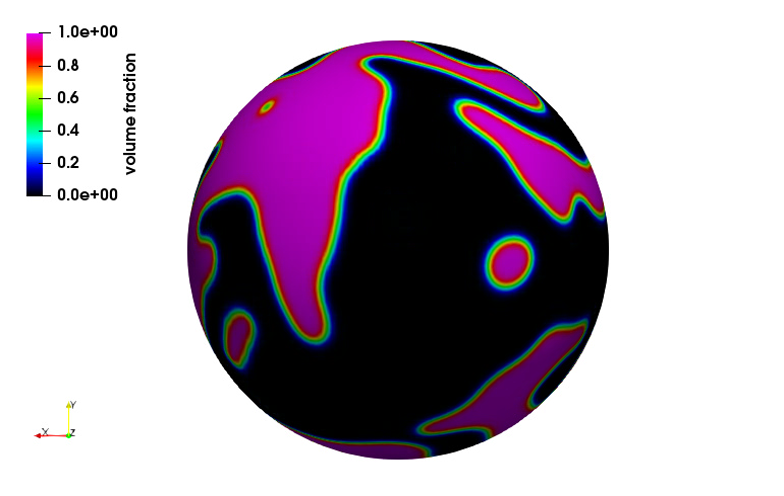}
         \put(32,100){\small{$t = 100$}}
\end{overpic}
\begin{overpic}[width=.20\textwidth, viewport=170 20 630 470, clip,grid=false]{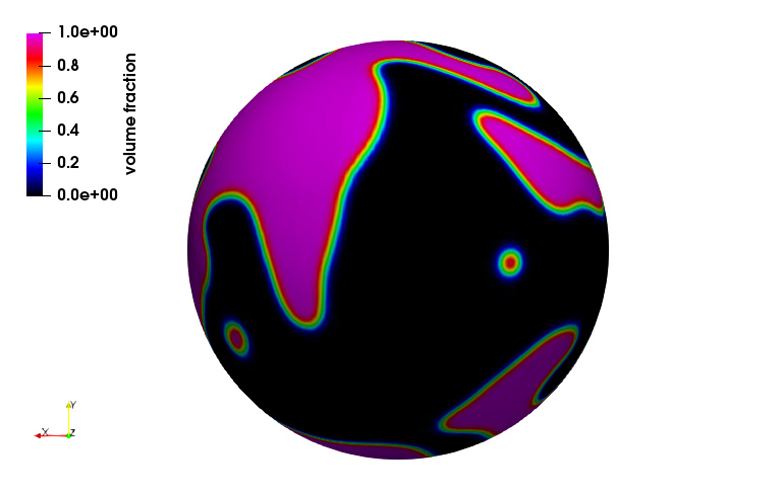}
         \put(32,100){\small{$t = 150$}}
\end{overpic}
\vskip .3cm
\begin{overpic}[width=.20\textwidth, viewport=170 20 630 470, clip,grid=false]{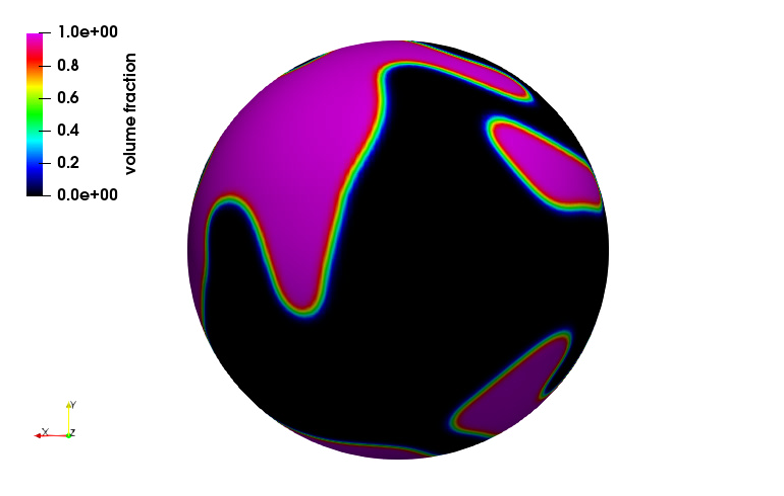}
         \put(32,100){\small{$t = 200$}}
\end{overpic}
\begin{overpic}[width=.20\textwidth, viewport=170 20 630 470, clip,grid=false]{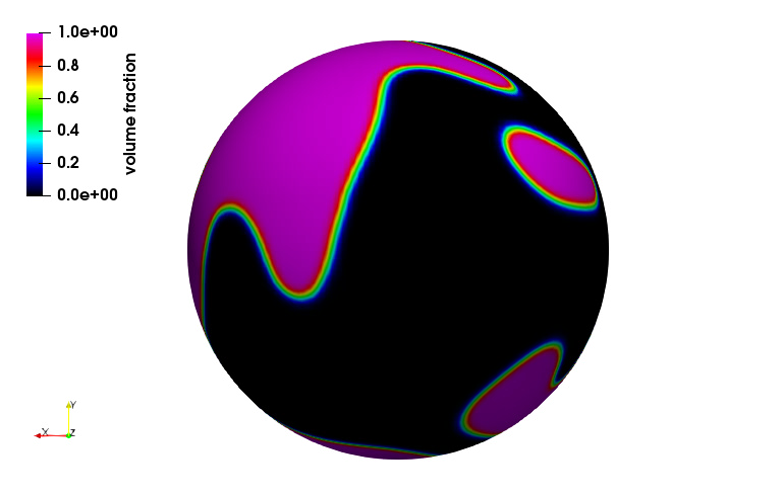}
         \put(32,100){\small{$t = 300$}}
\end{overpic}
\begin{overpic}[width=.20\textwidth, viewport=170 20 630 470, clip,grid=false]{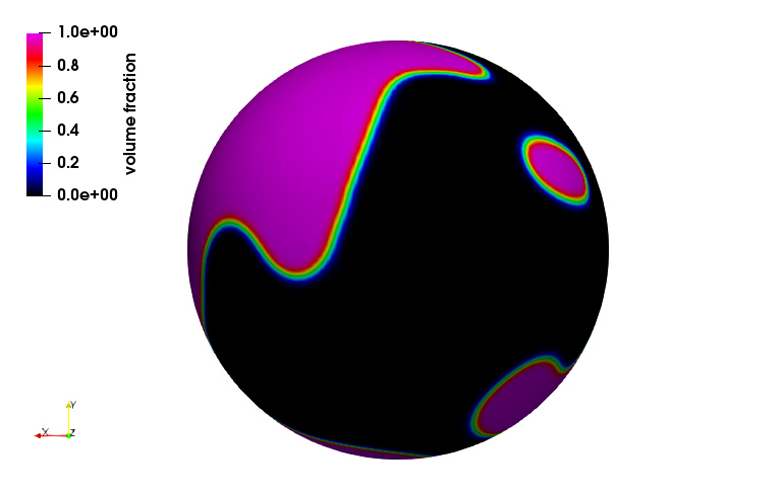}
         \put(32,100){\small{$t = 400$}}
\end{overpic}\hspace{1.2cm}
\begin{overpic}[width=.10\textwidth,grid=false]{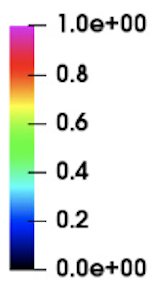}
\end{overpic}
\end{center}
\caption{\label{sphere_0_4}
Evolution of the numerical solution of the Allen--Cahn equation on the sphere for $t \in (0, 400]$
computed with mesh $\ell = 6$ and time step as specified in Table \ref{tab:dt_sphere}.
The simulation was started with a random initial condition depicted in the top left panel. View: $xy$-plane.}
\end{figure}

\begin{figure}
\begin{center}
\begin{overpic}[width=.20\textwidth,viewport=120 20 830 770, clip,grid=false]{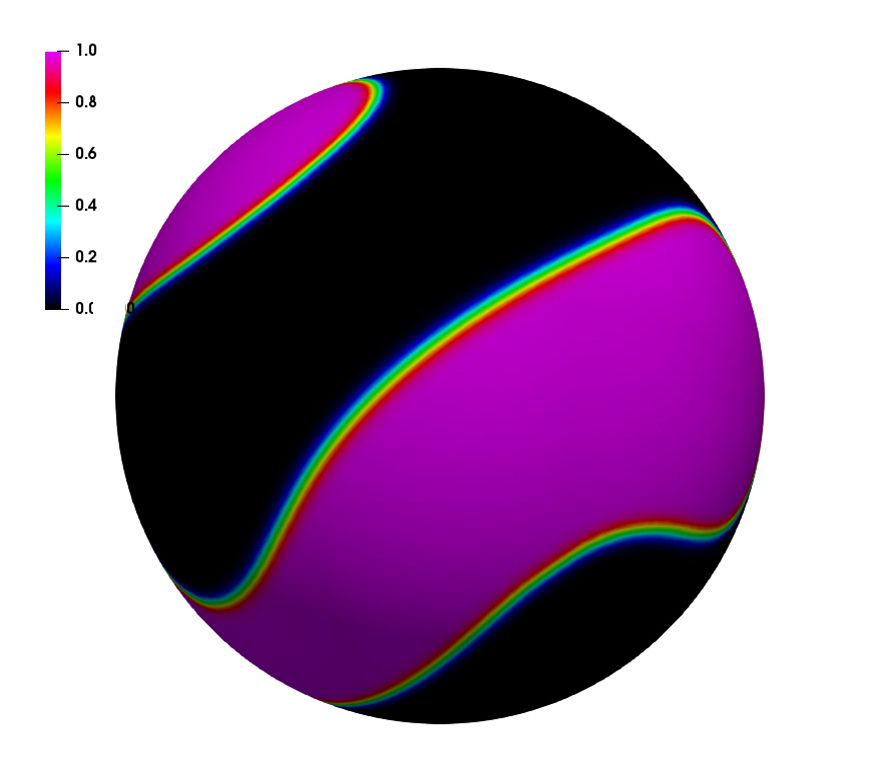}
        \put(21,100){\small{$t = 1060$}}
\end{overpic}
\begin{overpic}[width=.20\textwidth, viewport=120 20 830 770, clip,grid=false]{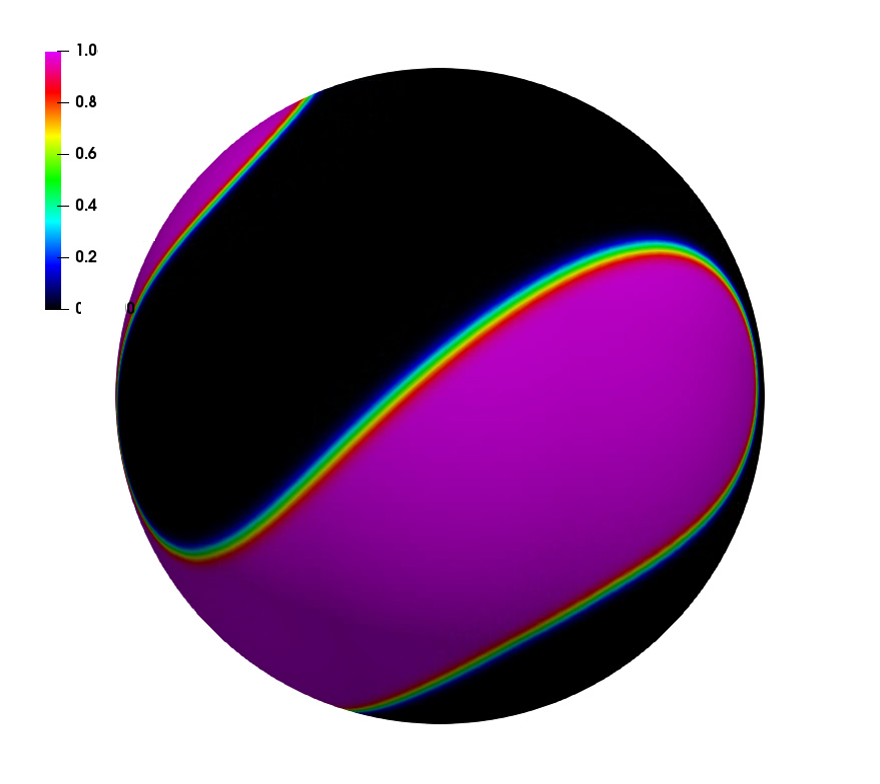}
         \put(21,100){\small{$t = 2060$}}
\end{overpic}
\begin{overpic}[width=.20\textwidth, viewport=120 20 830 770, clip,grid=false]{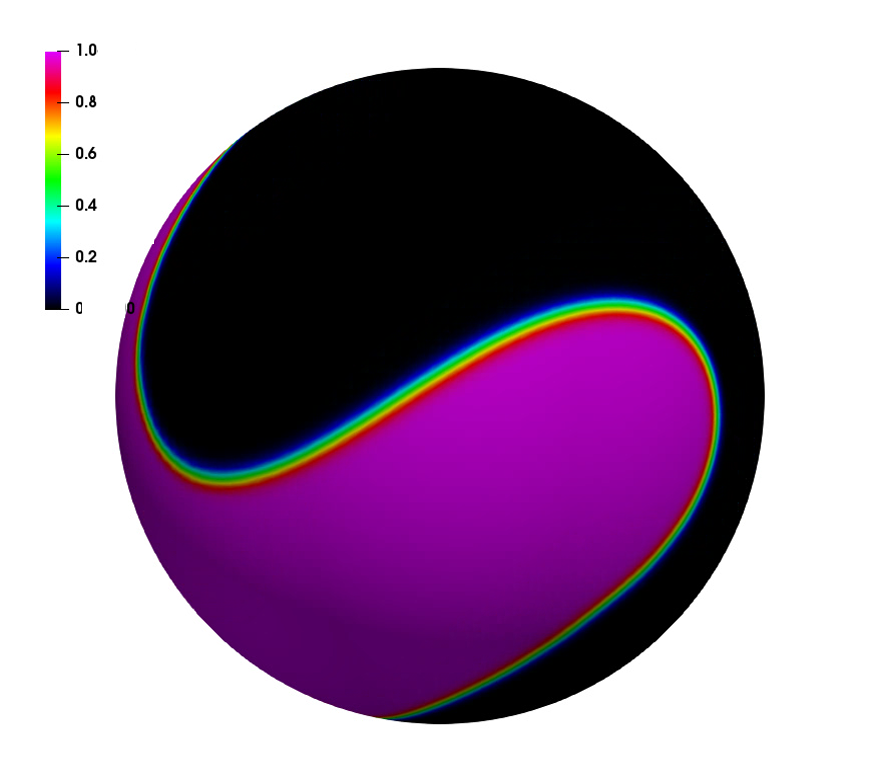}
         \put(21,100){\small{$t = 3060$}}
\end{overpic}
\begin{overpic}[width=.20\textwidth, viewport=120 20 830 770, clip,grid=false]{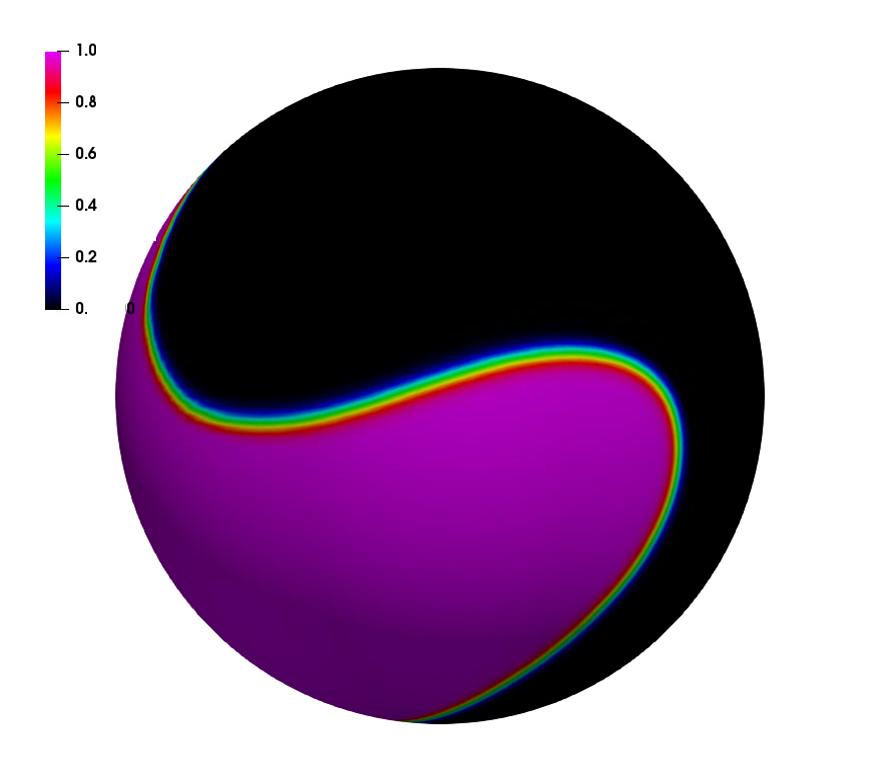}
         \put(21,100){\small{$t = 4060$}}
\end{overpic}
\vskip .3cm
\begin{overpic}[width=.20\textwidth, viewport=120 20 830 770, clip,grid=false]{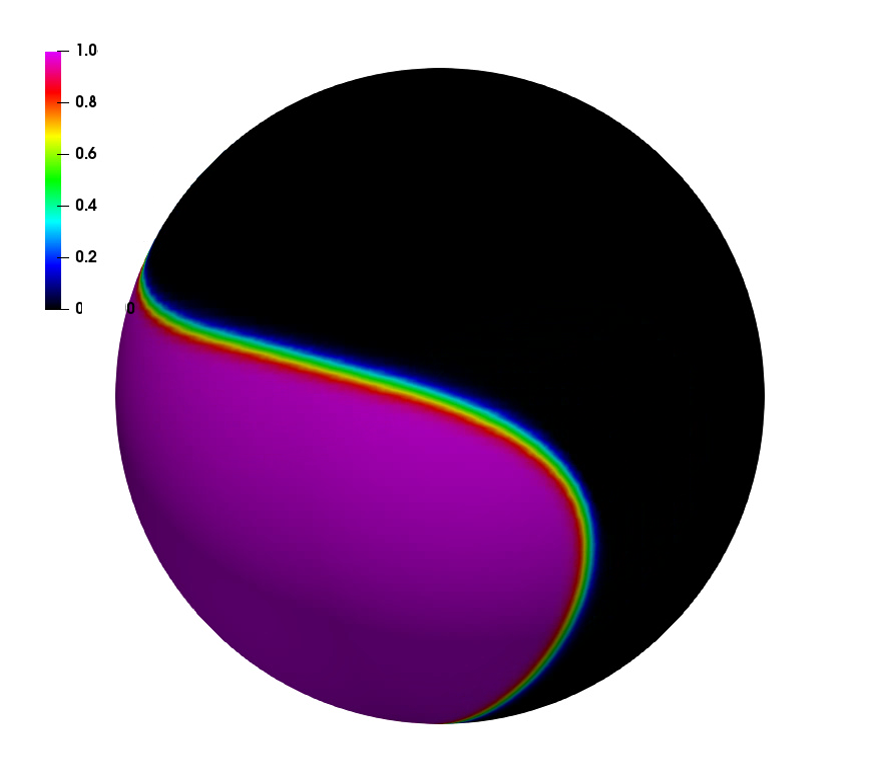}
         \put(21,100){\small{$t = 6060$}}
\end{overpic}
\begin{overpic}[width=.20\textwidth, viewport=120 20 830 770, clip,grid=false]{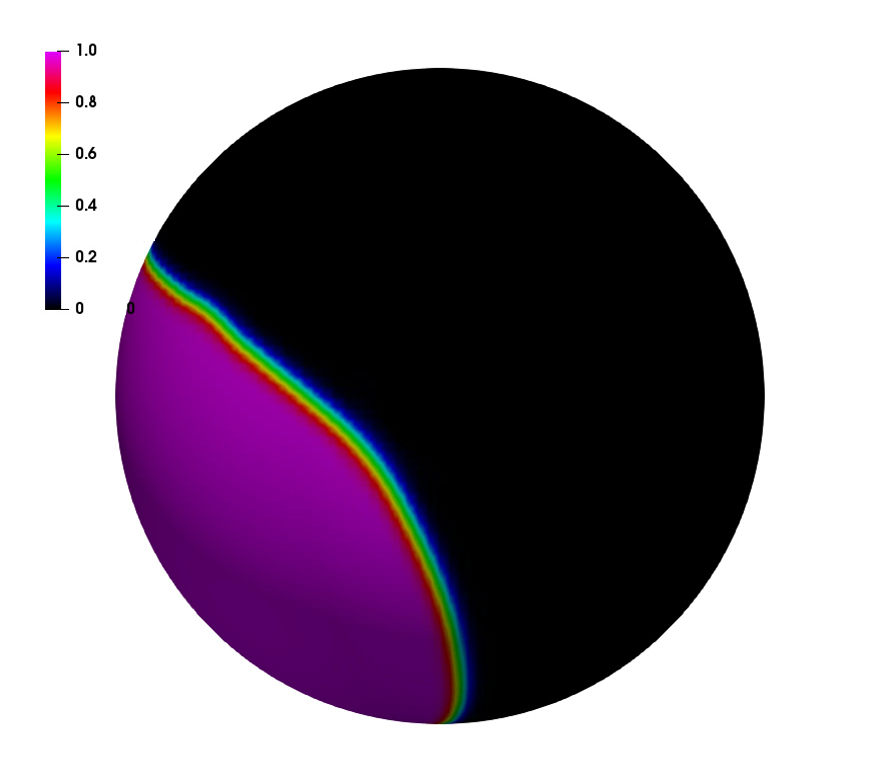}
         \put(21,100){\small{$t = 9060$}}
\end{overpic}
\begin{overpic}[width=.20\textwidth, viewport=120 20 830 770, clip,grid=false]{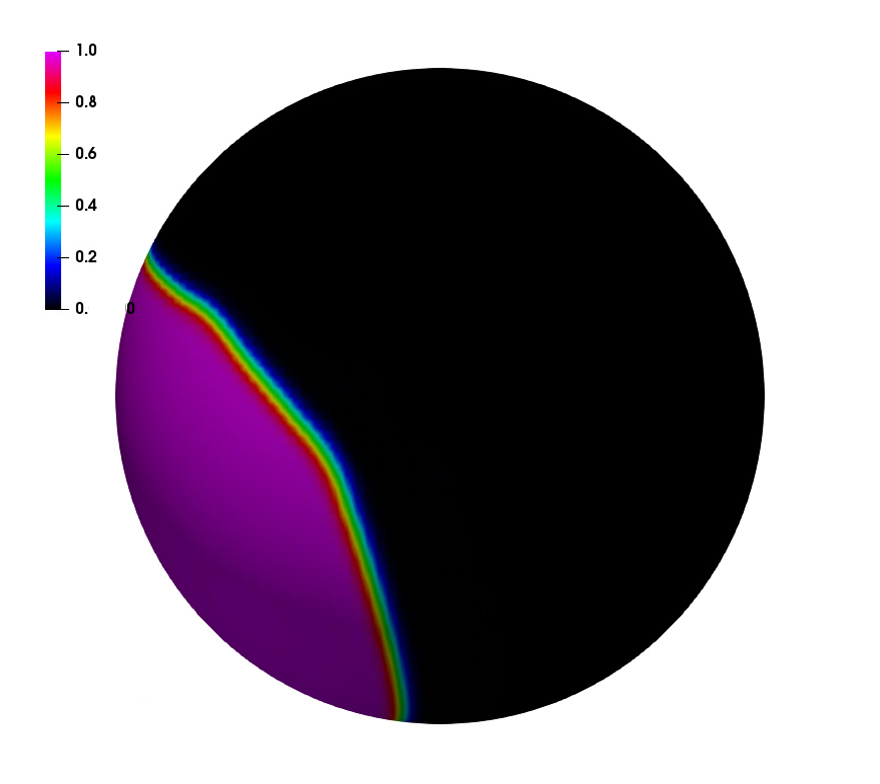}
         \put(19,100){\small{$t = 12060$}}
\end{overpic}
\hspace{1cm}
\begin{overpic}[width=.10\textwidth,grid=false]{figures_AC_sphere_legend.png}
\end{overpic}
\end{center}
\caption{\label{sphere_10_120}
Evolution of the numerical solution of the Allen--Cahn equation  on the sphere for $t \in [1060, 12060]$,
i.e.~after the separation into two regions has occurred,
computed with mesh $\ell = 6$  and time step as specified in Table \ref{tab:dt_sphere}.
View: $yz$-plane.}
\end{figure}

\begin{figure}
\begin{center}
\begin{overpic}[width=.20\textwidth,viewport=110 20 830 770, clip,grid=false]{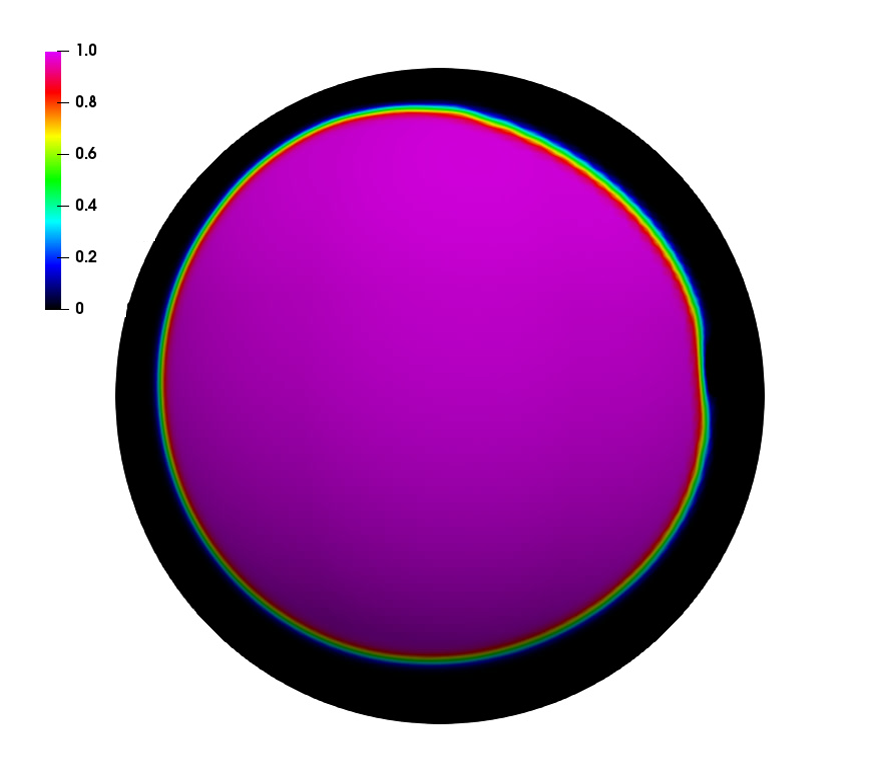}
        \put(18,100){\small{$t = 13560$}}
\end{overpic}
\begin{overpic}[width=.20\textwidth, viewport=110 20 830 770, clip,grid=false]{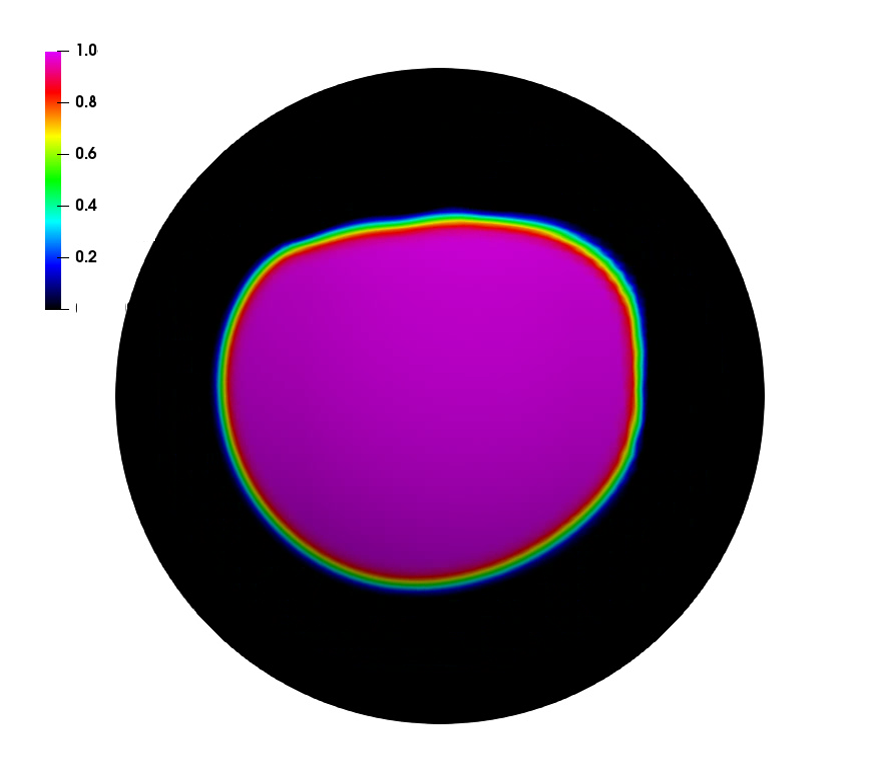}
         \put(18,100){\small{$t = 17560$}}
\end{overpic}
\begin{overpic}[width=.20\textwidth, viewport=110 20 830 770, clip,grid=false]{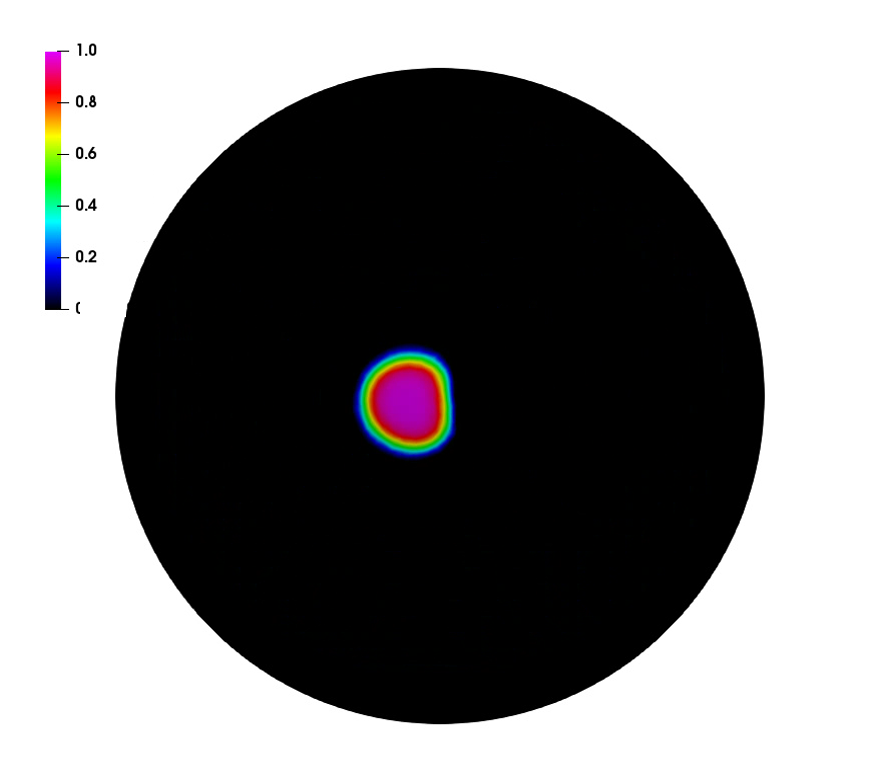}
         \put(24,100){\small{$t = 21560$}}
\end{overpic}
\begin{overpic}[width=.20\textwidth, viewport=110 20 830 770, clip,grid=false]{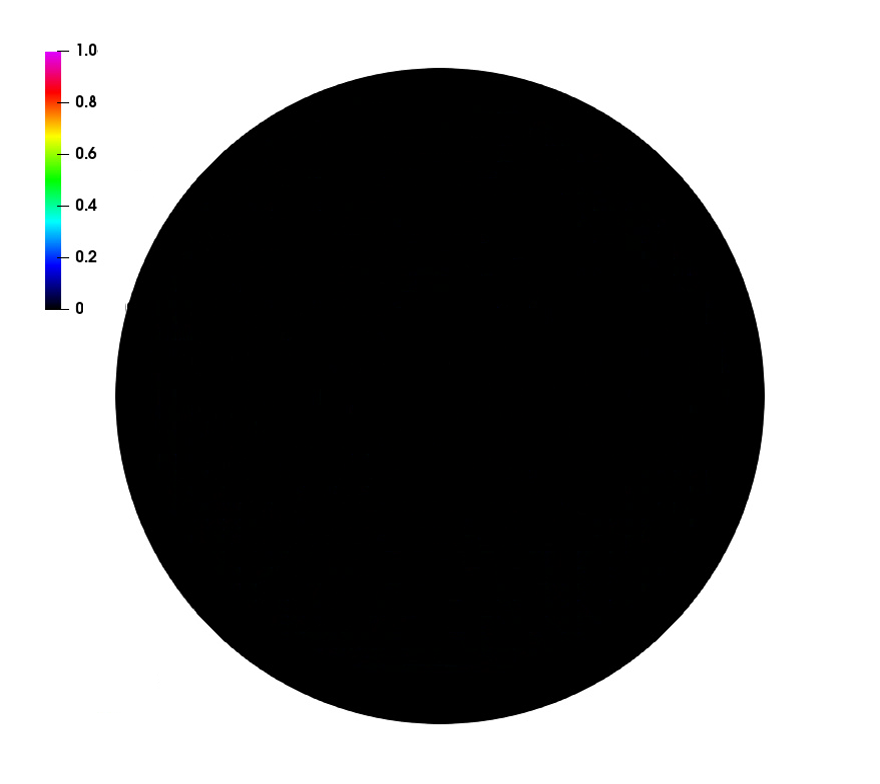}
         \put(24,100){\small{$t = 22560$}}
\end{overpic}
\hspace{.3cm}
\begin{overpic}[width=.10\textwidth,grid=false]{figures_AC_sphere_legend.png}
\end{overpic}
\end{center}
\caption{\label{sphere_135_225}
Evolution of the numerical solution of the Allen--Cahn equation on the sphere for $t \in [13560, 22560]$
computed with mesh $\ell = 6$ and time step as specified in Table \ref{tab:dt_sphere}.
Angled view of the $xz$-plane.}
\end{figure}

The evolution of the numerical solution to the Allen--Cahn equation
for $t \in (0, 400]$ is shown in Fig.~\ref{sphere_0_4}.
We see that after the initial fast stage, which ends in a little more than
ten time units, there is a considerable slow down in the evolution of the solution.
The separation into two regions, one pink region with $\eta = 1$ and one black region
with $\eta = 0$, happens around $t=800$.
Fig.~\ref{sphere_10_120} displays the evolution of the numerical solution of the Allen--Cahn equation
for $t \in [1060, 12060]$, i.e. after the separation into two regions has occurred. We notice
a further deceleration in the process of dissipation of the interfacial energy. Finally, Fig.~\ref{sphere_135_225}
shows how the solution evolves until one phase (the pink one) disappears at around $t = 22560$.
We recall that the Allen--Cahn equation describes the evolution of a non-conserved order parameter
during phase transformation. Thus it is expected that one phase vanishes on the sphere, since
it is not possible to trace a curve of minimal length on its surface, i.e.~there is no minimal
length interface (cf. Sec.~\ref{sec:scell}).
As we will see next, this is not the case when phase separation on the sphere
is modeled by the Cahn--Hilliard equation.

Fig.~\ref{sphere_0_4}, \ref{sphere_10_120}, and \ref{sphere_135_225} clearly show
the different time scales involved in the phase separation process modeled by the
Allen--Cahn equation. In order to follow it, we prescribe different time steps for the
different stages. Table \ref{tab:dt_sphere} reports the time step assigned to each time interval.
It would be less intrusive to use some time-adaptivity strategy \cite{guillen2014second,shen_et_al2016},
which will be addressed in a forthcoming paper.

\begin{table}[h!]
\begin{center}
  \begin{tabular}{ | c | c | c |  c |  c | c |  c |}
    \hline
    interval & $(0, 10]$ & $(10, 60]$ & $(60, 1060]$ & $(1060, 3560]$ & $(3560, 13560]$ & $(13560, 22560]$  \\
    \hline
    $\Delta t$ & 1 & 5 & 10 & 50 & 100 & 200 \\
    \hline
  \end{tabular}
  \caption{\label{tab:dt_sphere}
Time steps used for the different time intervals to obtain the results in Fig.~\ref{sphere_0_4}, \ref{sphere_10_120}, and \ref{sphere_135_225}.}
\end{center}
\end{table}

\subsubsection{Cahn--Hilliard model}\label{sec:sphere_CH}

\begin{figure}
\begin{center}
\begin{overpic}[width=.20\textwidth,viewport=160 140 660 670, clip,grid=false]{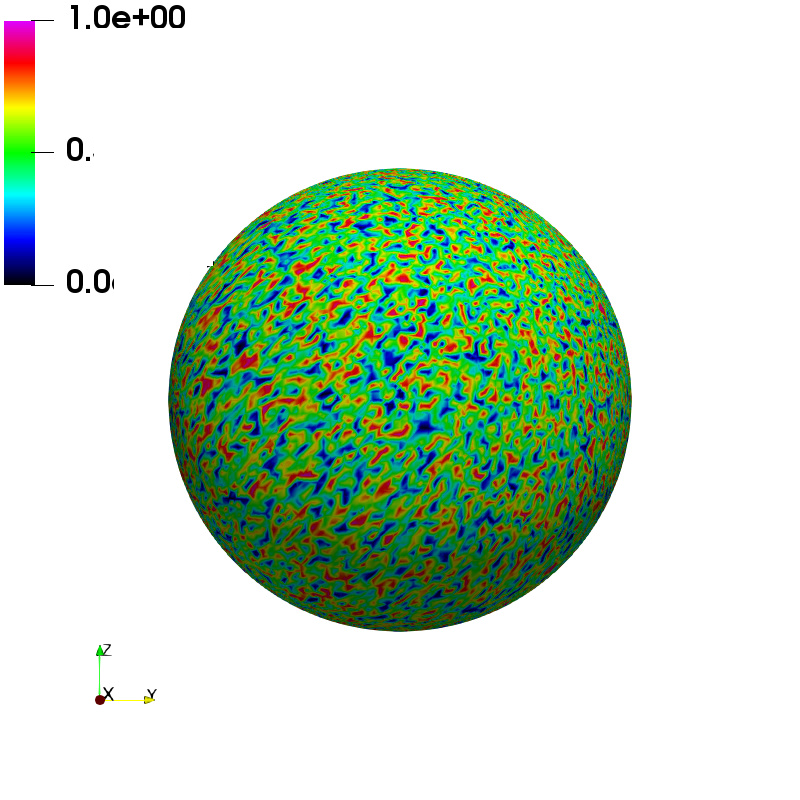}
        \put(32,100){\small{$t = 0$}}
\end{overpic}
\begin{overpic}[width=.20\textwidth, viewport=160 140 660 670, clip,grid=false]{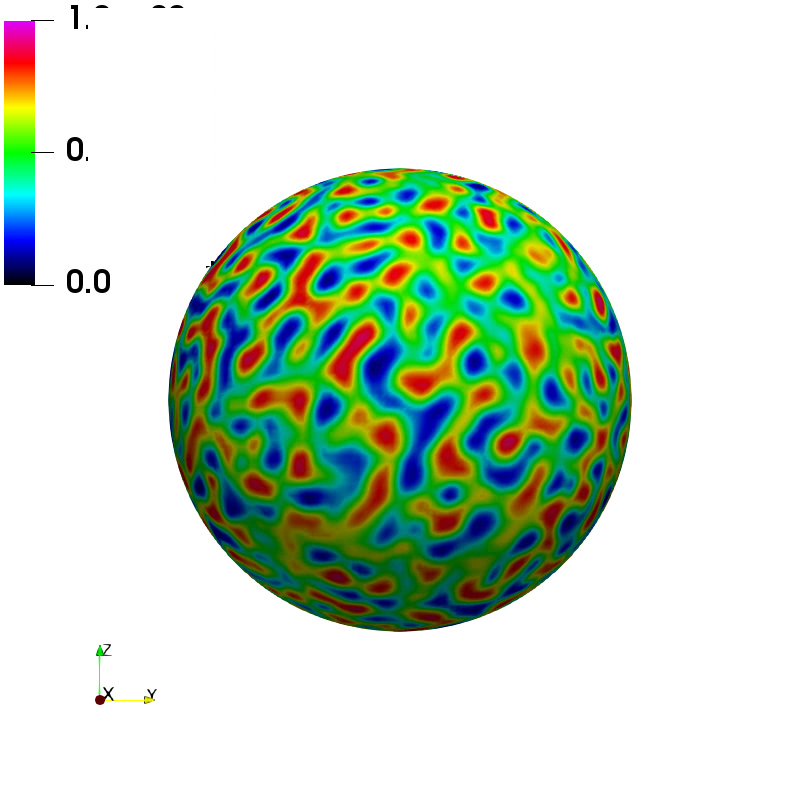}
         \put(28,100){\small{$t = 0.1$}}
\end{overpic}
\begin{overpic}[width=.20\textwidth, viewport=160 140 660 670, clip,grid=false]{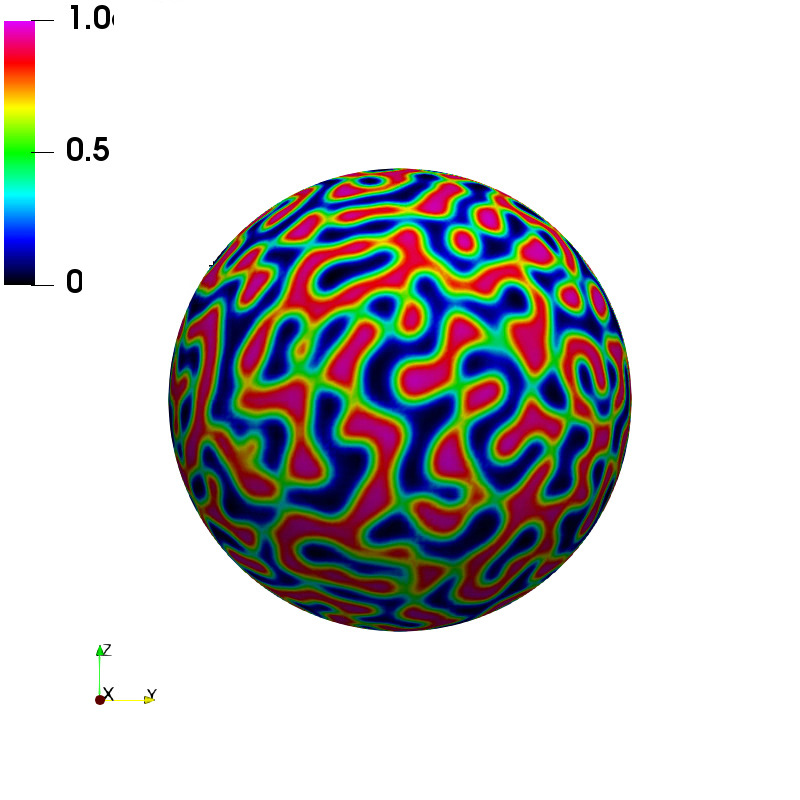}
         \put(28,100){\small{$t = 0.2$}}
\end{overpic}
\begin{overpic}[width=.20\textwidth, viewport=160 140 660 670, clip,grid=false]{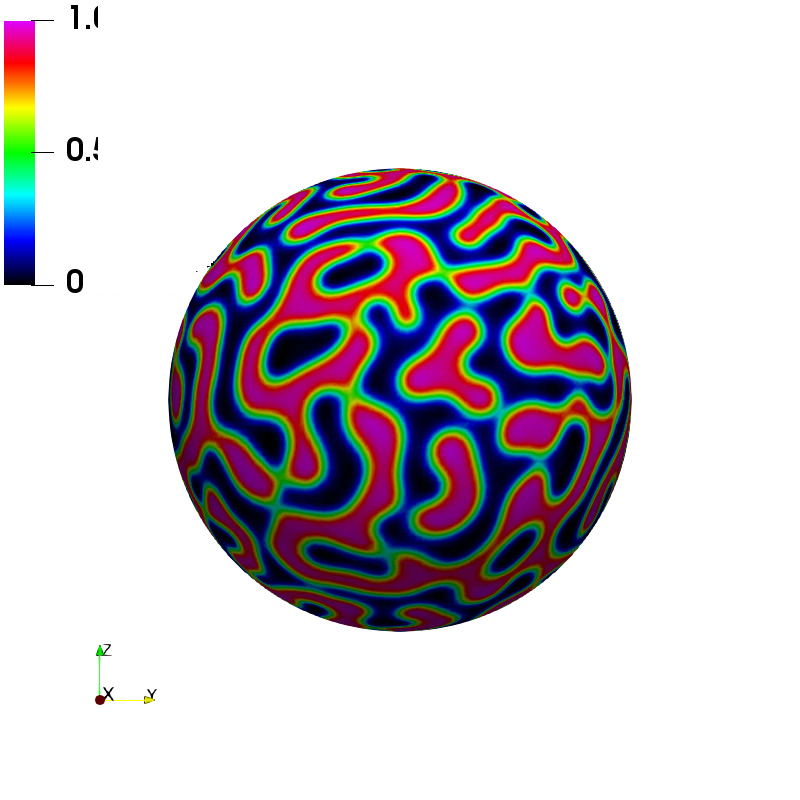}
         \put(28,100){\small{$t = 0.5$}}
\end{overpic}
\vskip .3cm
\begin{overpic}[width=.20\textwidth,viewport=160 140 660 670, clip,grid=false]{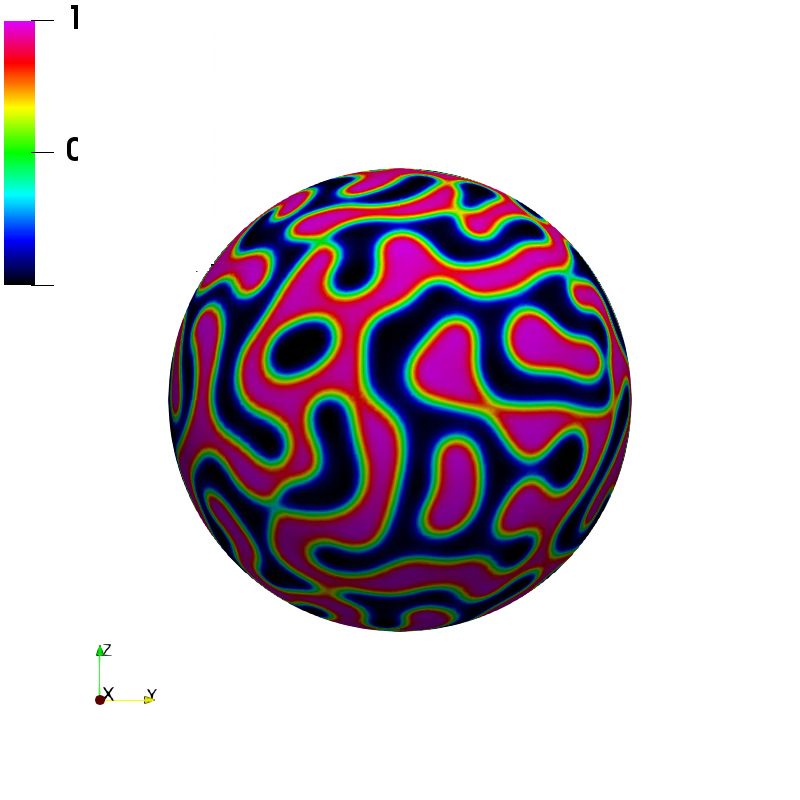}
        \put(30,100){\small{$t = 1$}}
\end{overpic}
\begin{overpic}[width=.20\textwidth,viewport=160 140 660 670, clip,grid=false]{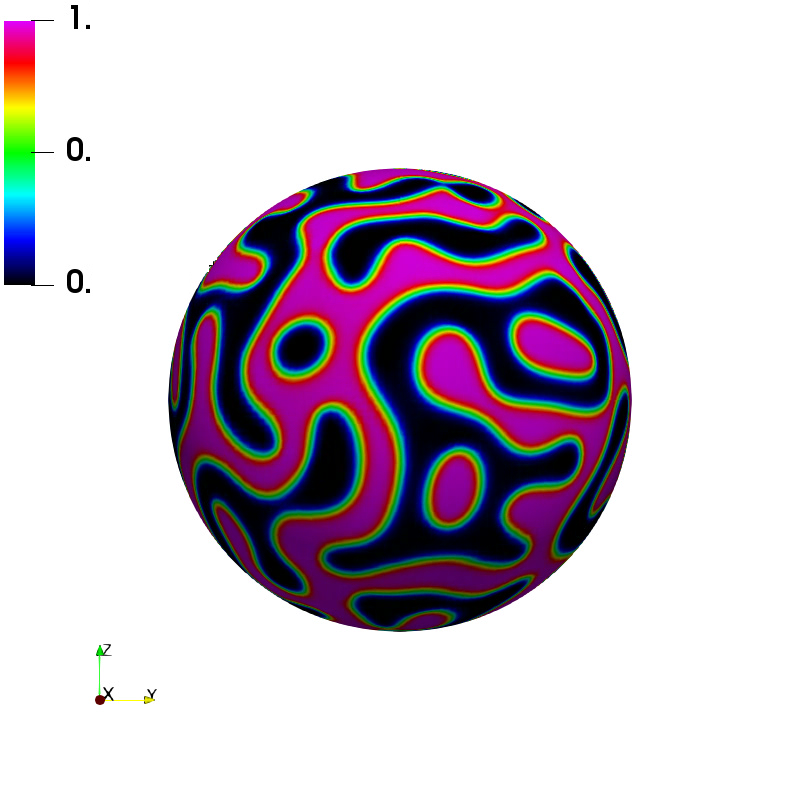}
        \put(30,100){\small{$t = 10$}}
\end{overpic}
\begin{overpic}[width=.20\textwidth,viewport=160 140 660 670, clip,grid=false]{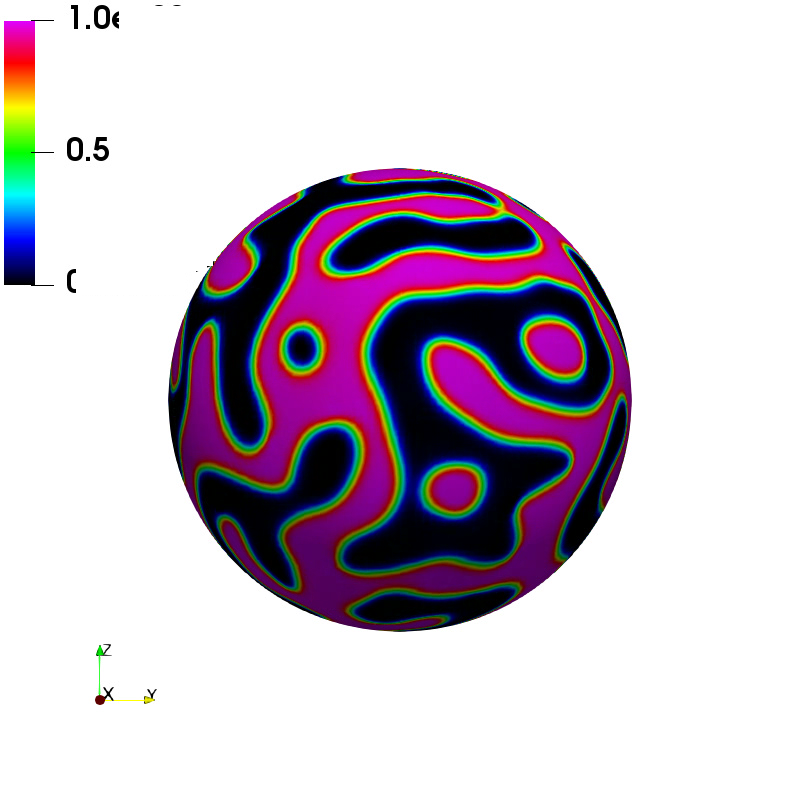}
        \put(30,100){\small{$t = 20$}}
\end{overpic}
\begin{overpic}[width=.20\textwidth,viewport=160 140 660 670, clip,grid=false]{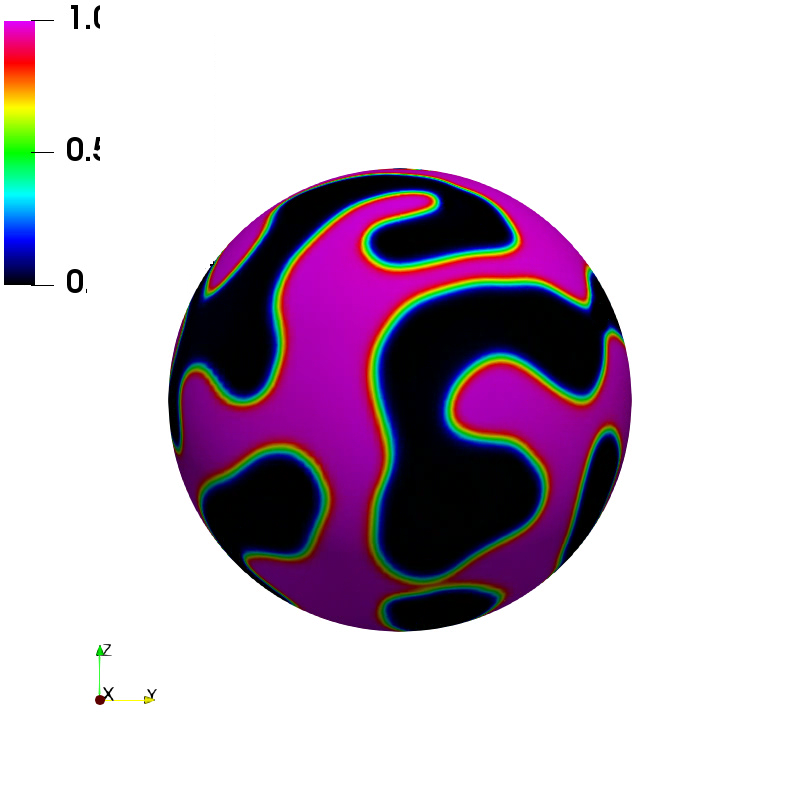}
        \put(30,100){\small{$t = 50$}}
\end{overpic}
\vskip .3cm
\begin{overpic}[width=.20\textwidth,viewport=160 140 660 670, clip,grid=false]{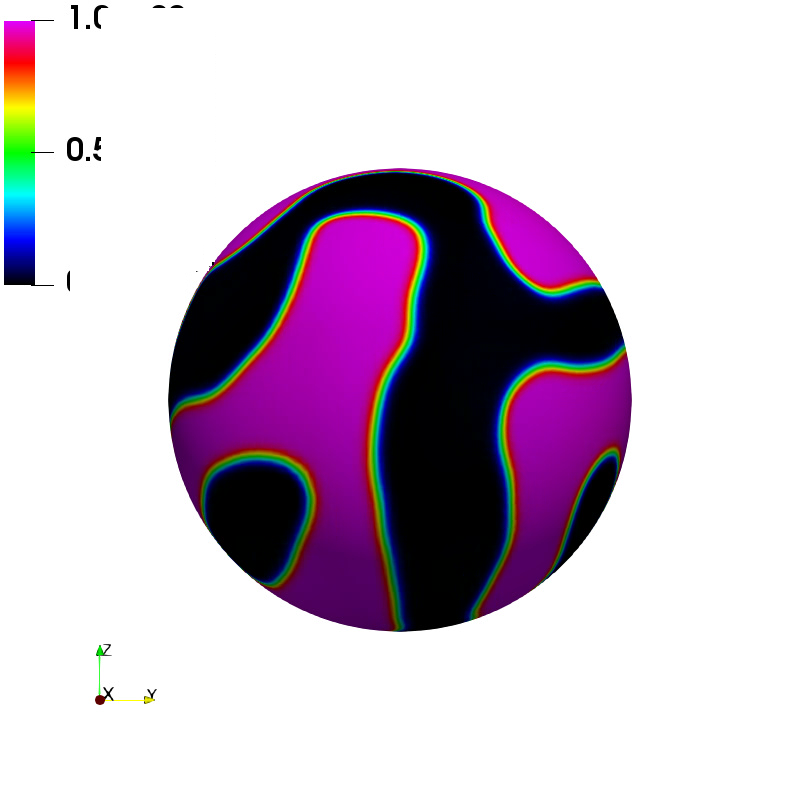}
        \put(30,100){\small{$t = 100$}}
\end{overpic}
\begin{overpic}[width=.20\textwidth,viewport=160 140 660 670, clip,grid=false]{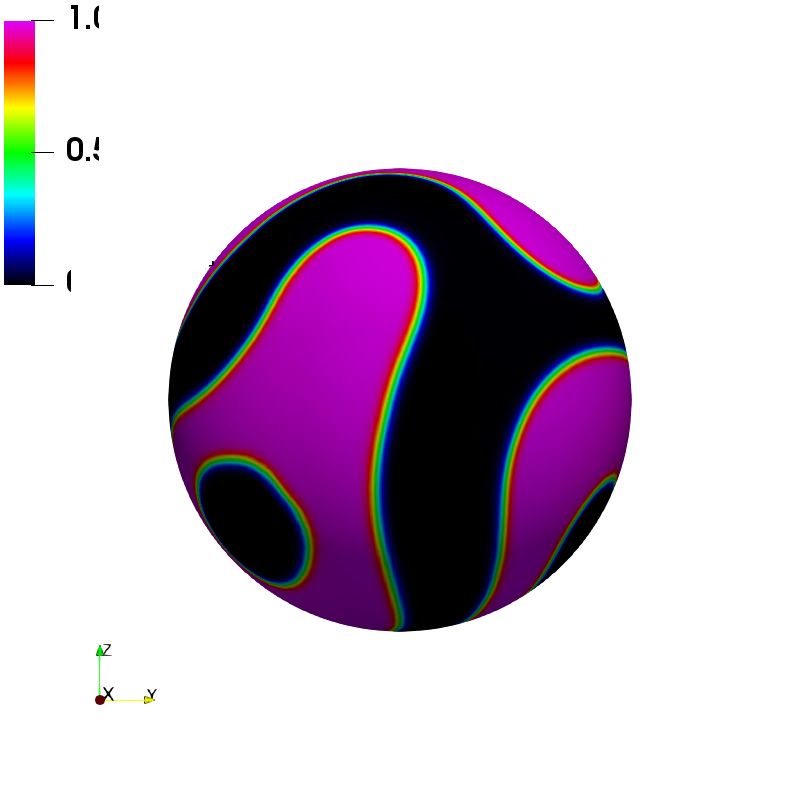}
        \put(30,100){\small{$t = 200$}}
\end{overpic}
\begin{overpic}[width=.20\textwidth,viewport=160 140 660 670, clip,grid=false]{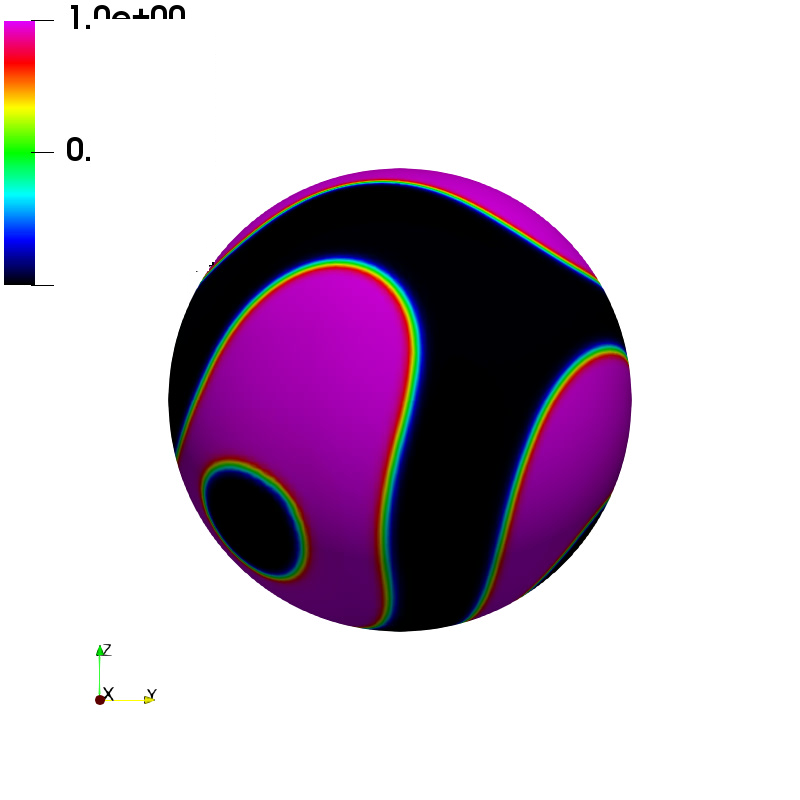}
        \put(30,100){\small{$t = 500$}}
\end{overpic}
\hspace{1cm}
\begin{overpic}[width=.10\textwidth,grid=false]{figures_AC_sphere_legend.png}
\end{overpic}
\end{center}
\caption{\label{sphereCH_0_5}
Evolution of the numerical solution of the Cahn--Hilliard equation on the sphere for $t \in (0, 500]$
computed with mesh $\ell = 6$ .
Time step is $\Delta t =0.01$ for $t \in (0, 1]$ and $\Delta t =1$
for $t \in (1, 500]$.
View: $yz$-plane.}
\end{figure}

Fig.~\ref{sphereCH_0_5} shows
the evolution of the numerical solution to the Cahn--Hilliard problem
for $t \in (0, 500]$. We see that during the initial stage a pattern forms
much faster than in the evolution modeled by the Allen--Cahn equation. In fact
from Fig.~\ref{sphereCH_0_5} we observe the emergence of a pattern already at
time $t = 0.5$, while Fig.~\ref{sphere_0_4} shows no pattern yet at $t = 5$.
We recall that we have used the same free energy $f_0$ for both the Allen--Cahn
and Cahn--Hilliard models. 
Because of this difference in the early stage, Fig.~\ref{sphere_0_4}
and \ref{sphereCH_0_5} display solutions computed at different times.
Also in the evolution given by the Cahn--Hilliard equation,
we observe a considerable slow down after the initial fast stage, which ends around $t = 1$.
After $t = 50$, the process of dissipation of the interfacial energy seems to slow down even further.

Fig.~3 (b) in \cite{VEATCH20033074} shows the spinodal decomposition observed experimentally in a vesicle
that has a 1:1 concentration of DOPC (an unsaturated lipid)/DPPC (a saturated phospholipid)
with 35\% cholesterol. This ternary mixture at a certain temperature gives rise to separation into
two phases: a dark, liquid-ordered phase that is rich in DPPC and
cholesterol, and a bright, less ordered phase that is rich in DOPC.
The panels in Fig.~\ref{sphereCH_0_5} associated to $t = 0.5, 1, 10$ resemble
the images in Fig.~3 (b) in \cite{VEATCH20033074}, indicating that the Cahn--Hilliard equation
provides a {possible} model for spinodal decomposition in vesicles.

\begin{figure}
\begin{center}
\begin{overpic}[width=.18\textwidth,viewport=160 140 660 670, clip,grid=false]{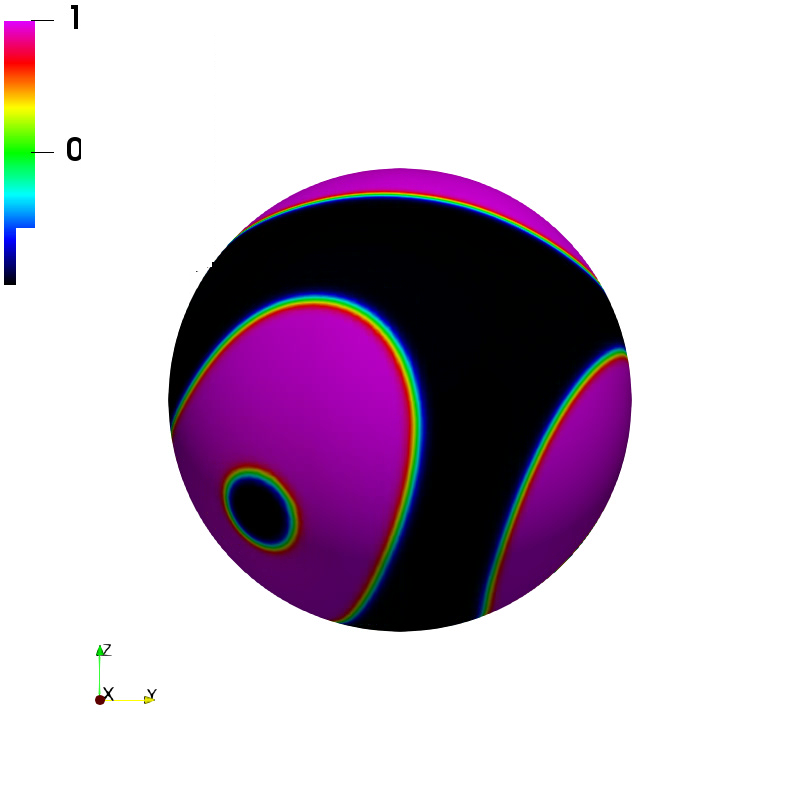}
        \put(26,100){\small{$t = 1000$}}
\end{overpic}
\begin{overpic}[width=.18\textwidth,viewport=160 140 660 670, clip,grid=false]{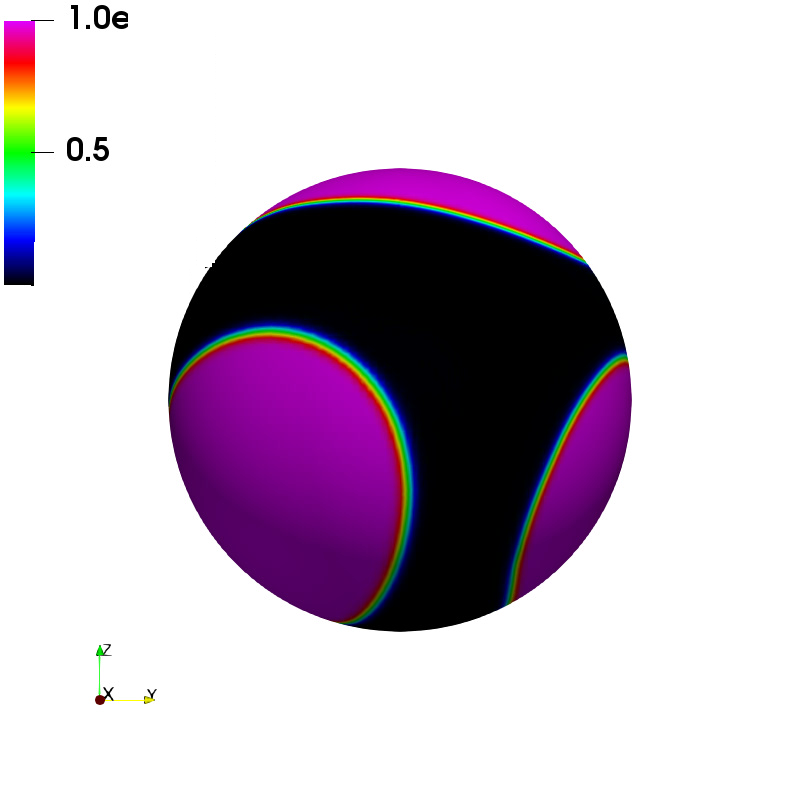}
        \put(26,100){\small{$t = 2000$}}
\end{overpic}
\begin{overpic}[width=.18\textwidth,viewport=160 140 660 670, clip,grid=false]{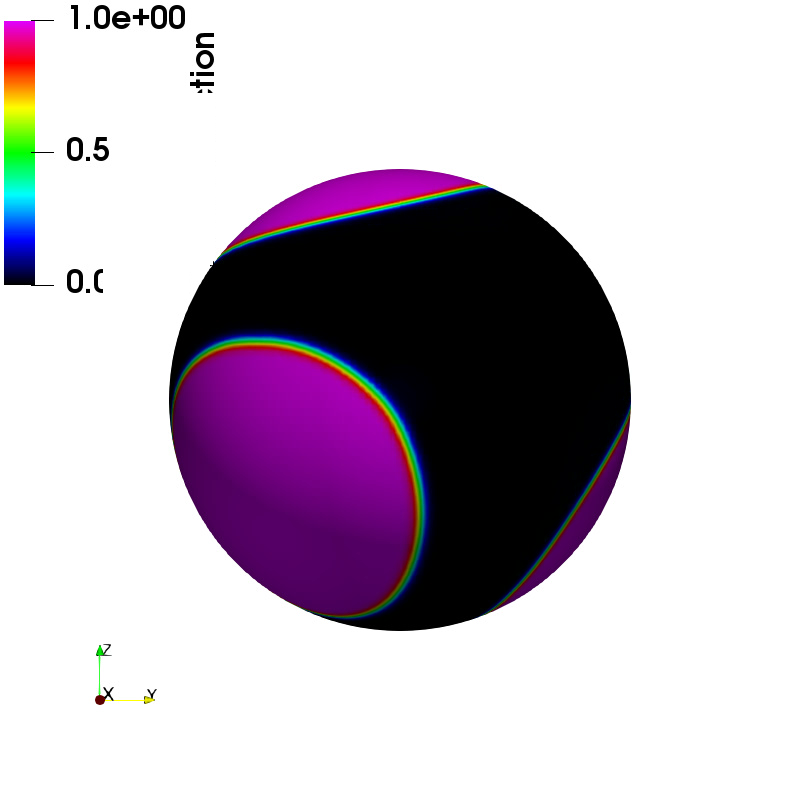}
        \put(26,100){\small{$t = 8000$}}
\end{overpic}
\begin{overpic}[width=.18\textwidth,viewport=160 140 660 670, clip,grid=false]{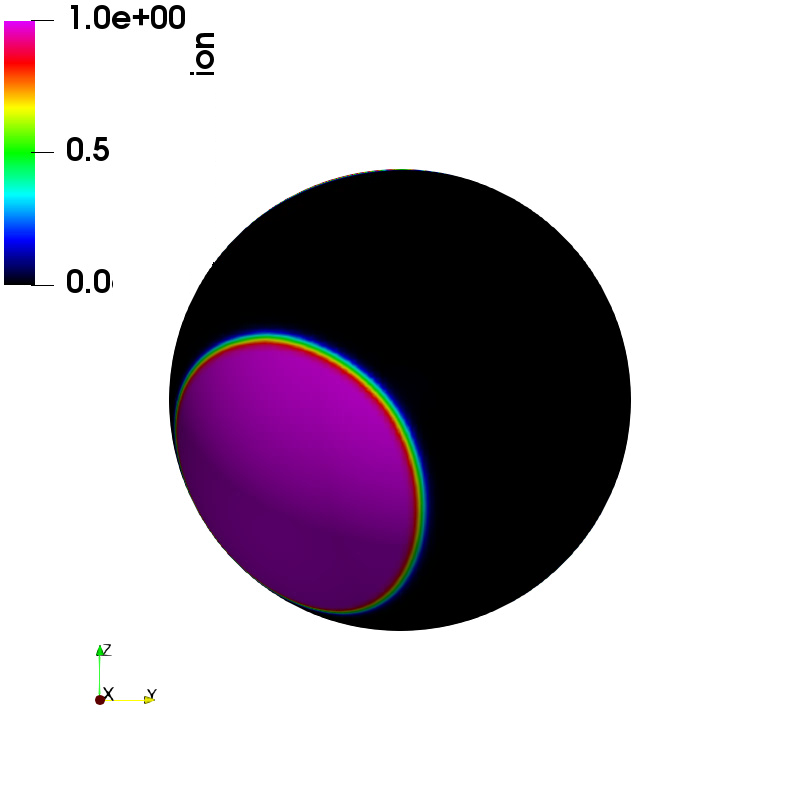}
        \put(26,100){\small{$t = 15000$}}
\end{overpic}
\begin{overpic}[width=.18\textwidth,viewport=160 140 660 670, clip,grid=false]{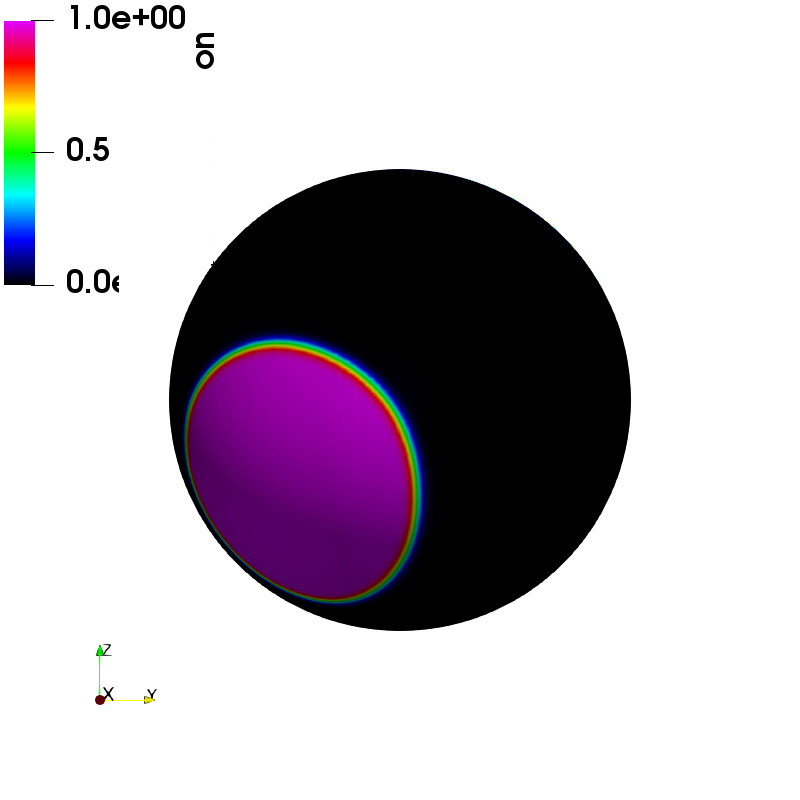}
        \put(26,100){\small{$t = 25000$}}
\end{overpic}
\end{center}
\caption{\label{sphereCH_10_250}
Evolution of the numerical solution of the Cahn--Hilliard equation on the sphere for $t \in [1000, 25000]$
computed with mesh $\ell = 6$ and time step $\Delta t =1$.
View: $yz$-plane. The legend is as in Fig.~\ref{sphereCH_0_5}}
\end{figure}

\begin{figure}
\begin{center}
\begin{overpic}[width=.20\textwidth,viewport=160 140 660 670, clip,grid=false]{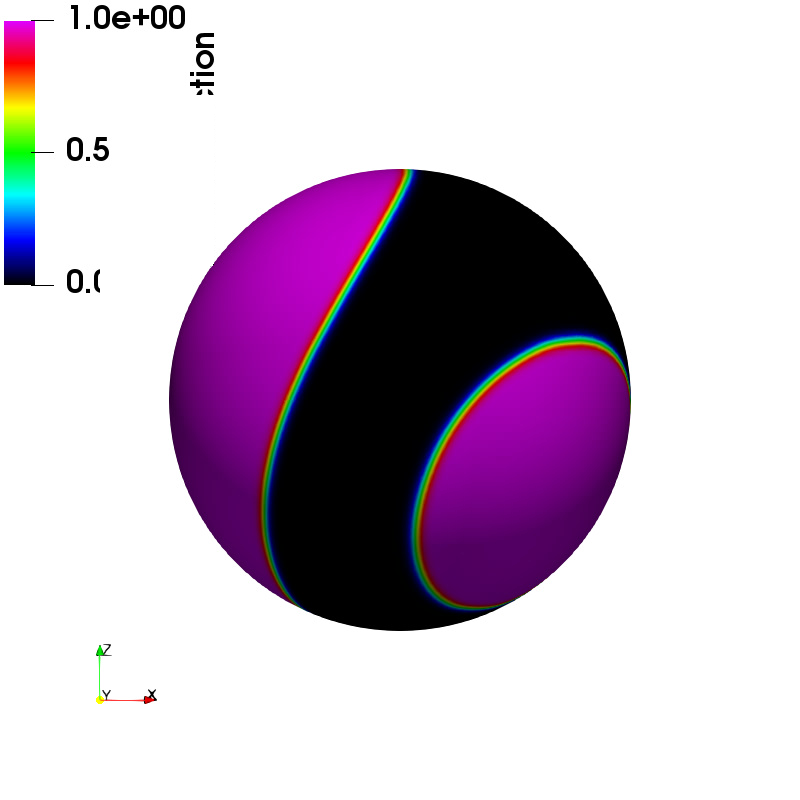}
        \put(26,100){\small{$t = 16000$}}
\end{overpic}
\begin{overpic}[width=.20\textwidth,viewport=160 140 660 670, clip,grid=false]{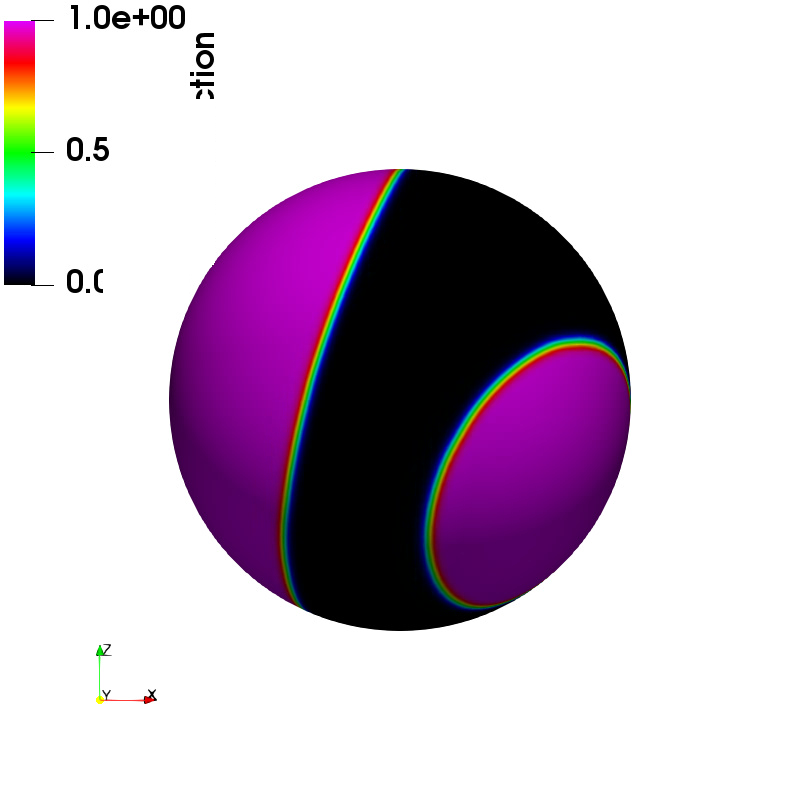}
        \put(26,100){\small{$t = 19000$}}
\end{overpic}
\begin{overpic}[width=.20\textwidth,viewport=160 140 660 670, clip,grid=false]{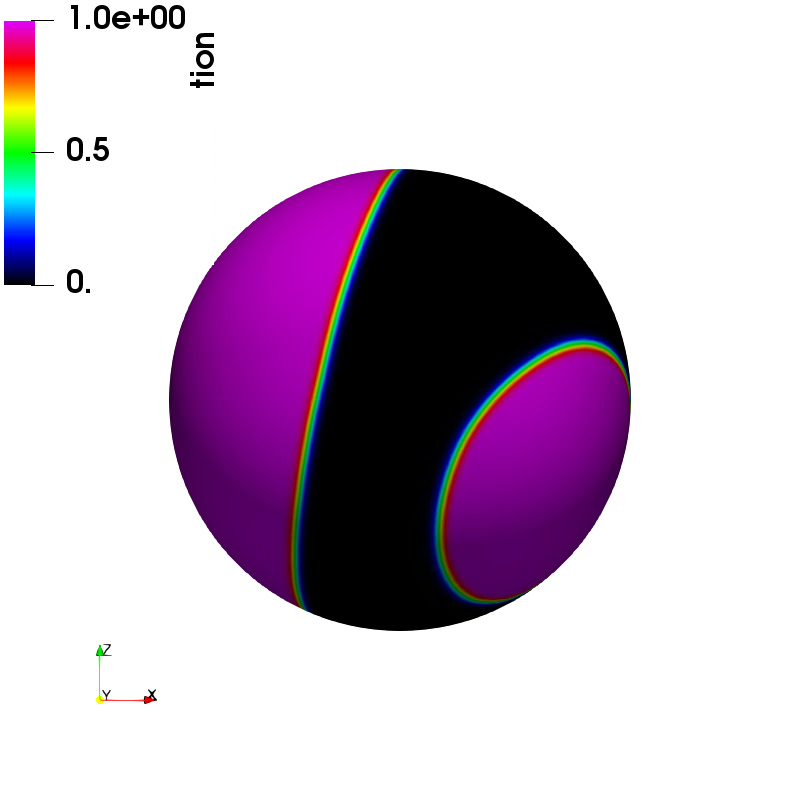}
        \put(26,100){\small{$t = 22000$}}
\end{overpic}
\begin{overpic}[width=.20\textwidth,viewport=160 140 660 670, clip,grid=false]{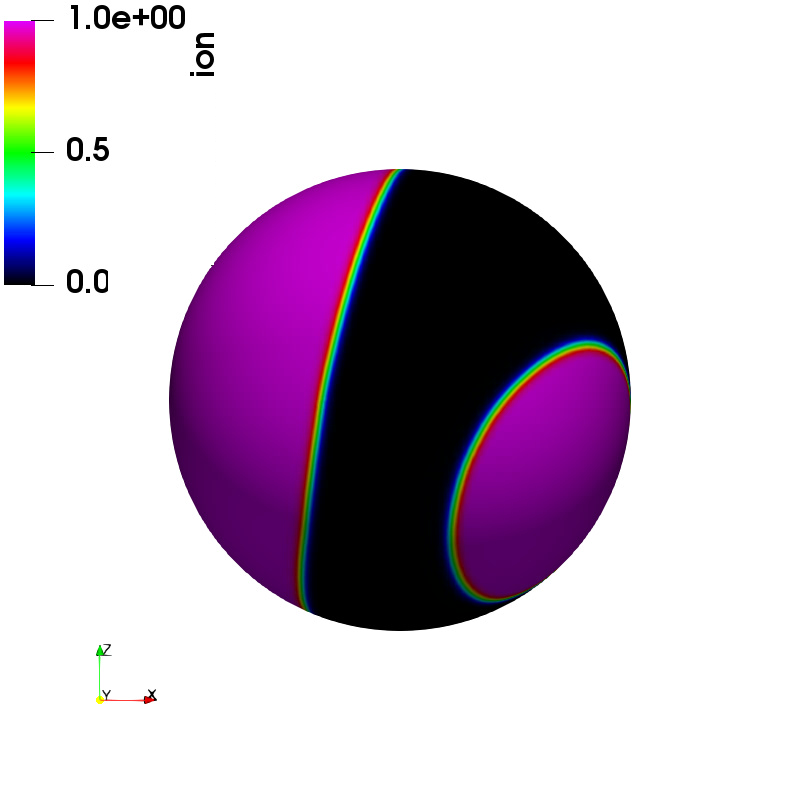}
        \put(26,100){\small{$t = 25000$}}
\end{overpic}
\hspace{.3cm}
\begin{overpic}[width=.10\textwidth,grid=false]{figures_AC_sphere_legend.png}
\end{overpic}
\end{center}
\caption{\label{sphereCH_160_250xz}
Evolution of the numerical solution of the Cahn--Hilliard equation on the sphere for $t \in [16000, 25000]$
computed with mesh $\ell = 6$ and time step $\Delta t =1$.
View: $xz$-plane.}
\end{figure}

The evolution of the numerical solution to the Cahn--Hilliard problem
for $t \in [1000, 25000]$ is shown in Fig.~\ref{sphereCH_10_250}. We see
that little changes in the solution occur between $t = 15000$ and $t = 25000$.
After $t = 25000$, there is no visible change in the position of the interface
between phases. In order to give a better idea of the equilibrium,
we show another view of the solution for $t \in [16000, 25000]$ in Fig.~\ref{sphereCH_160_250xz}.
Close to the steady state, we observe one large and one small pink domain (i.e., $c = 1$).
We recall that the specific free energy \eqref{eq:f0} for the Cahn--Hilliard equation
is not convex, this being a necessary condition to have phase separation.
Thus, more than one stable equilibrium might exist {solely due to the definition
of the specific free energy, not considering surface symmetries}.

In order to follow the initial fast phase, we set $\Delta t = 0.01$ for $t \in (0, 1]$. For the remaining
part of the considered time interval  $\Delta t$ was set to 1. Obviously,
also for the Cahn--Hilliard problem a time-adaptivity strategy can be adopted. See, e.g.,  \cite{WODO20116037,GOMEZ20115310,zhang_qiao_2012,guillen2014second}

\subsubsection{{On stabilization parameter $\beta_s$}}\label{sec:beta_s}

The fully discrete problems \eqref{eq:AC_FE} and \eqref{eq:CH_FE} feature
user defined stabilization parameter $\beta_s$. As a rule of thumb on how to
set its value, in \cite{Shen_Yang2010} it is suggested that $\beta_s$ should not be ``too small''
in order to relax the severity of the restriction for time step
that comes from the analysis.
See also Sec.~\ref{sec:num_method}. In this section, we report
some numerical results to further support the suggestion in \cite{Shen_Yang2010}.

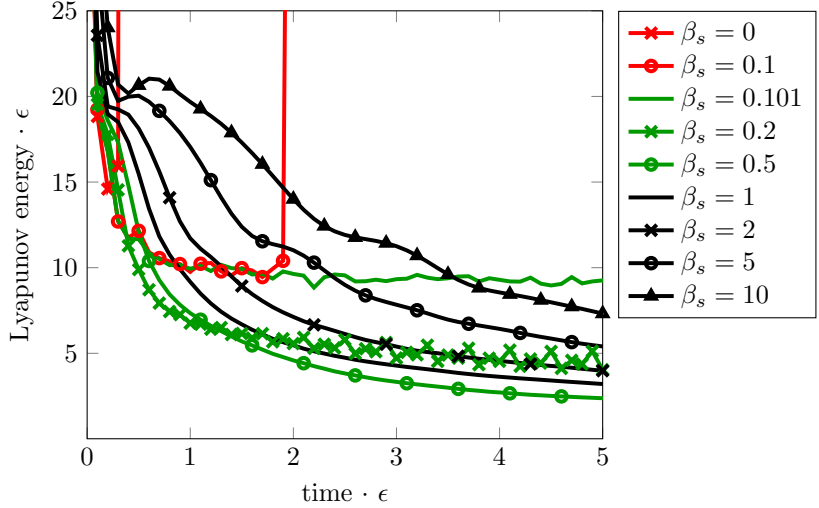
\begin{figure}[h]
		\begin{center}
			\def\vara{10}
			\def\varb{5}
			\begin{tikzpicture}
			\begin{axis}[ xlabel={time $\cdot$ $\epsilon$}, ylabel={Lyapunov energy $\cdot$ $\epsilon$}, ymin=1E-2, xmax=5,xmin=0,
			, ymax=25, legend style={ cells={anchor=west}, legend pos=outer north east} ]
			
			\addplot+[red,mark=x,solid,mark size=3pt,line width=1.5pt] table[x=Time, y=Lyapunov_energy] {energy_AC_S=0_Plot_test9_l=6_BDF2=0.100000.txt};

			\addplot+[red,mark=o,mark repeat={2},solid,mark size=2pt,line width=1.5pt] table[x=Time, y=Lyapunov_energy] {energy_AC_S=0.1_Plot_test9_l=6_BDF2=0.100000.txt};
			
			\addplot+[black!40!green,mark=t,mark repeat={1},solid,mark size=3pt,line width=1.5pt] table[x=Time, y=Lyapunov_energy] {energy_AC_S=0.101_Plot_test9_l=6_BDF2=0.100000.txt};

			\addplot+[black!40!green,mark=x,solid,mark size=3pt,line width=1.5pt] table[x=Time, y=Lyapunov_energy] {energy_AC_S=0.2_Plot_test9_l=6_BDF2=0.100000.txt};
			
			\addplot+[black!40!green,mark=o,solid,mark size=2pt,mark repeat={5},line width=1.5pt] table[x=Time, y=Lyapunov_energy] {energy_AC_S=0.5_Plot_test9_l=6_BDF2=0.100000.txt};

			\addplot+[black,mark=t,solid,mark size=1pt, mark repeat={9},line width=1.5pt] table[x=Time, y=Lyapunov_energy] {energy_AC_S=1_Plot_test9_l=6_BDF2=0.100000.txt};
			
			\addplot+[black,mark=x,solid,mark repeat={7},mark size=3pt,line width=1.5pt] table[x=Time, y=Lyapunov_energy] {energy_AC_S=2_Plot_test9_l=6_BDF2=0.100000.txt};
			
			\addplot+[black,mark=o,solid, mark repeat={5},mark size=2pt,line width=1.5pt] table[x=Time, y=Lyapunov_energy] {energy_AC_S=5_Plot_test9_l=6_BDF2=0.100000.txt};

			\addplot+[black,mark=triangle,solid,mark repeat={3},mark size=2pt,line width=1.5pt] table[x=Time, y=Lyapunov_energy] {energy_AC_S=10_Plot_test9_l=6_BDF2=0.100000.txt};

			\legend{
				$\beta_s=0$ ,
				$\beta_s=0.1$ ,
				$\beta_s=0.101$ ,
				$\beta_s=0.2$ ,							$\beta_s=0.5$ ,
				$\beta_s=1$ ,
				$\beta_s=2$ ,
				$\beta_s=5$ ,
				$\beta_s=10$
			}
			\end{axis}
			\end{tikzpicture}
			\caption{{Evolution of the discrete Lyapunov energy functional \eqref{eq:Lyapunov_E}
			for the Allen--Cahn problem \eqref{eq:AC_FE} posed on the unit sphere with random initial data, $\epsilon=0.01$, $\Delta t=10$, and different values of stabilization parameter $\beta_s$.}}
			\label{fig:E_AC}
		\end{center}
	\end{figure}

\begin{figure}
	\begin{center}
		\def\vara{10}
		\def\varb{5}
		\begin{tikzpicture}
		\begin{axis}[ xlabel={time $\cdot$ $\epsilon$}, ylabel={Lyapunov energy $\cdot$ $\epsilon$}, ymin=1E-2, xmax=5,xmin=0,
		, ymax=25, legend style={ cells={anchor=west}, legend pos=outer north east} ]
		
		\addplot+[red,mark=t,solid,mark size=4pt,line width=1.5pt] table[x=Time, y=Lyapunov_energy] {energy_CH_S=0_Plot_test9_l=6_BDF2=0.010000.txt};

		\addplot+[red,mark=x,mark repeat={2},solid,mark size=3pt,line width=1.5pt] table[x=Time, y=Lyapunov_energy] {energy_CH_S=0.1_Plot_test9_l=6_BDF2=0.010000.txt};
		
		\addplot+[red,mark=o,mark repeat={2},solid,mark size=2pt,line width=1.5pt] table[x=Time, y=Lyapunov_energy] {energy_CH_S=0.28_Plot_test9_l=6_BDF2=0.010000.txt};
			
		\addplot+[black!40!green,mark=t,solid,mark size=1.5pt,line width=1.5pt,mark repeat={5}] table[x=Time, y=Lyapunov_energy] {energy_CH_S=0.29_Plot_test9_l=6_BDF2=0.010000.txt};
			
		\addplot+[black!40!green,mark=x,solid,mark size=3pt, mark repeat={9},line width=1.5pt] table[x=Time, y=Lyapunov_energy] {energy_CH_S=0.5_Plot_test9_l=6_BDF2=0.010000.txt};
		
		\addplot+[black!40!green,mark=o,solid,mark size=2pt, mark repeat={9},line width=1.5pt] table[x=Time, y=Lyapunov_energy] {energy_CH_S=1_Plot_test9_l=6_BDF2=0.010000.txt};
				
		\addplot+[black,mark=t,solid,mark repeat={3},mark size=1.5pt,line width=1.5pt] table[x=Time, y=Lyapunov_energy] {energy_CH_S=10_Plot_test9_l=6_BDF2=0.010000.txt};

		\legend{
			$\beta_s=0$ ,
			$\beta_s=0.1$ ,
			$\beta_s=0.28$ ,		
			$\beta_s=0.29$ ,		
			$\beta_s=0.5$ ,		
			$\beta_s=1$ ,
			$\beta_s=10$
		}
		\end{axis}
		\end{tikzpicture}
		\caption{{Evolution of the discrete Lyapunov energy functional \eqref{eq:Lyapunov_E}
			for the Cahn--Hilliard problem \eqref{eq:CH_FE} posed on the unit sphere with random initial data, $\epsilon=0.01$, $\Delta t=1$, and different values of stabilization parameter $\beta_s$.}}
		\label{fig:E_CH}
	\end{center}
\end{figure}
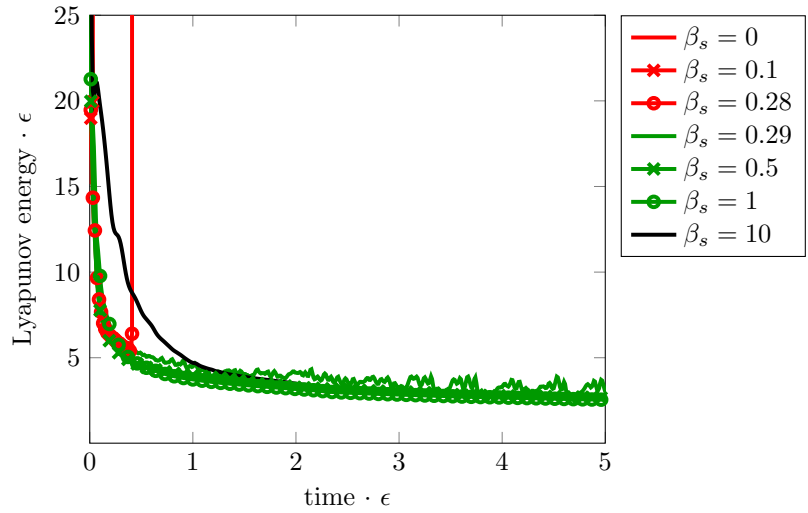

Let
\begin{equation}\label{eq:Lyapunov_E}
E^L_h(u_h) = \int_{\Gamma_h} f(u_h) ds = \int_{\Gamma_h} \left( f_0(u_h) + \frac{1}{2} \epsilon^2 | \gradG u_h |^2 \right)ds
\end{equation}
be the discrete Lyapunov energy functional for the Allen-Cahn ($u_h = \eta_h$) and Cahn--Hilliard ($u_h = c_h$)
problems. We consider again both problems posed on the unit sphere, with $\epsilon = 0.01$,
and initial condition given by $\text{rand}(\bx)$. We let the value of $\beta_s$ range from 0 to 10.
Fig.~\ref{fig:E_AC} shows the evolution of $E^L_h(\eta_h)$ over time for $\Delta t = 10$,
while Fig.~\ref{fig:E_CH} displays the evolution of $E^L_h(c_h)$ over time for $\Delta t = 1$.
We see that the discrete Lyapunov energy functional blows up when $\beta_s$ is small
with respect to the time step, meaning $\beta_s \leq 0.1$ for the Allen-Cahn problem
and $\beta_s \leq 0.28$ for the Cahn--Hilliard problem.

We remind the stability conditions in \cite{Shen_Yang2010} are independent of $\beta_s$. The
suggestion for the value of $\beta_s$ comes from experimental practice. The numerical experiments 
in \cite{Shen_Yang2010} refer to problems \eqref{eq:AC_FE} and \eqref{eq:CH_FE} in planar domains.
The Cahn--Hilliard problem in \cite{Shen_Yang2010} is simplified by setting mobility to 1. The results reported in Fig.~\ref{fig:E_AC} and \ref{fig:E_CH} show that the rule of thumb suggested
in \cite{Shen_Yang2010} holds also for non-planar surfaces and concentration dependent mobility.

\subsection{Phase separation on a spindle surface}\label{sec:spindle}

The surface $\Gamma$ is given by the zero of the level set function
\begin{equation*}
\phi(\bx) = \left\{\begin{split}
                      \cos^2 \left(\frac{x_1\pi}{10}\right)-(4x_2)^2-(4x_3)^2 &\quad\text{if}~x_1<5  \\
                       -25  & \quad\text{if}~x_1>5
                   \end{split} \right. .
\end{equation*}
This corresponds to the surface of a spindle with maximum radius 0.25.
We consider this surface because it resembles that of the bacteria studied in
\cite{Bramkamp_Lopez2015}. The surface $\Gamma$ is embedded in an
outer domain $\Omega=  [-5,5] \times [-4/3, 4/3] \times [-4/3, 4/3]$.
A computational mesh $\T_{h}^\Gamma$ is generated with mesh size $h = 0.026$.
This is a level of refinement comparable to mesh $\ell = 6$ in Sec.~\ref{sec:NumVal},
which we have seen is appropriate for interface thickness $\epsilon = 0.01$ or bigger.

In \cite{Bramkamp_Lopez2015}, it is reported that bacteria organize many
processes in functional membrane microdomains equivalent to the lipid rafts.
Thus, we simulate phase separation on the spindle
using both models, i.e. eq.~\eqref{eq:AC_eq} and \eqref{eq:CH}, with free energy per unit surface \eqref{eq:f0},
$\alpha = 1$, and $\epsilon = 0.01$.
The initial condition is $\eta_0(\bx)=\text{rand}(\bx)$, where $\text{rand}(\bx)$ is a uniformly distributed
random number between 0 and 1.

\subsubsection{Allen--Cahn model}

\begin{figure}
\begin{center}
\begin{overpic}[width=.07\textwidth, viewport=410 0 530 780, clip,grid=false]{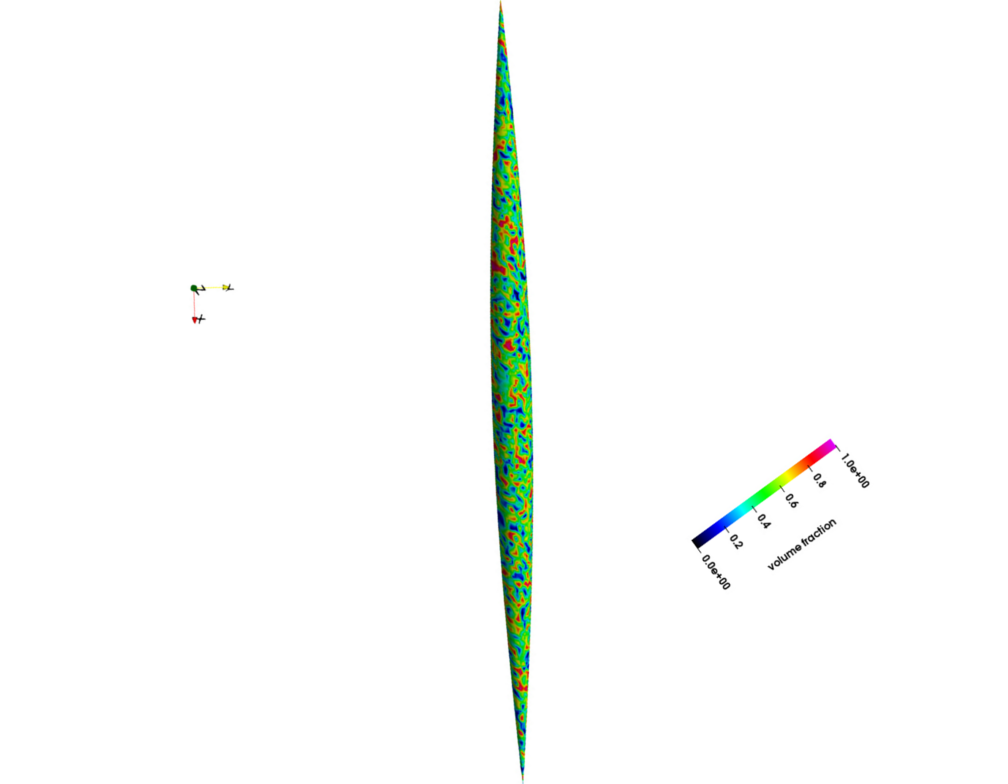}
         \put(5,103){\small{$t = 0$}}
\end{overpic}
\begin{overpic}[width=.07\textwidth, viewport=410 0 530 780, clip,grid=false]{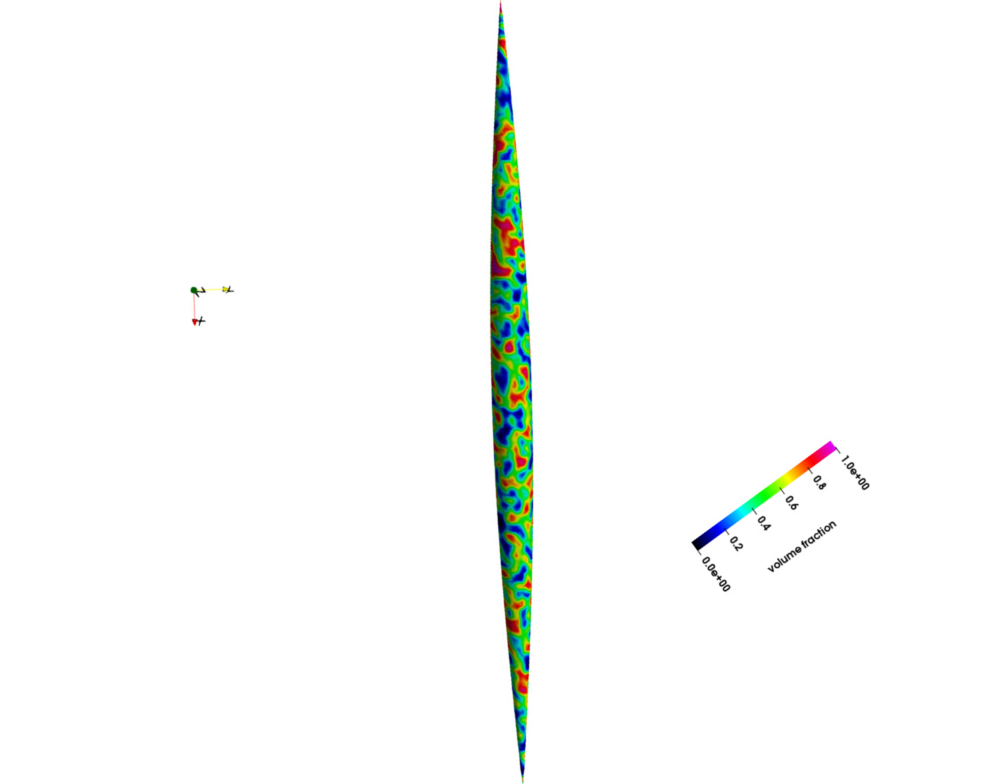}
         \put(3,103){\small{$t = 5$}}
\end{overpic}
\begin{overpic}[width=.07\textwidth, viewport=410 0 530 780, clip,grid=false]{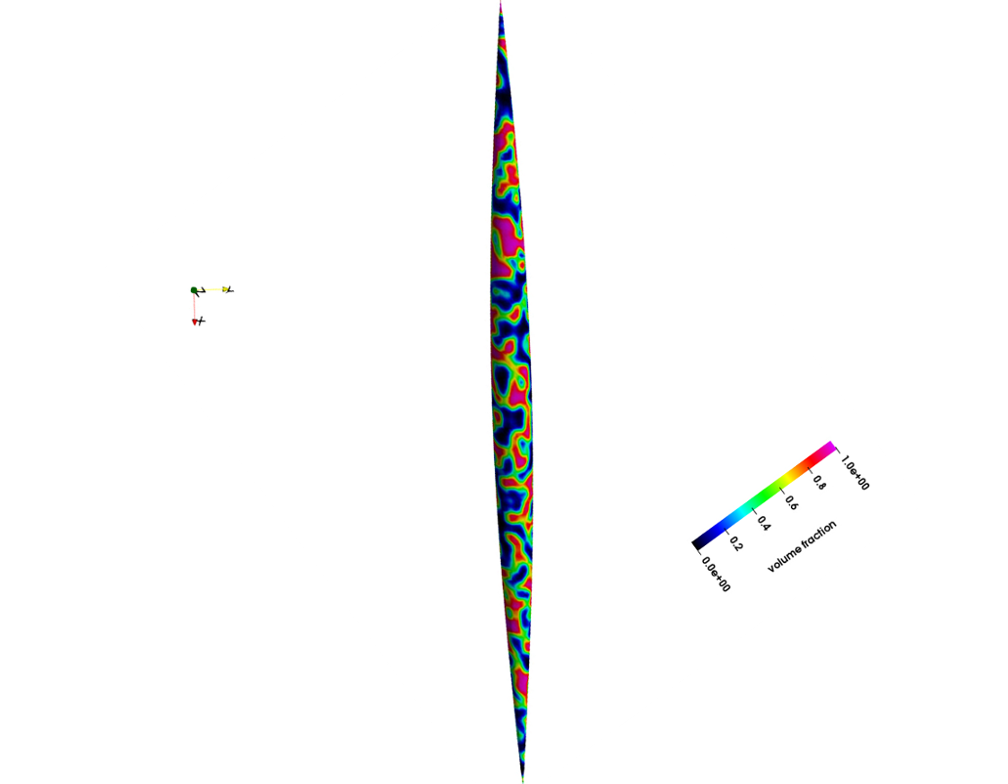}
         \put(3,103){\small{$t = 10$}}
\end{overpic}
\begin{overpic}[width=.07\textwidth, viewport=410 0 530 780, clip,grid=false]{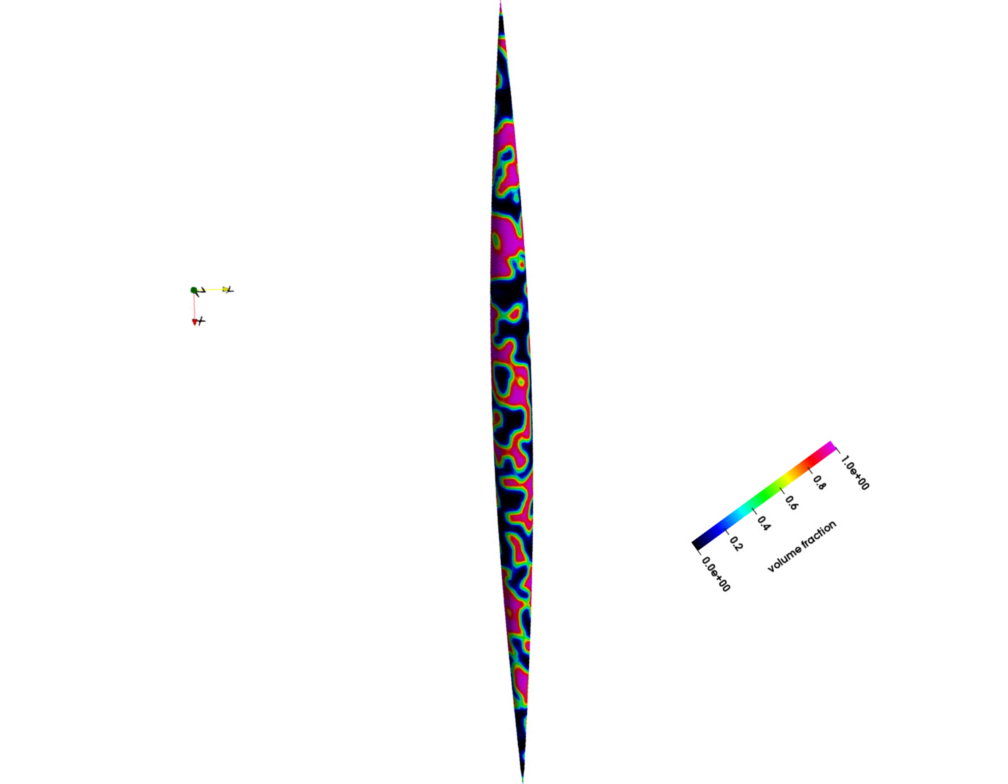}
         \put(3,103){\small{$t = 25$}}
\end{overpic}
\begin{overpic}[width=.07\textwidth, viewport=410 0 530 780, clip,grid=false]{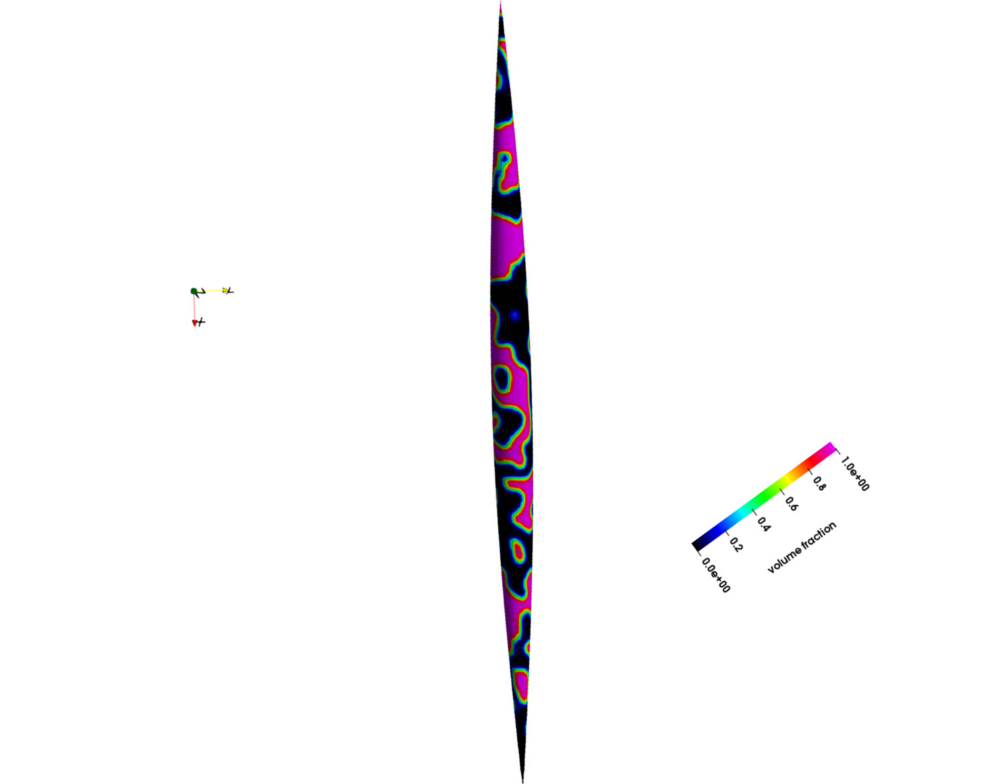}
         \put(3,103){\small{$t = 50$}}
\end{overpic}
\begin{overpic}[width=.07\textwidth, viewport=410 0 530 780, clip,grid=false]{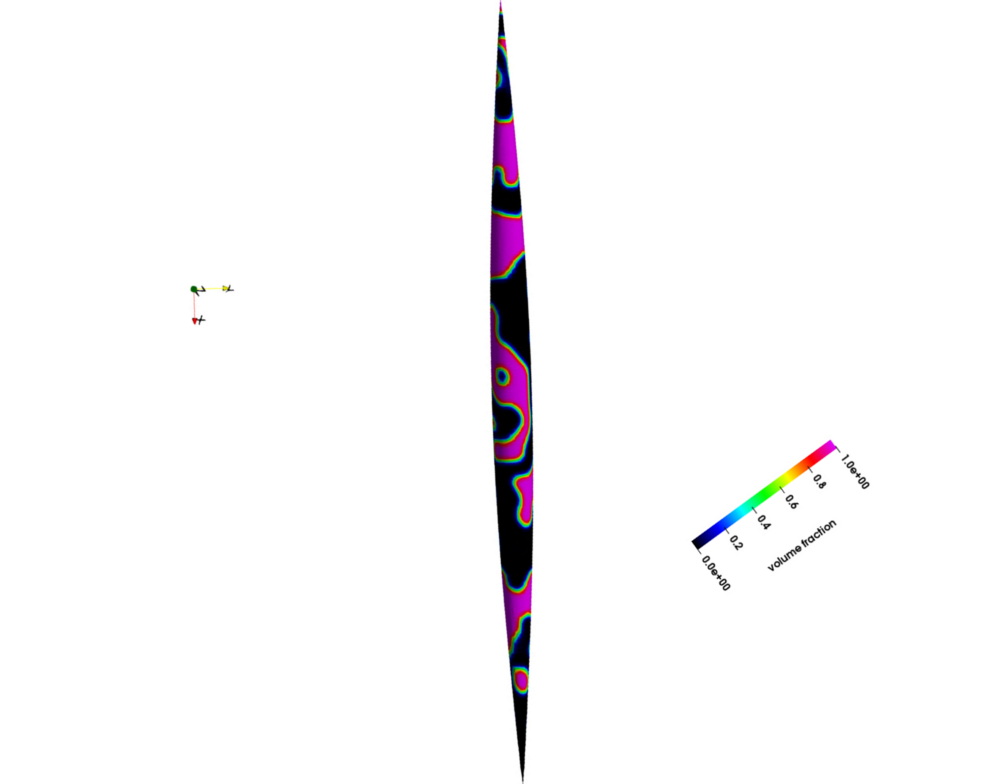}
         \put(3,103){\small{$t = 100$}}
\end{overpic}~
\begin{overpic}[width=.07\textwidth, viewport=410 0 530 780, clip,grid=false]{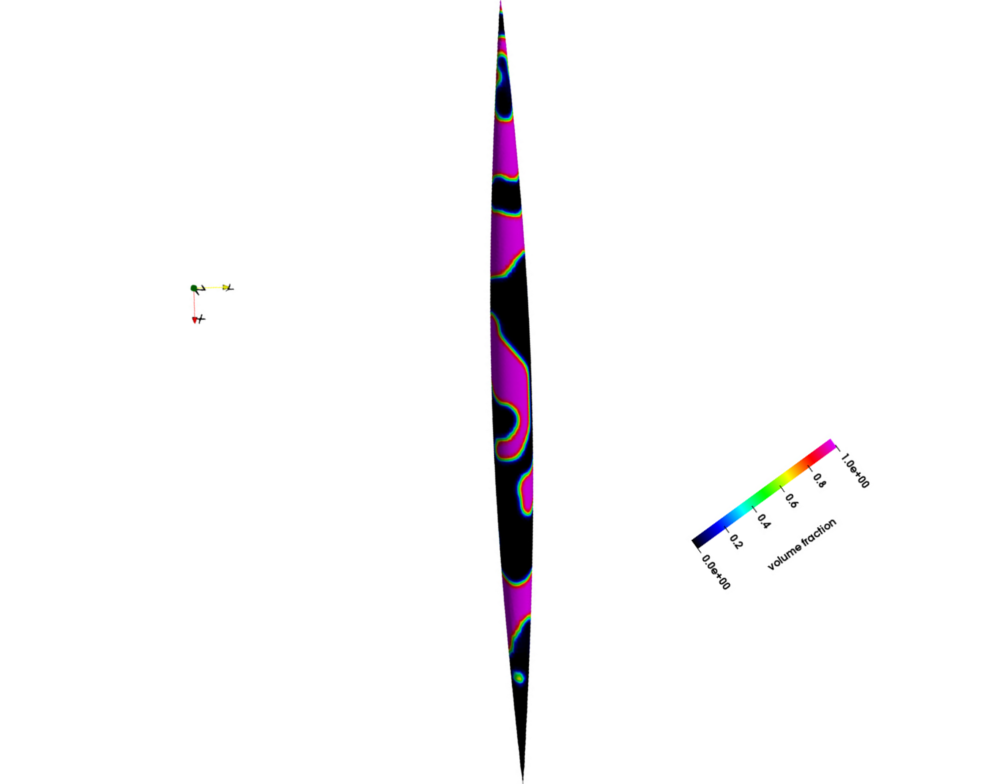}
         \put(3,103){\small{$t = 200$}}
\end{overpic}~
\begin{overpic}[width=.07\textwidth, viewport=410 0 530 780, clip,grid=false]{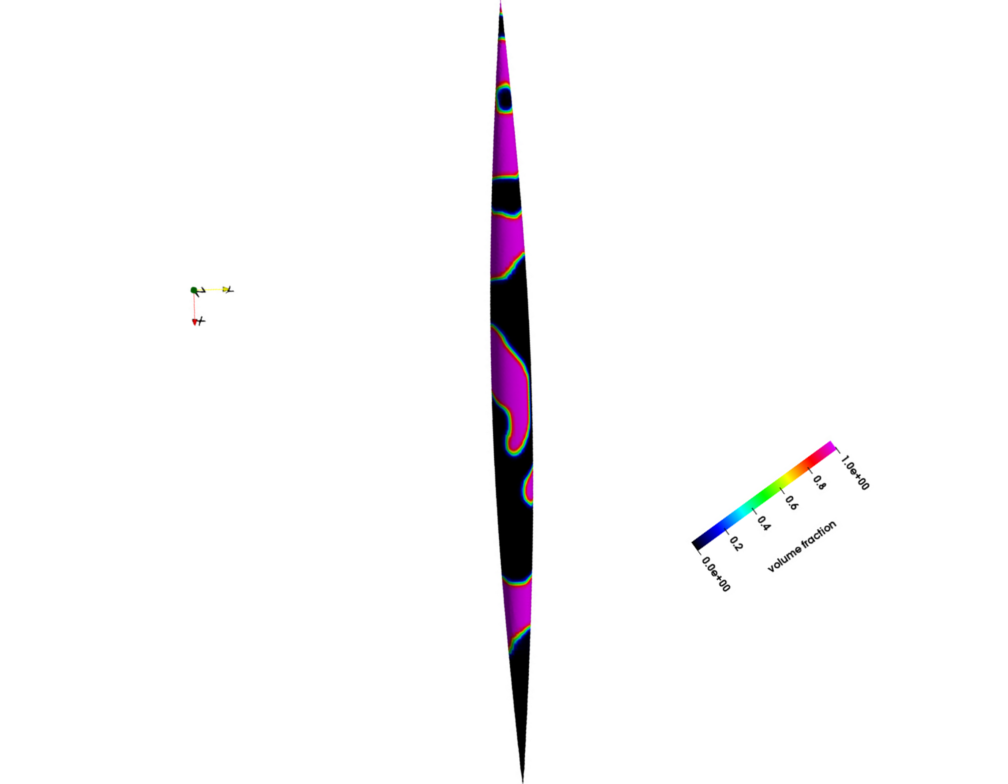}
         \put(3,103){\small{$t = 400$}}
\end{overpic}\hspace{.5cm}
\begin{overpic}[width=.10\textwidth,grid=false]{figures_AC_sphere_legend.png}
\end{overpic}
\end{center}
\caption{\label{spindleAC_0_4}
Evolution of the numerical solution of the Allen--Cahn equation on the spindle for $t \in (0, 400]$
computed with time step $\Delta t = 1$ and on a mesh with size $h =0.026$.
The simulation was started with a random initial condition depicted in the leftmost panel.
View: $xy$-plane.}
\end{figure}

Fig.~\ref{spindleAC_0_4} shows the evolution
of the numerical solution of the Allen--Cahn equation computed with the mesh
described above and $\Delta t = 1$. Again, we observe that
the initial fast stage ends after a little after 10 time units and then the
evolution of the solution slows down considerably. Fig.~1 in \cite{Bramkamp_Lopez2015}
shows a heterogeneous distribution in the cytoplasmic membrane of bacillus subtilis, displaying a
punctate pattern along the entire cell which qualitatively resembles the rightmost panel in Fig.~\ref{spindleAC_0_4}.
We would like to stress that the rightmost panel in Fig.~\ref{spindleAC_0_4} does not correspond to a steady state
of the solution to the Allen--Cahn equation. In fact, if we were to let the simulation continue, it would
evolve towards the disappearance of one phase as observed on the sphere in Sec.~\ref{sec:sphere}.

\subsubsection{Cahn--Hilliard model}\label{sec:CH_spindle}

\begin{figure}
\begin{center}
\hspace{-1cm}
\begin{overpic}[width=.07\textwidth, viewport=410 0 530 790, clip,grid=false]{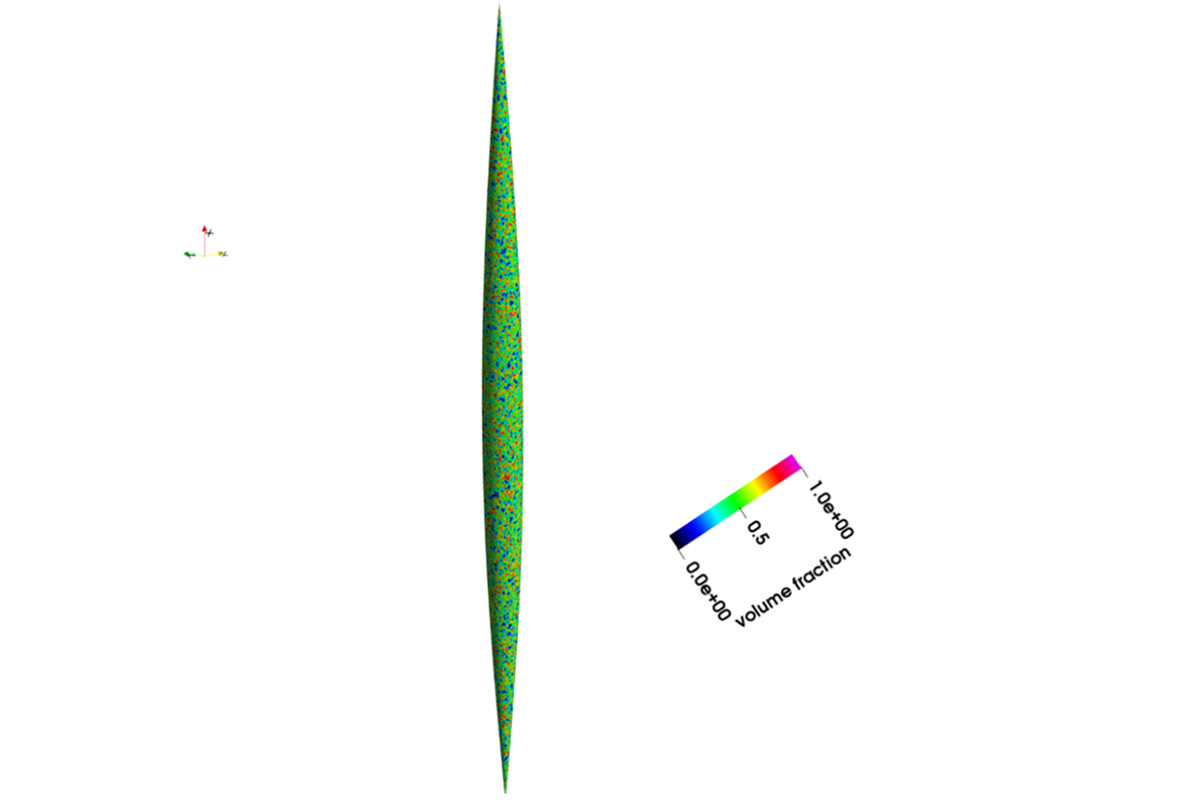}
         \put(5,103){\small{$t = 0$}}
\end{overpic}
\begin{overpic}[width=.07\textwidth, viewport=410 0 530 790, clip,grid=false]{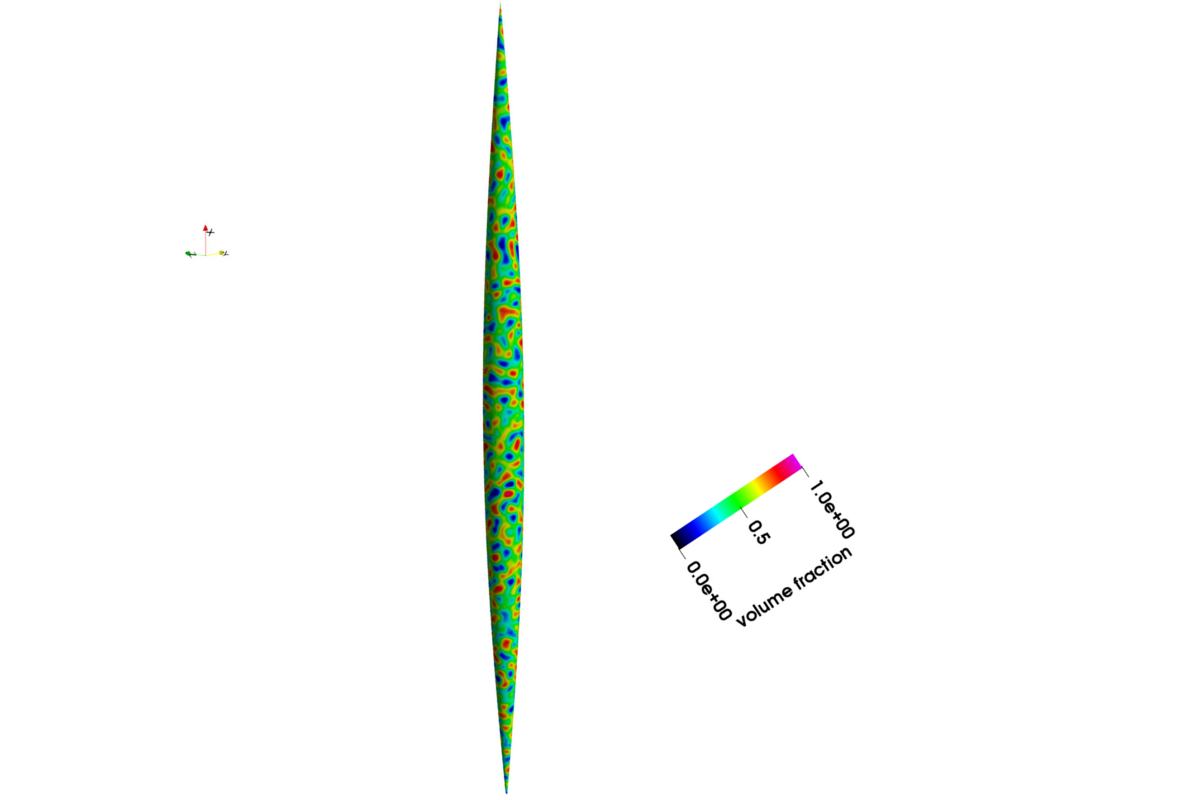}
         \put(3,103){\small{$t = 0.1$}}
\end{overpic}~
\begin{overpic}[width=.07\textwidth, viewport=410 0 530 790, clip,grid=false]{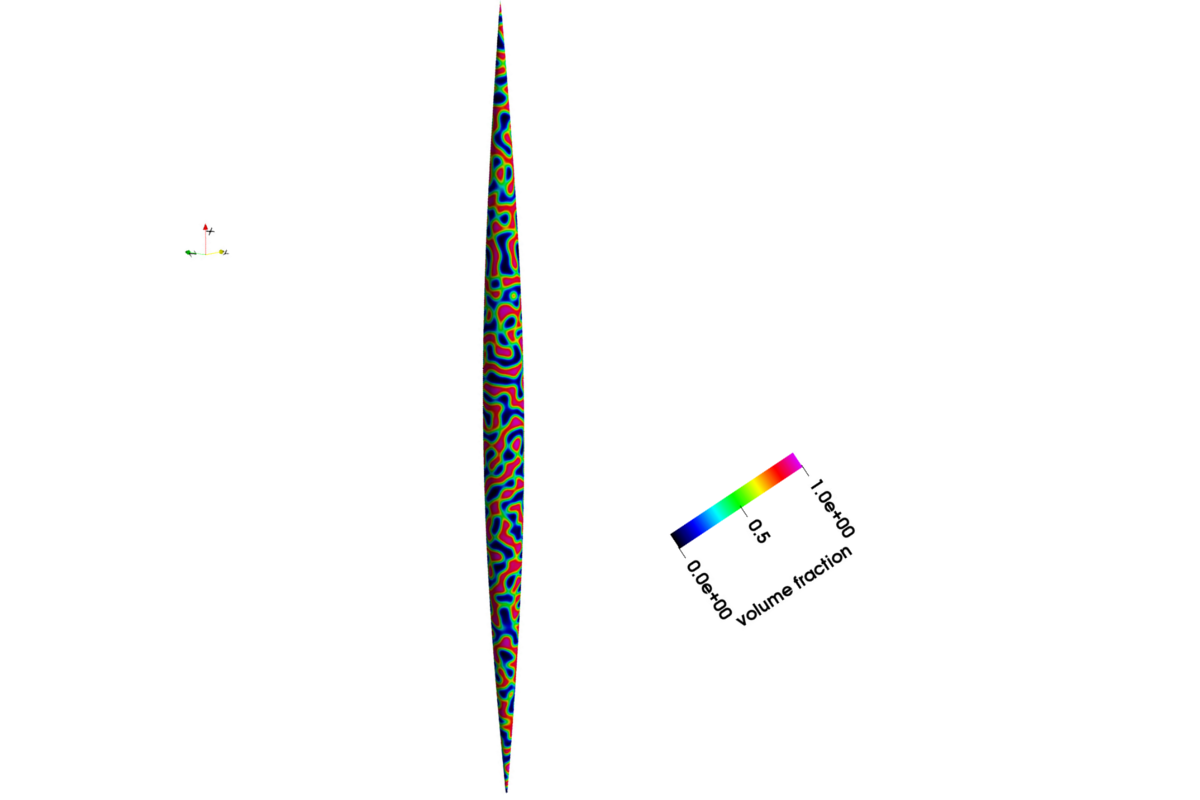}
         \put(3,103){\small{$t = 0.2$}}
\end{overpic}~
\begin{overpic}[width=.07\textwidth, viewport=410 0 530 790, clip,grid=false]{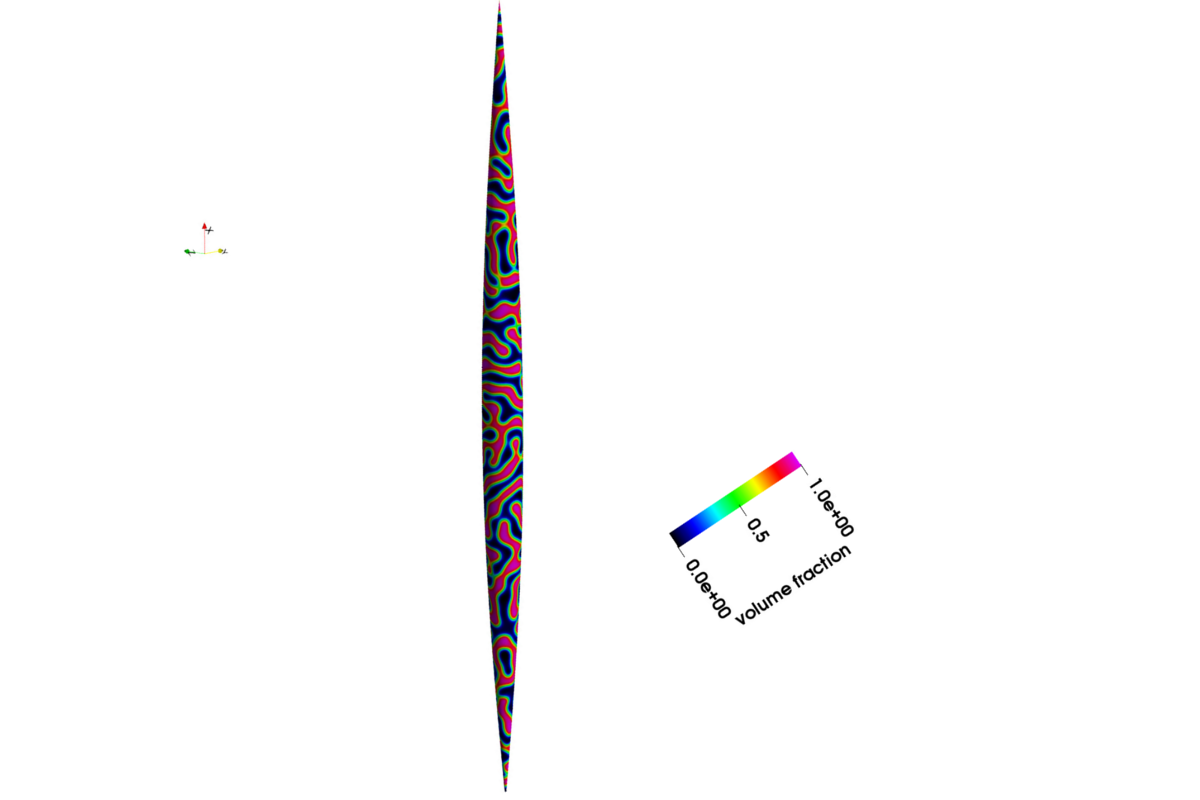}
         \put(3,103){\small{$t = 0.5$}}
\end{overpic}
\begin{overpic}[width=.07\textwidth, viewport=410 0 530 790, clip,grid=false]{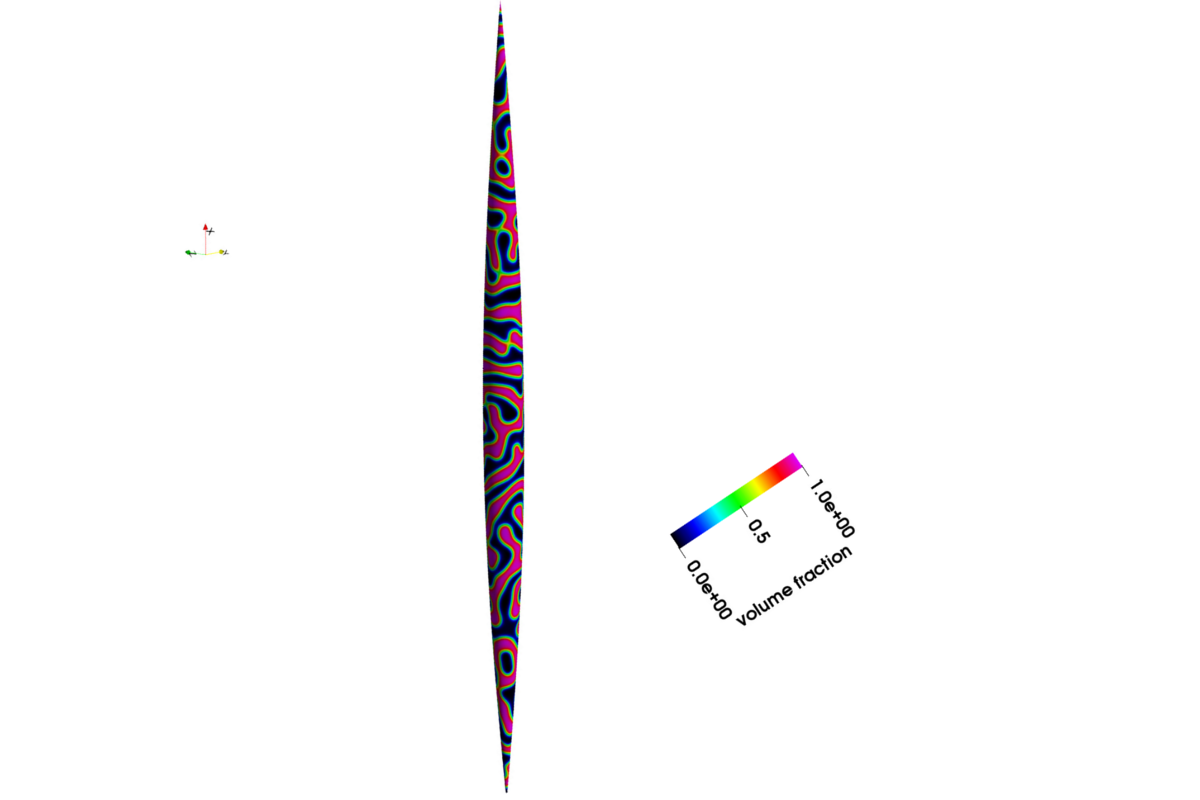}
         \put(6,103){\small{$t = 1$}}
\end{overpic}
\begin{overpic}[width=.07\textwidth, viewport=410 0 530 790, clip,grid=false]{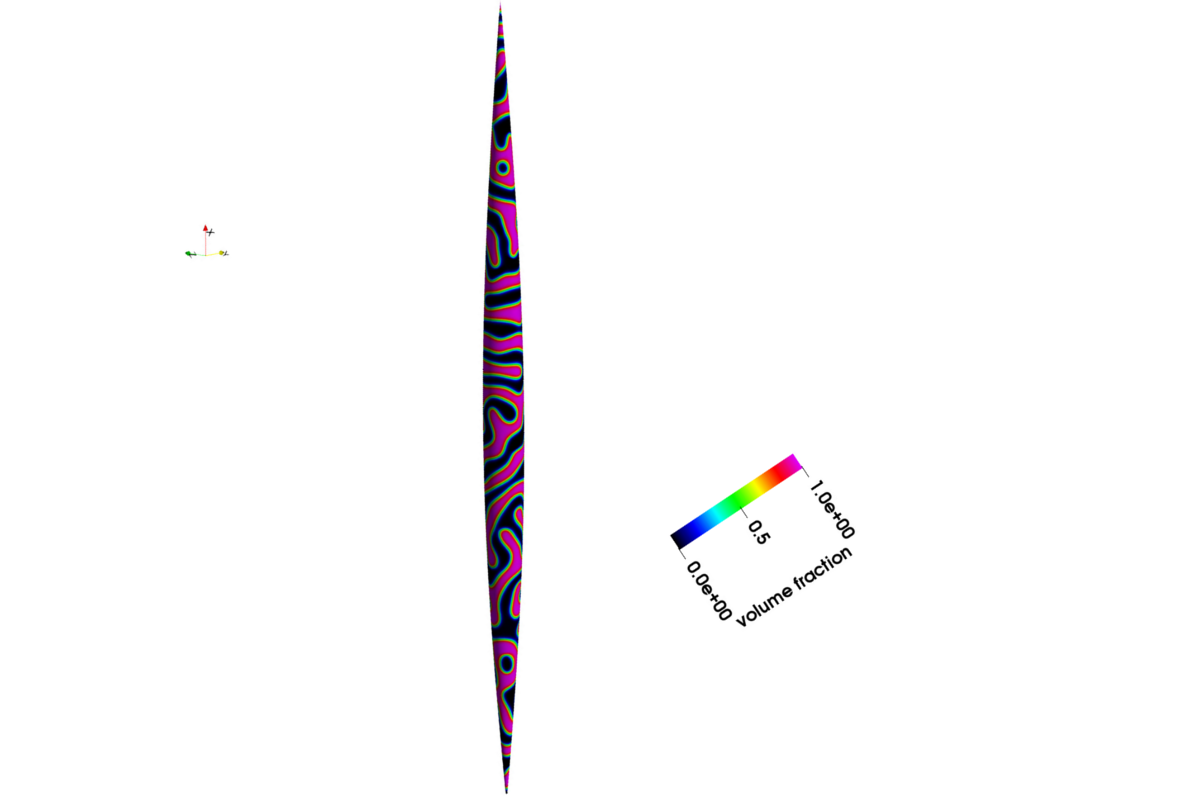}
         \put(4,103){\small{$t = 10$}}
\end{overpic}
\begin{overpic}[width=.07\textwidth, viewport=410 0 530 790, clip,grid=false]{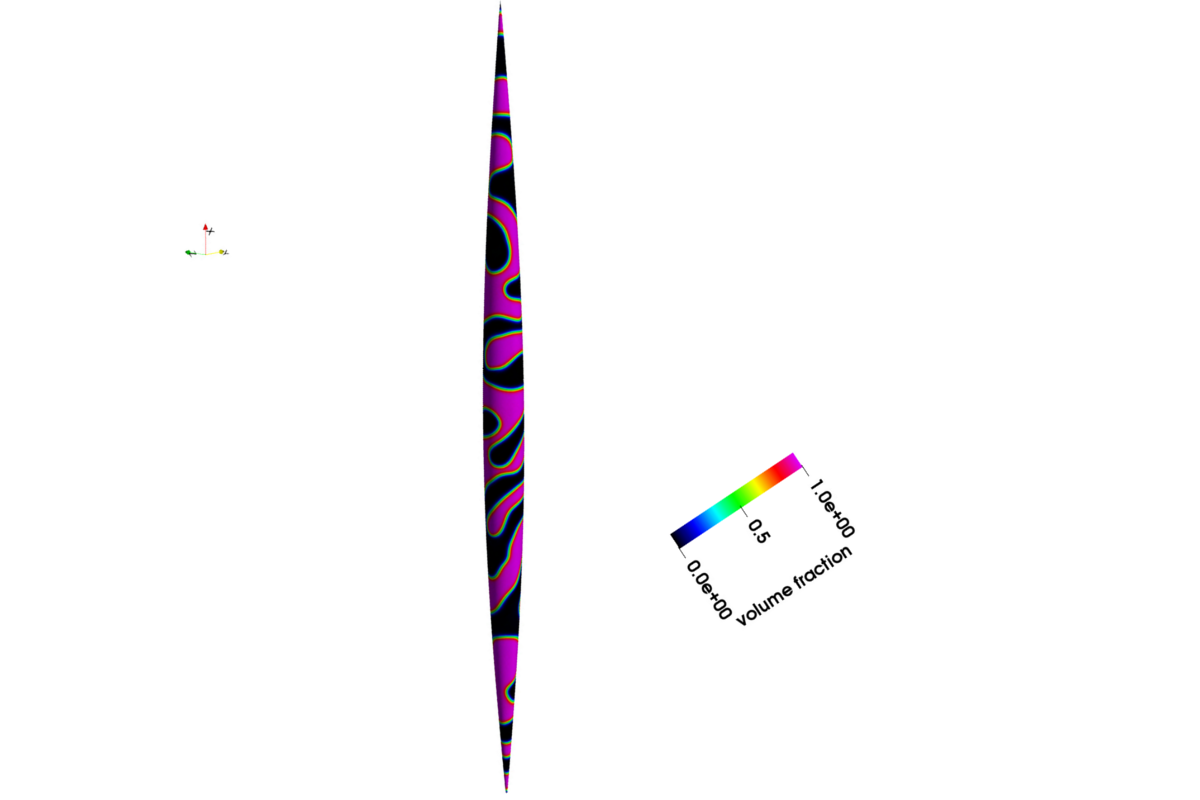}
         \put(4,103){\small{$t = 50$}}
\end{overpic}
\begin{overpic}[width=.07\textwidth, viewport=410 0 530 790, clip,grid=false]{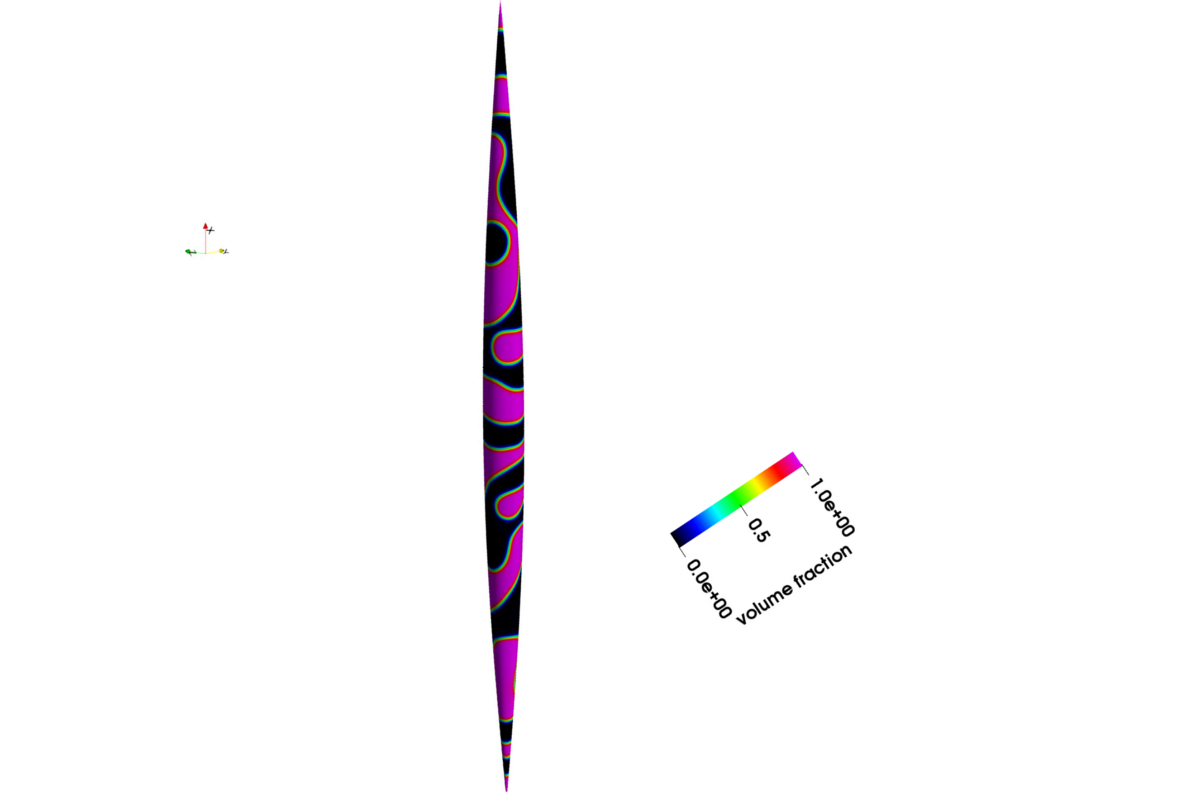}
         \put(3,103){\small{$t = 100$}}
\end{overpic}~
\begin{overpic}[width=.07\textwidth, viewport=410 0 530 790, clip,grid=false]{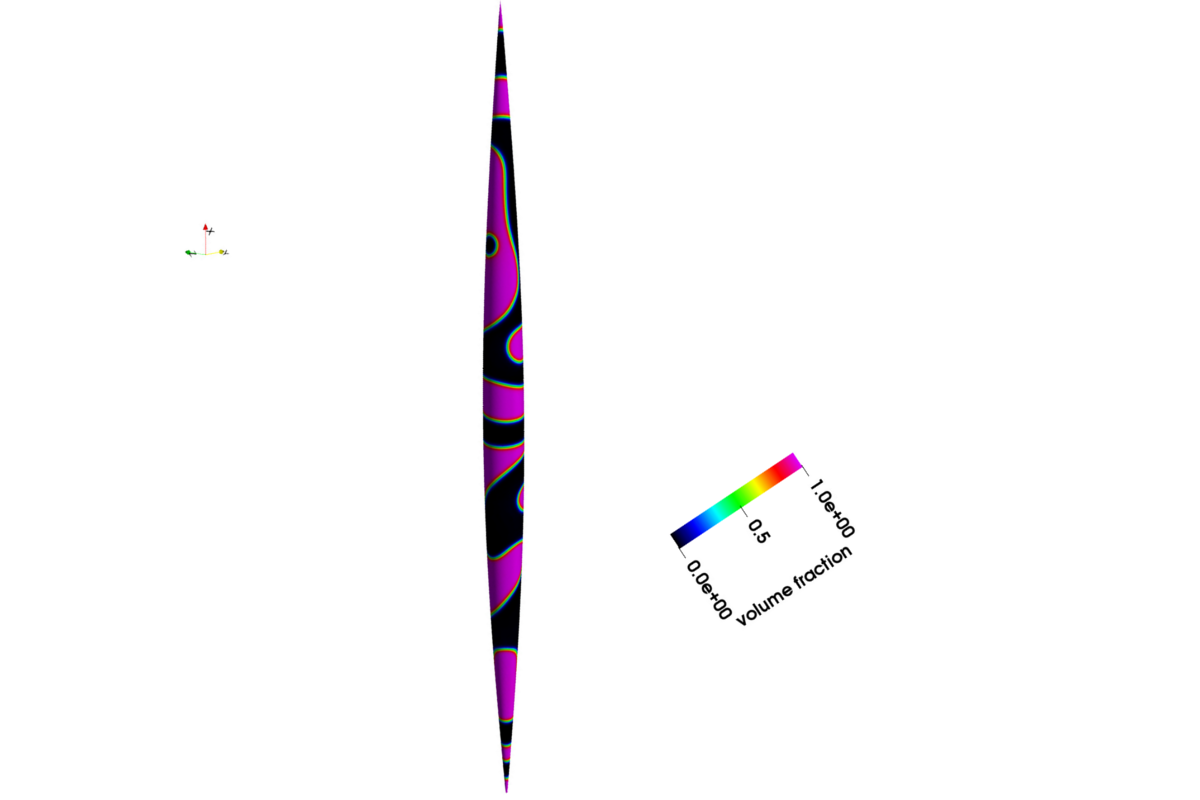}
         \put(3,103){\small{$t = 200$}}
\end{overpic}~
\begin{overpic}[width=.07\textwidth, viewport=410 0 530 790, clip,grid=false]{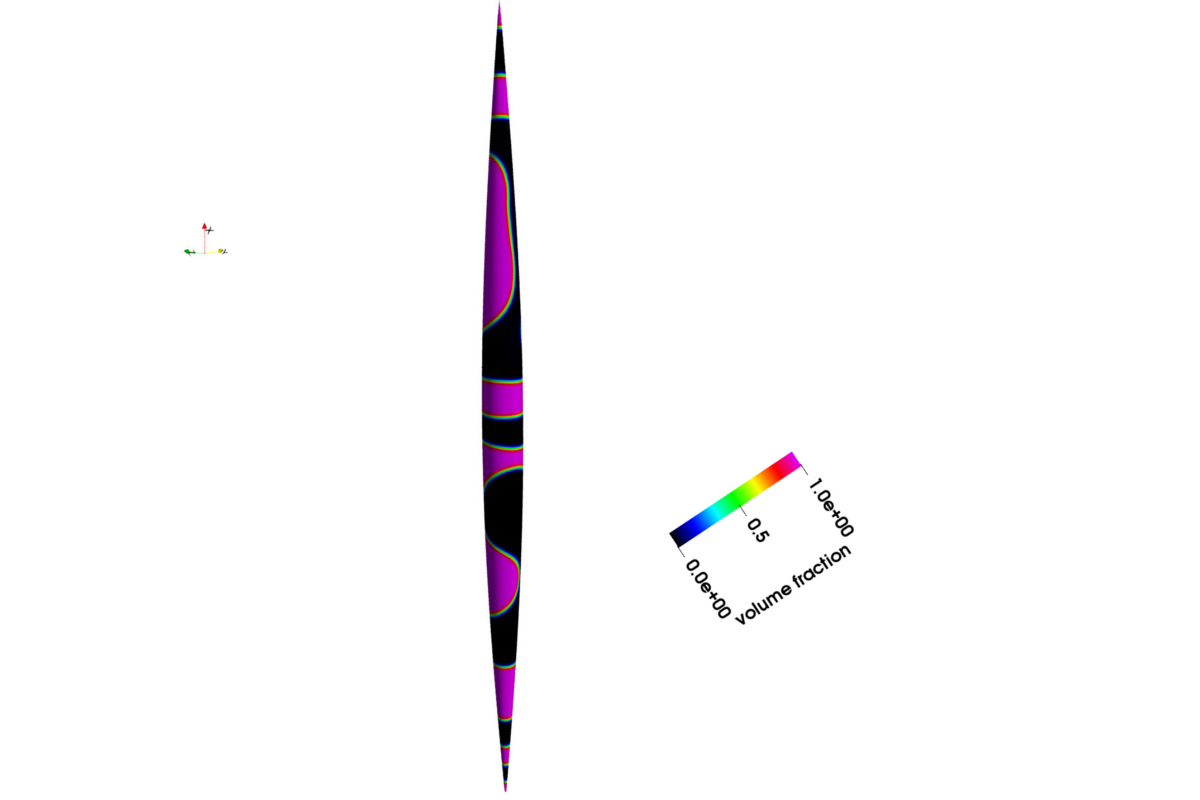}
         \put(3,103){\small{$t = 400$}}
\end{overpic}
\hspace{.5cm}
\begin{overpic}[width=.10\textwidth,grid=false]{figures_AC_sphere_legend.png}
\end{overpic}
\end{center}
\caption{\label{spindleCH_0_4}
Evolution of the numerical solution of the Cahn--Hilliard equation on the spindle for $t \in (0, 400]$
computed on a mesh with size $h =0.026$. Time step is
$\Delta t =0.01$ for $t \in (0, 1]$ and $\Delta t =1$ for $t \in (1, 400]$.
The simulation was started with a random initial condition depicted in the leftmost panel.
View: $xy$-plane.}
\end{figure}

The evolution of the numerical solution to the Cahn--Hilliard equation
for $t \in (0, 400]$ is shown in Fig.~\ref{spindleCH_0_4}. Just like on the sphere,
a pattern emerges as early as $t = 0.5$. In order to follow this fast initial phase, we
choose the same time steps as for the sphere, i.e. $\Delta t =0.01$ for $t \in [0, 1]$ and
$\Delta t = 1$ for $t > 1$.

In order to compare phase separation on the spindle modeled by the Allen--Cahn and Cahn--Hilliard equations,
we point out that Fig.~\ref{spindleAC_0_4}
and \ref{spindleCH_0_4} display solutions computed at different times. Although we see a
faster and more ordered separation of phases from the Cahn--Hilliard model, the solutions computed
with both models at time $t = 400$ are not significantly different in terms of surface area occupied
by each phase. Compare the rightmost panels in Fig.~\ref{spindleAC_0_4}
and \ref{spindleCH_0_4}. So, from these preliminary results it is not possible to conclude which
model is better suited to simulate the formation of microdomains in the cytoplasmic membrane of bacillus subtilis.
Further investigation is needed. In addition, the role played by the initial condition
that has to be understood, since in our tests on the spindle we have not tried
anything other than random initial condition.

\subsection{Phase separation on  an idealized cell}\label{sec:scell}

\begin{figure}
  \centering
    \subfloat[Side view]{\includegraphics[width=.5\textwidth]{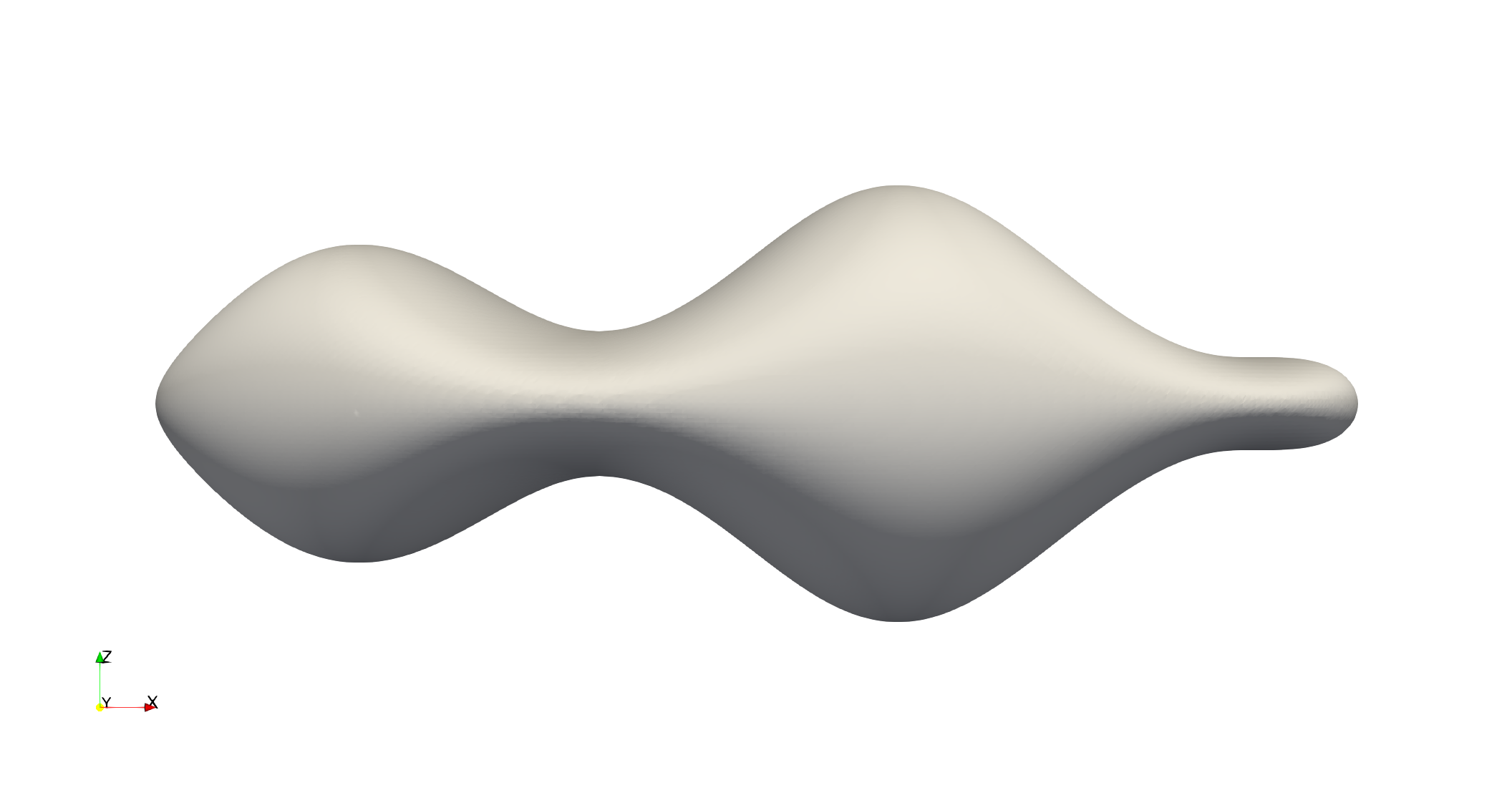}}~
  \subfloat[Angled view with surface mesh]{\includegraphics[width=.5\textwidth]{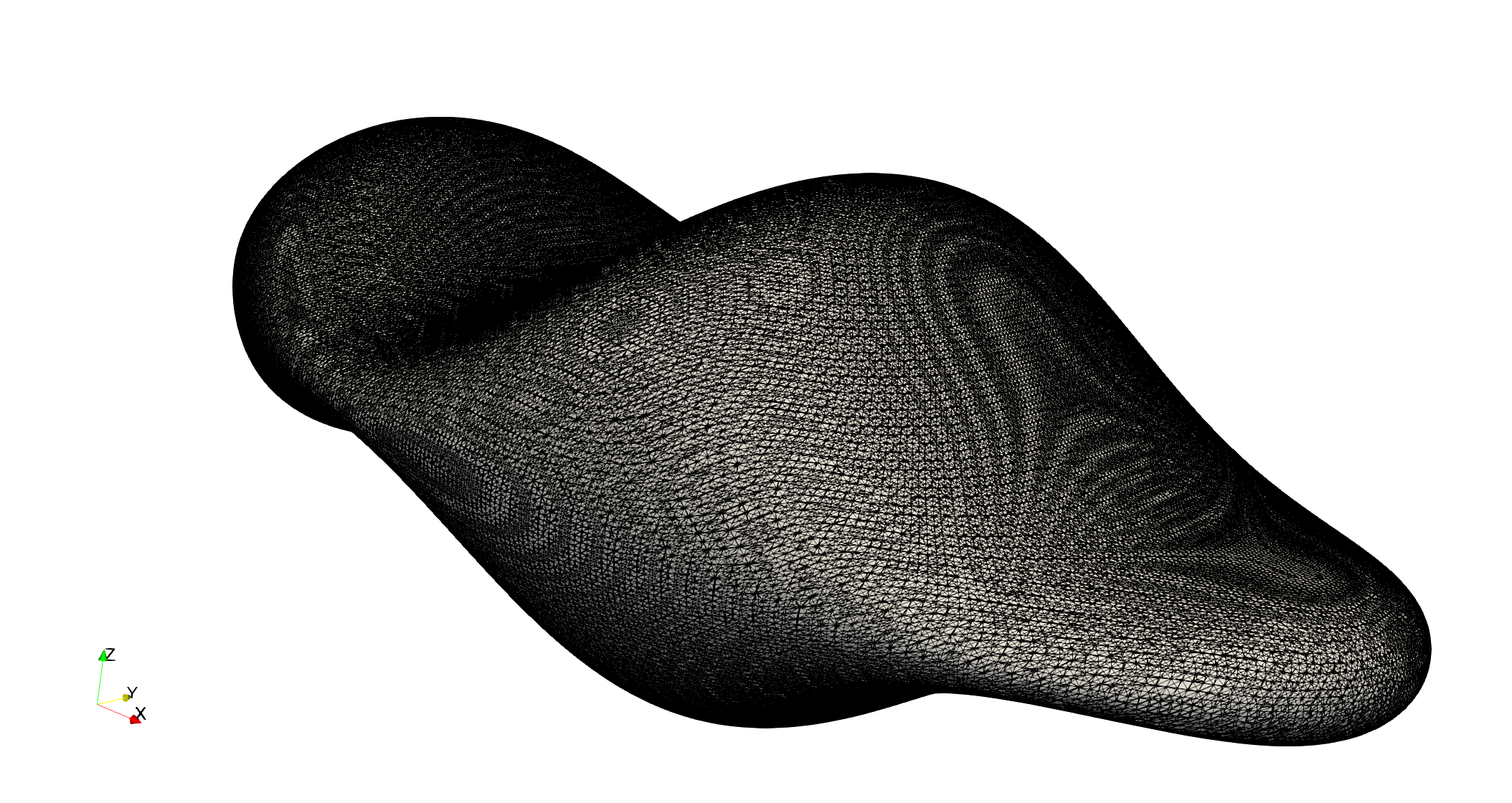}}
  \hfill
  \caption{Illustration of the surface of the idealized cell: (a) side view and (b) angled view with the trace surface for the
  mesh in use.}
  \label{fig:amorphous}
\end{figure}

In this example, the surface $\Gamma$  is given by the zero of the following level set function taken from \cite{dziuk2013finite}:
\begin{align}
\phi(\bx) = \frac{1}{4} x_1^2 + x_2^2 + \frac{4 x_3^2}{(1 + \frac{1}{2} \sin(\pi x_1))^2}-1. \el
\end{align}
The surface is illustrated in Fig.~\ref{fig:amorphous} and can be viewed as an
idealized cell.
As mentioned in Sec.~\ref{sec:num_method}, our long term goal is to
simulate phase separation on biological membranes that exhibit shape transitions and shape instabilities.
Thus, we need to go beyond simple surfaces like the sphere and the spindle.
Besides its more complex, amorphous shape, what distinguishes the surface of the idealized cell
is that it is possible to trace a curve of minimal length on it.
We shall see that this geometric property critically defines the equilibrium stage of the non-conservative evolution of lateral  phase separation process.

\begin{figure}
  \centering
    \subfloat[$\times$ 4 zoom]{\includegraphics[width=.31\textwidth]{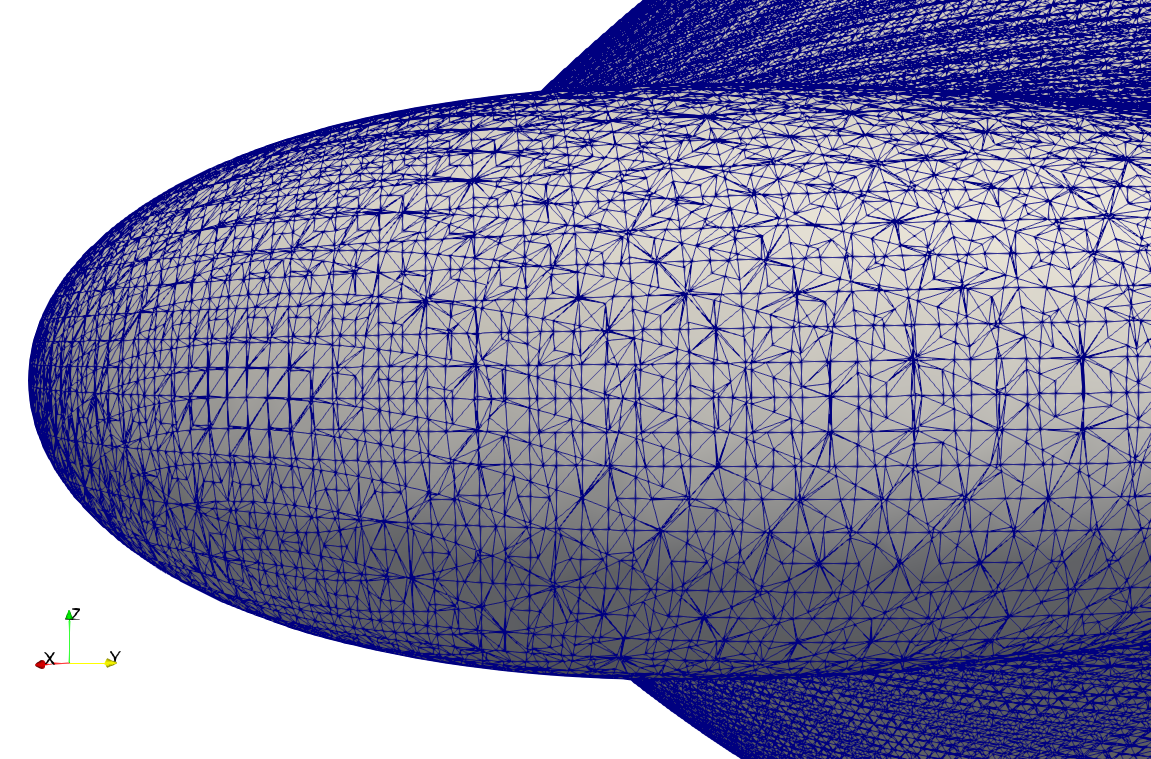}}~~
  \subfloat[$\times$ 8 zoom]{\includegraphics[width=.31\textwidth]{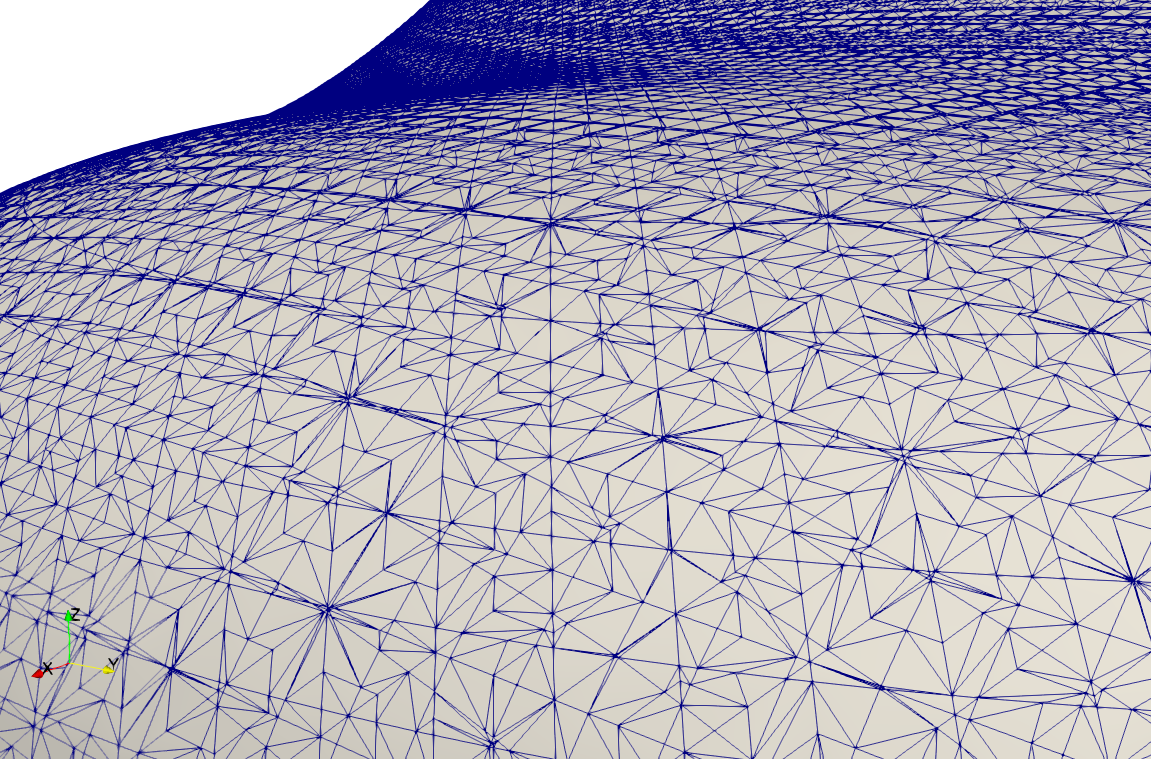}}~~
  \subfloat[$\times$ 24 zoom]{\includegraphics[width=.31\textwidth]{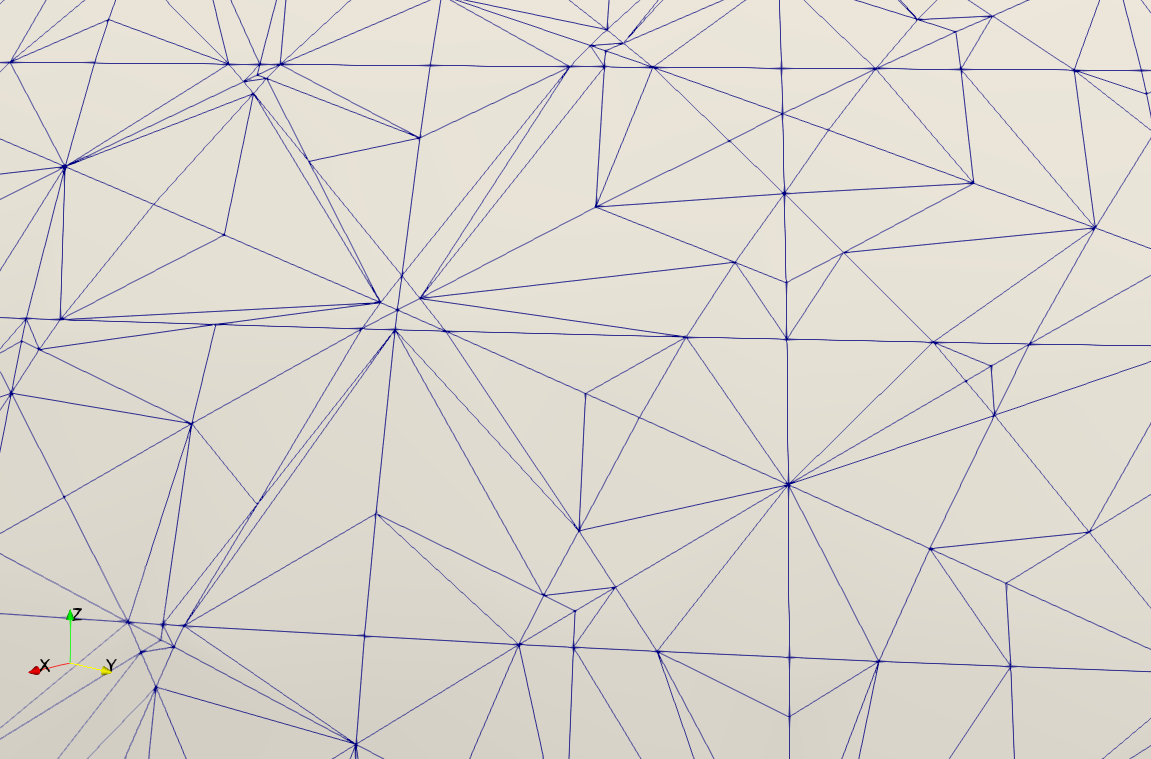}}
  \hfill
  \caption{Zoomed views the surface triangulation induced by embedding of the surface into the bulk tetrahedral mesh.}
  \label{fig:mesh}
\end{figure}

The surface is embedded in an outer domain $\Omega= [-2,2] \times [-4/3, 4/3] \times [-4/3, 4/3]$.
The computational mesh $\T_{h}^\Gamma$ of $\Omega$ has mesh size $h =0.031$. This
tetrahedral mesh is generated in the same way as
the meshes for the sphere and the spindle are generated, i.e. by diving $\Omega$
into cubes and then diving the cubes into tetrahedra. The elements cut by the surface
are further refined.
This mesh has a level of refinement comparable to mesh $\ell = 6$ in Sec.~\ref{sec:NumVal}.
We recall that the discrete surface $\Gamma_h$ is implicitly
defined by the zero set of $\phi_{h/2}=I_{h/2}(\phi)$ and does not require
any explicit parametrization of $\Gamma$. This makes our numerical
approach particularly suitable for dealing with complex shapes. The resulting surface mesh, which is
illustrated in Fig.~\ref{fig:amorphous} (b) for the entire $\Gamma$ and
in Fig.~\ref{fig:mesh} for the `beak' part of $\Gamma$,
is used for numerical integration only, while test and trial finite element functions are defined on shape regular bulk tetrahedra.

We simulate phase separation on the idealized cell surface
using the Allen--Cahn and Cahn-Hillard models, i.e.~eq.~\eqref{eq:AC_eq} and eq.~\eqref{eq:CH}, with free energy per unit surface \eqref{eq:f0}, $\alpha = 1$, and $\epsilon = 0.01$.
Like for the sphere and the spindle, the initial condition is $\eta_0(\bx)=\text{rand}(\bx)$.

\subsubsection{Allen--Cahn model}\label{sec:cell_AC}

\begin{figure}
\begin{center}
\begin{overpic}[width=.20\textwidth,viewport=85 70 880 720, clip,grid=false]{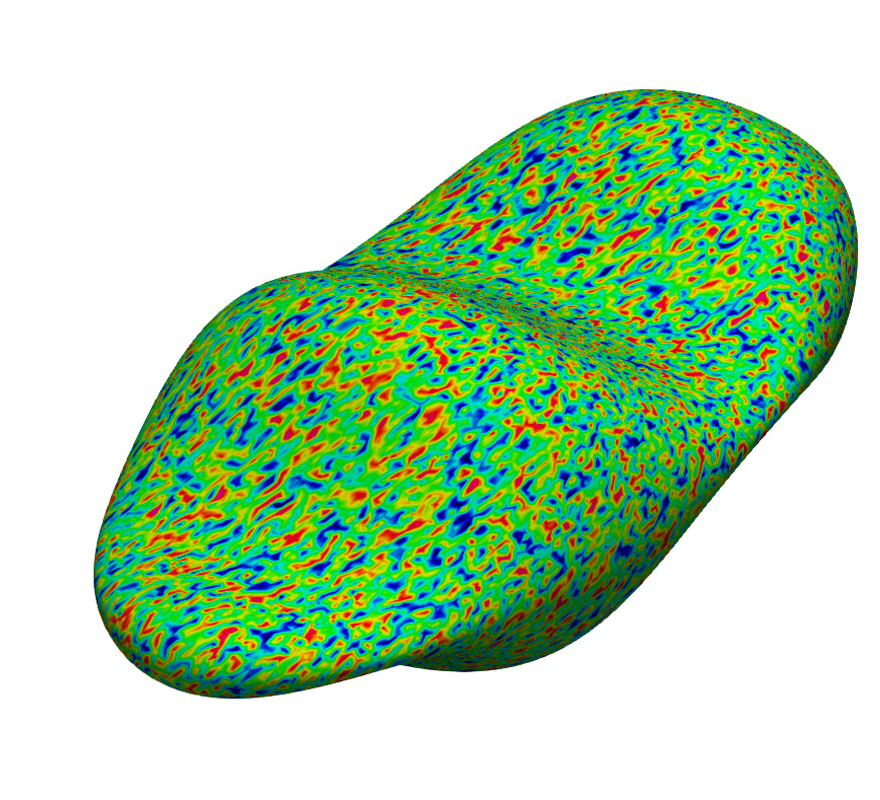}
        \put(35,85){\small{$t = 0$}}
\end{overpic}
\begin{overpic}[width=.20\textwidth, viewport=85 70 880 720, clip,grid=false]{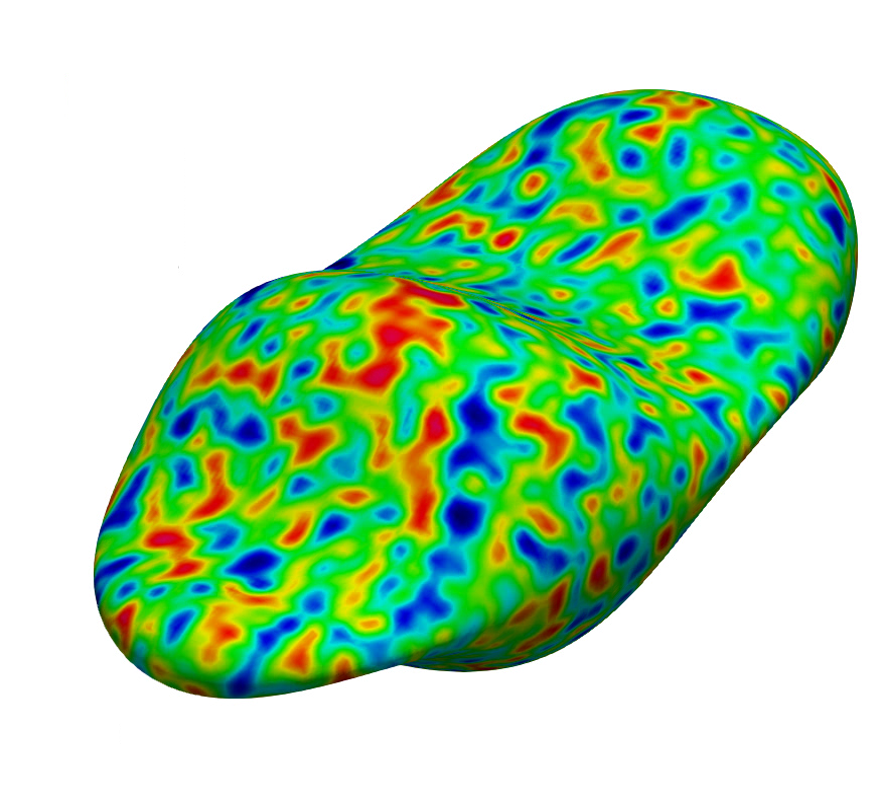}
        \put(35,85){\small{$t = 5$}}
\end{overpic}
\begin{overpic}[width=.20\textwidth, viewport=85 70 880 720, clip,grid=false]{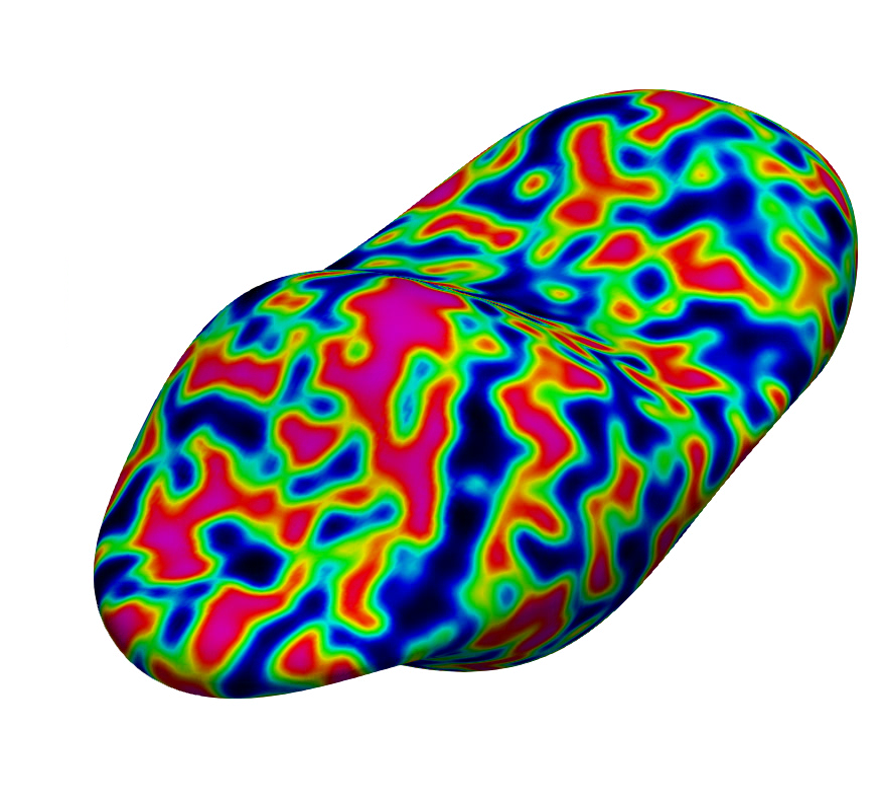}
        \put(35,85){\small{$t = 10$}}
\end{overpic}
\begin{overpic}[width=.20\textwidth, viewport=85 70 880 720, clip,grid=false]{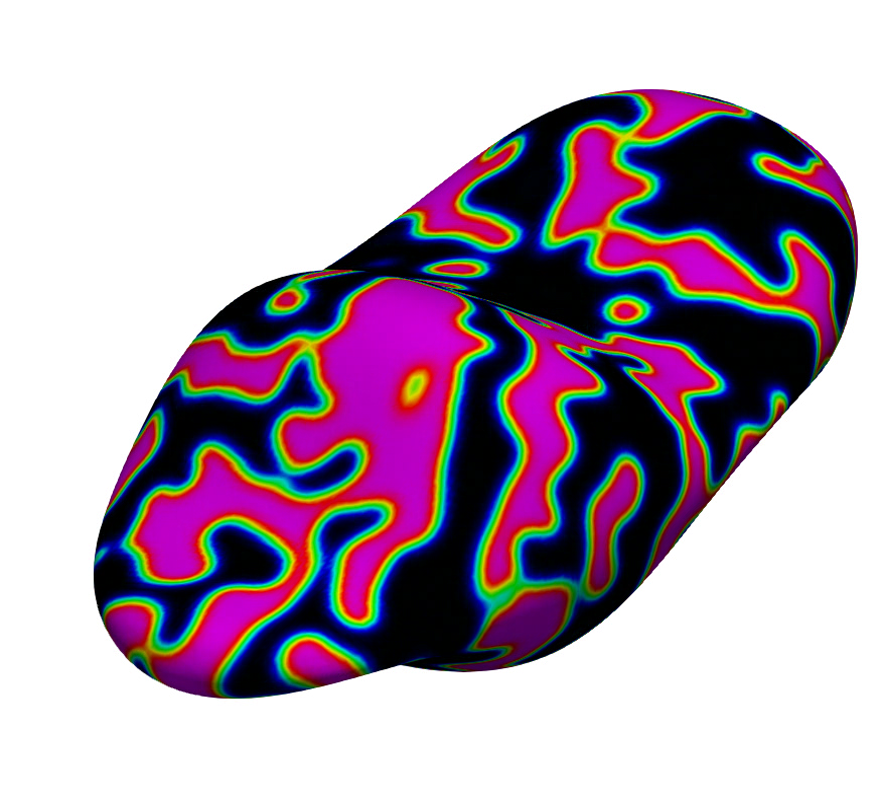}
         \put(33,85){\small{$t = 25$}}
\end{overpic}
\vskip .3cm
\begin{overpic}[width=.20\textwidth, viewport=85 70 880 720, clip,grid=false]{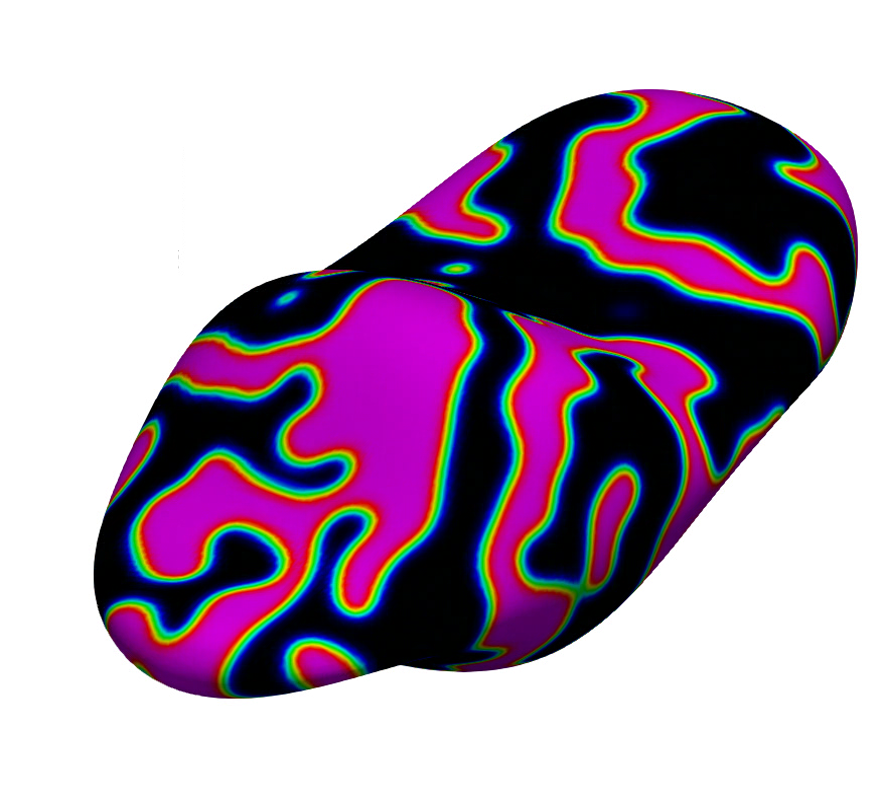}
        \put(33,85){\small{$t = 50$}}
\end{overpic}
\begin{overpic}[width=.20\textwidth, viewport=85 70 880 720, clip,grid=false]{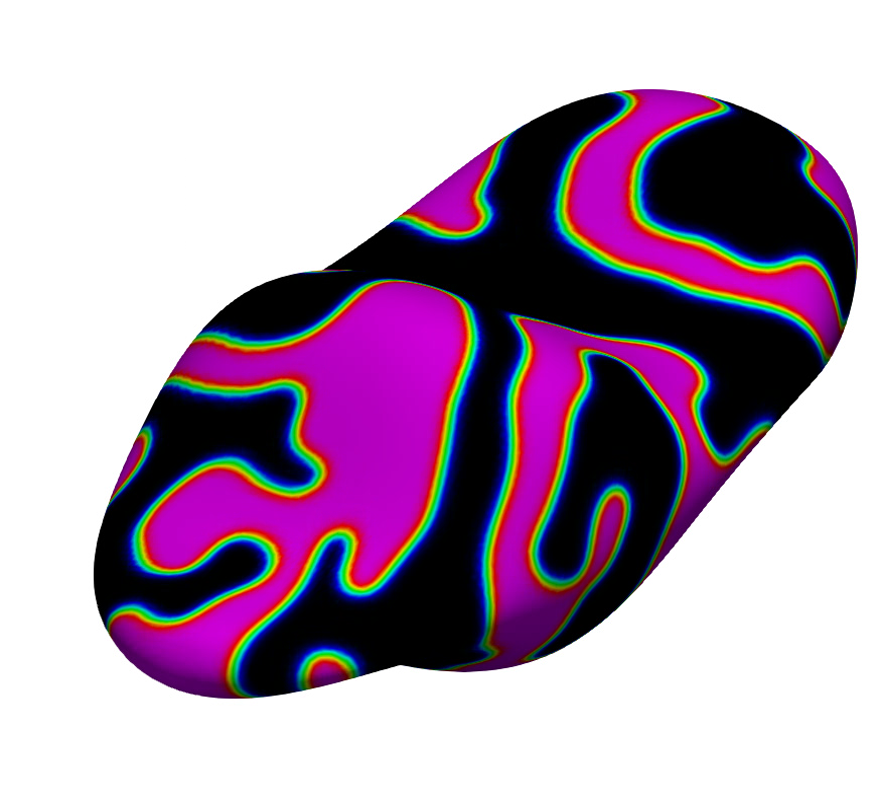}
        \put(33,85){\small{$t = 100$}}
\end{overpic}
\begin{overpic}[width=.20\textwidth, viewport=85 70 880 720, clip,grid=false]{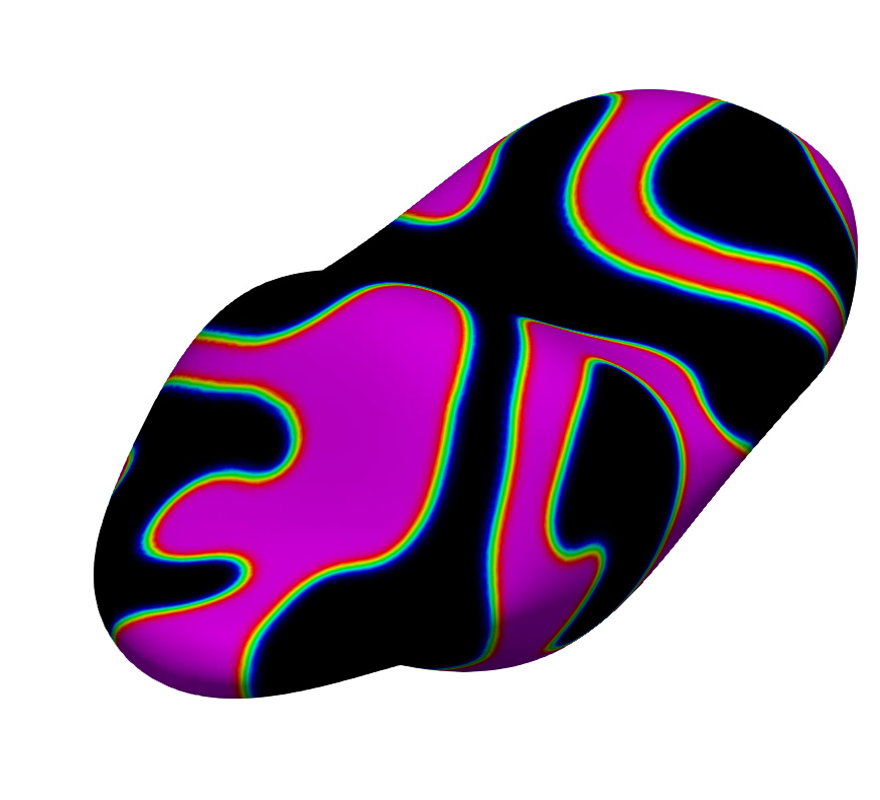}
        \put(33,85){\small{$t = 200$}}
\end{overpic}
\hspace{1.2cm}
\begin{overpic}[width=.08\textwidth,grid=false]{figures_AC_sphere_legend.png}
\end{overpic}
\end{center}
\caption{\label{cell_0_2}
Evolution of the numerical solution of the Allen--Cahn equation on the idealized cell surface for $t \in (0, 200]$
computed with time step as specified in Table \ref{tab:dt_cell} and on a mesh with size $h =0.031$.
The simulation was started with a random initial condition depicted in the top left panel.}
\end{figure}

\begin{figure}
\begin{center}
\begin{overpic}[width=.20\textwidth,viewport=90 190 680 600, clip,grid=false]{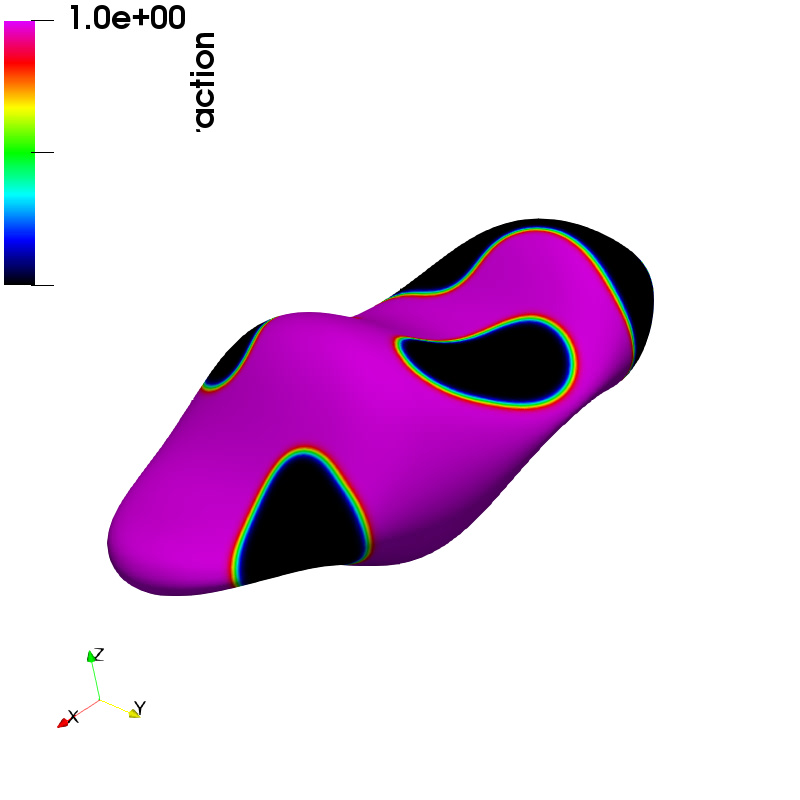}
        \put(33,70){\small{$t = 500$}}
\end{overpic}
\begin{overpic}[width=.20\textwidth, viewport=90 190 680 600, clip,grid=false]{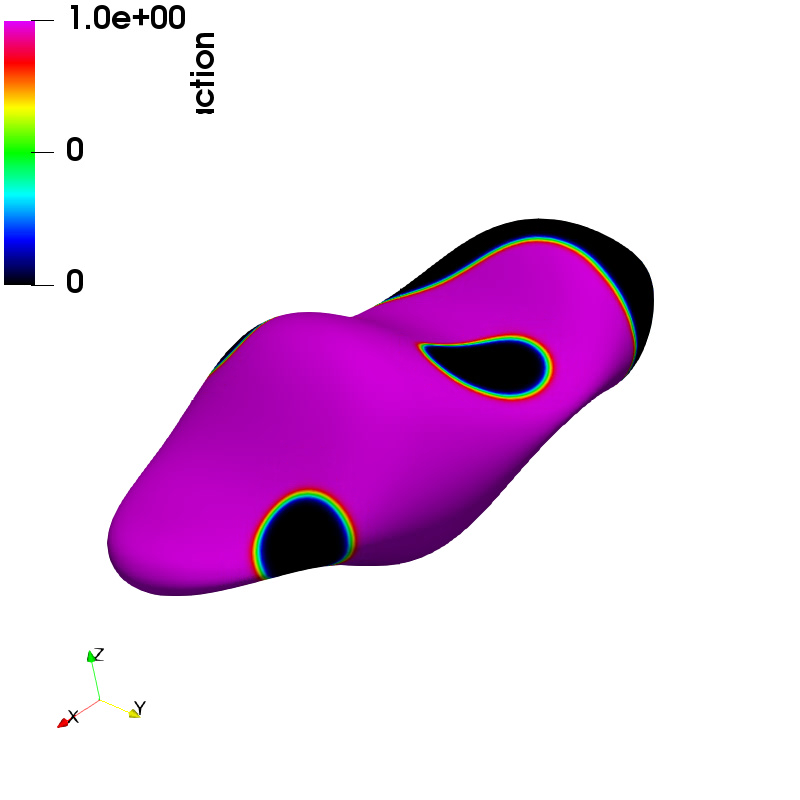}
        \put(33,70){\small{$t = 1000$}}
\end{overpic}
\begin{overpic}[width=.20\textwidth, viewport=90 190 680 600, clip,grid=false]{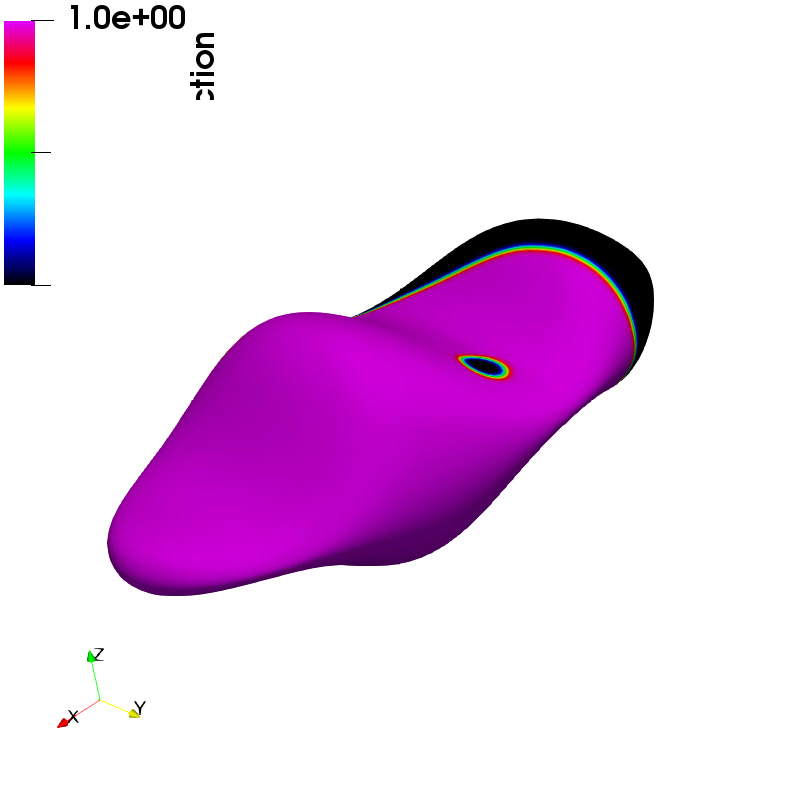}
        \put(33,70){\small{$t = 1500$}}
\end{overpic}
\begin{overpic}[width=.20\textwidth, viewport=90 190 680 600, clip,grid=false]{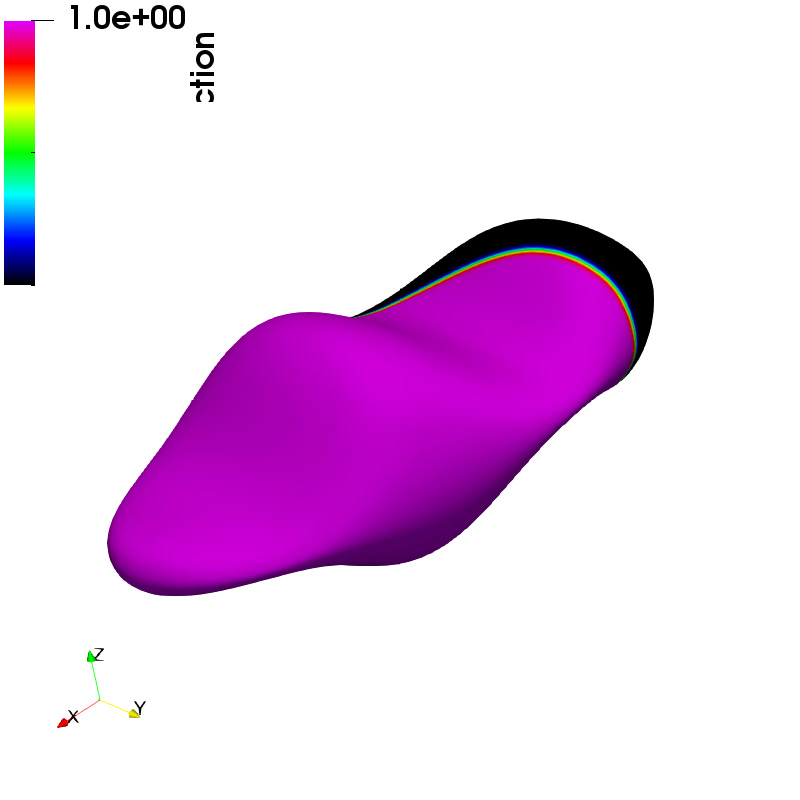}
        \put(33,70){\small{$t = 2000$}}
\end{overpic}
\begin{overpic}[width=.08\textwidth,grid=false]{figures_AC_sphere_legend.png}
\end{overpic}
\end{center}
\caption{\label{cell_5_20}
Evolution of the numerical solution of the Allen--Cahn equation on the idealized cell surface for $t \in [500, 2000]$
computed with time step as specified in Table \ref{tab:dt_cell} and
on a mesh with size $h =0.031$.
}
\end{figure}

\begin{figure}
\begin{center}
\begin{overpic}[width=.20\textwidth,viewport=90 190 680 600, clip,grid=false]{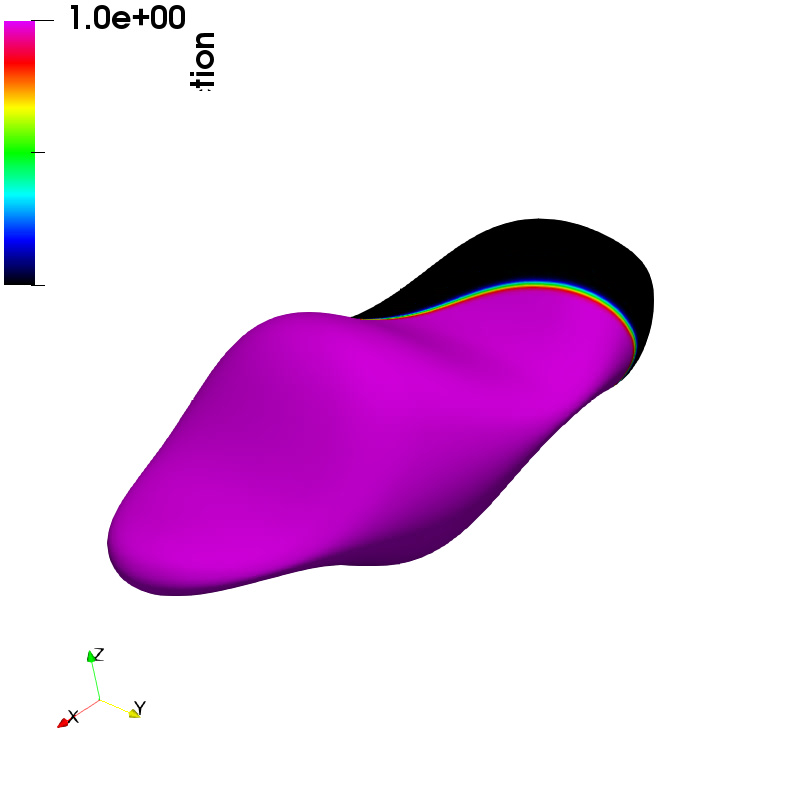}
        \put(33,70){\small{$t = 4000$}}
\end{overpic}
\begin{overpic}[width=.20\textwidth, viewport=90 190 680 600, clip,grid=false]{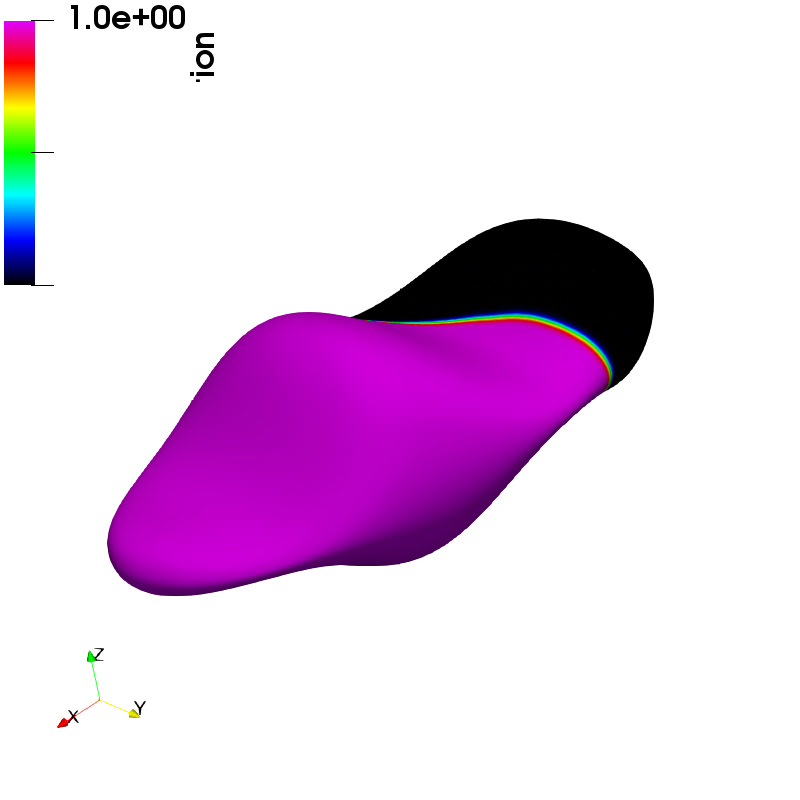}
        \put(33,70){\small{$t = 8000$}}
\end{overpic}
\begin{overpic}[width=.20\textwidth, viewport=90 190 680 600, clip,grid=false]{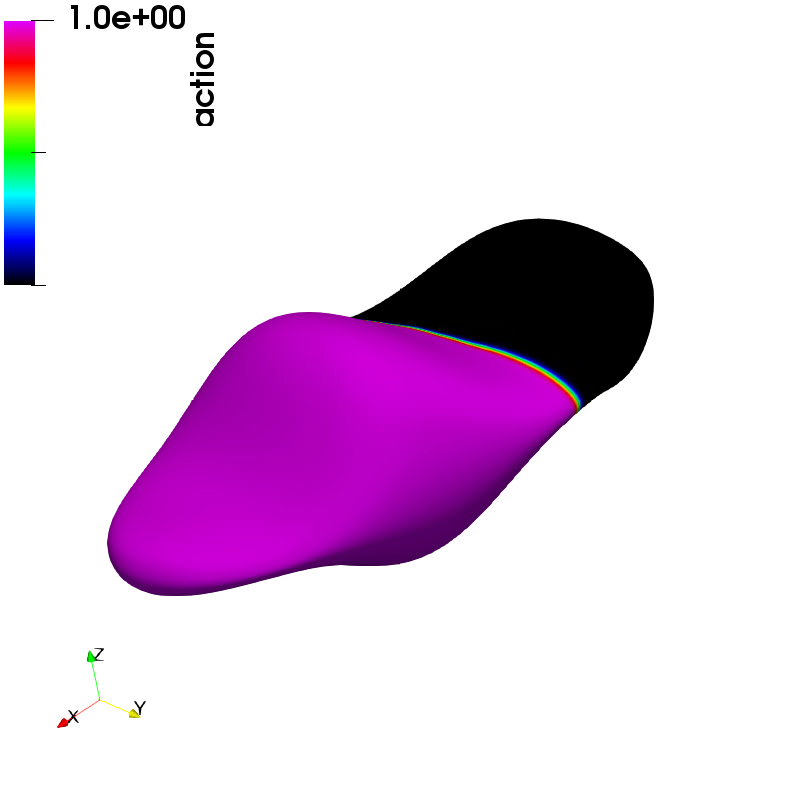}
        \put(33,70){\small{$t = 16000$}}
\end{overpic}
\begin{overpic}[width=.20\textwidth, viewport=90 190 680 600, clip,grid=false]{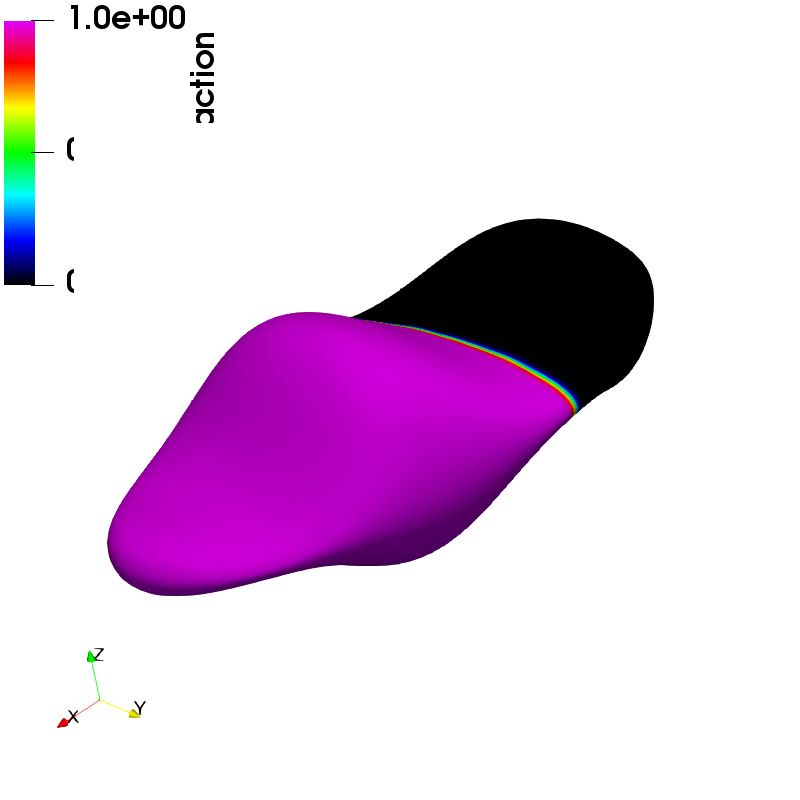}
        \put(33,70){\small{$t = 20000$}}
\end{overpic}
\begin{overpic}[width=.08\textwidth,grid=false]{figures_AC_sphere_legend.png}
\end{overpic}
\end{center}
\caption{\label{cell_40_200}
Evolution of the numerical solution of the Allen--Cahn equation on the idealized cell surface for $t \in [4000, 20000]$
computed with time step as specified in Table \ref{tab:dt_cell} and
on a mesh with size $h =0.031$.
}
\end{figure}

First we apply Allen--Cahn model to simulate non-conservative evolution of phase separation.
Fig.~\ref{cell_0_2} shows the initial fast evolution of the numerical solution of the Allen--Cahn
equation, followed by the beginning of the slowdown phase.
In Fig.~\ref{cell_5_20} we see the process of separation into two regions, one pink region
with $\eta = 1$ and one black region with $\eta = 0$. The separation itself occurs around $t=1600$.
Finally, Fig.~\ref{cell_40_200} shows the evolution towards the steady state.
After $t = 20000$ there is no visible change in the position of the interface
between phases, which virtually \emph{coincided with the curve of minimal length} on $\Gamma$.
The computed equilibrium state is consistent with the well-known limiting (as $\epsilon\to0$) behaviour of
the Allen--Cahn model, which defines the evolution of the sharp interface as the mean curvature motion~\cite{evans1992phase}. For the surface Allen--Cahn equation one obtains geodesic curvature motion of the interface in the asymptotic limit~\cite{elliott2010modeling} with any closed geodesic as an equilibrium state.
We note that for the given surface this equilibrium state is stable.
Therefore, unlike the cases of the sphere and spindle, on the idealized cell neither phase vanishes.

Like for the sphere, we prescribe different time steps for the different stages
of the phase separation. Table \ref{tab:dt_sphere} reports the time step assigned for each time intervals.

\begin{table}
\begin{center}
  \begin{tabular}{ | c | c | c |  c |}
    \hline
    interval & $(0, 200]$ & $(200, 500]$ & $(500, 23000]$   \\
    \hline
    $\Delta t$ & 1 & 5 & 10 \\
    \hline
  \end{tabular}
  \caption{\label{tab:dt_cell}
Time steps used for the different time intervals to obtain the results Fig.~\ref{cell_0_2}, \ref{cell_5_20}, and \ref{cell_40_200}.}
\end{center}
\end{table}


\subsubsection{Cahn--Hilliard model} \label{sec:CH_cell}

\begin{figure}
\begin{center}
\begin{overpic}[width=.20\textwidth,viewport=110 150 700 580, clip,grid=false]{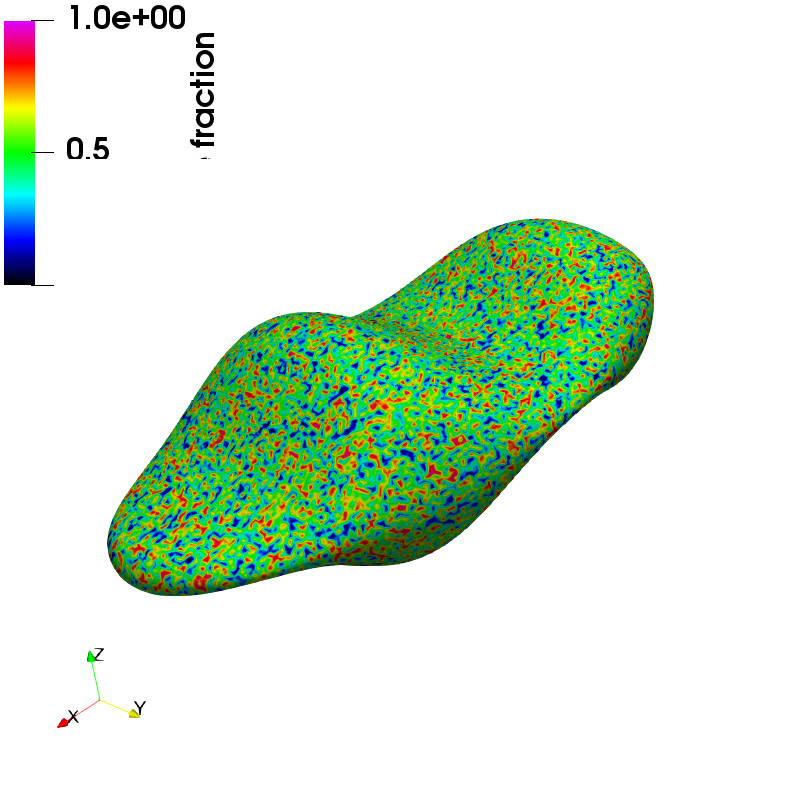}
        \put(35,80){\small{$t = 0$}}
\end{overpic}
\begin{overpic}[width=.20\textwidth,viewport=110 150 700 580, clip,grid=false]{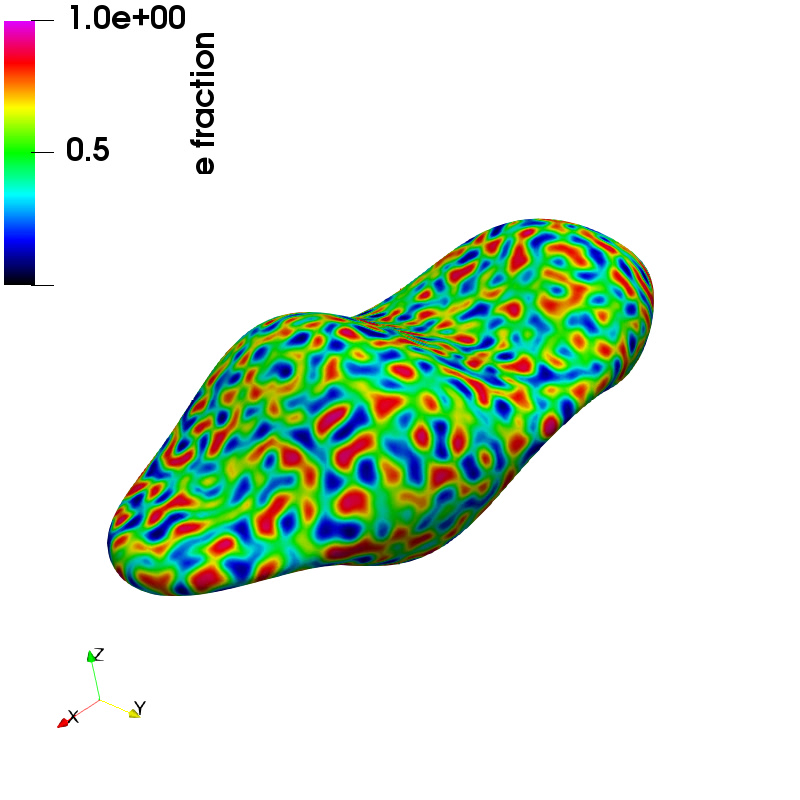}
        \put(35,80){\small{$t = 0.1$}}
\end{overpic}
\begin{overpic}[width=.20\textwidth,viewport=110 150 700 580, clip,grid=false]{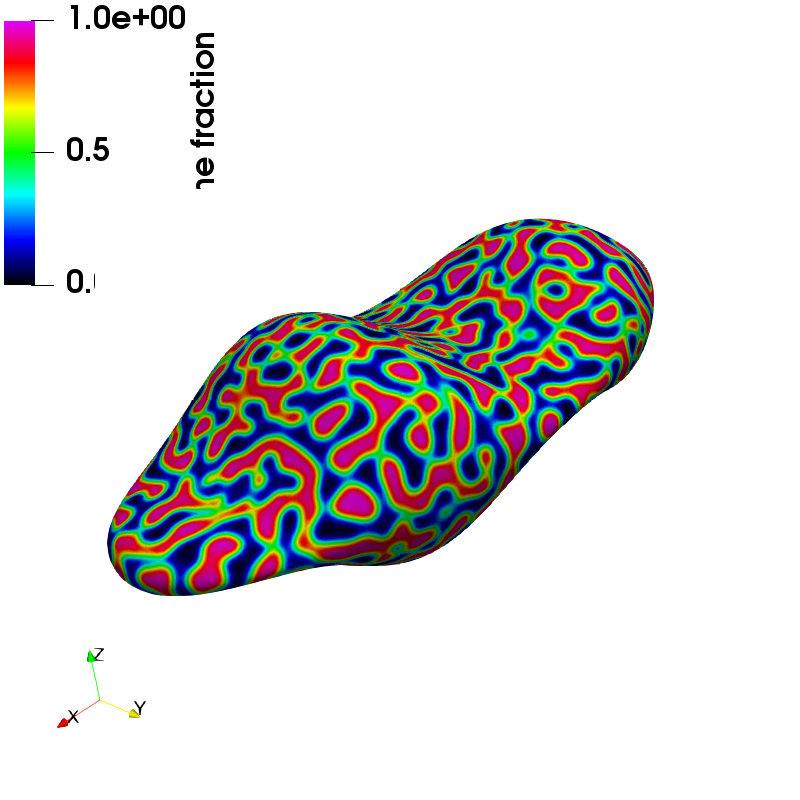}
        \put(35,80){\small{$t = 0.2$}}
\end{overpic}
\begin{overpic}[width=.20\textwidth,viewport=110 150 700 580, clip,grid=false]{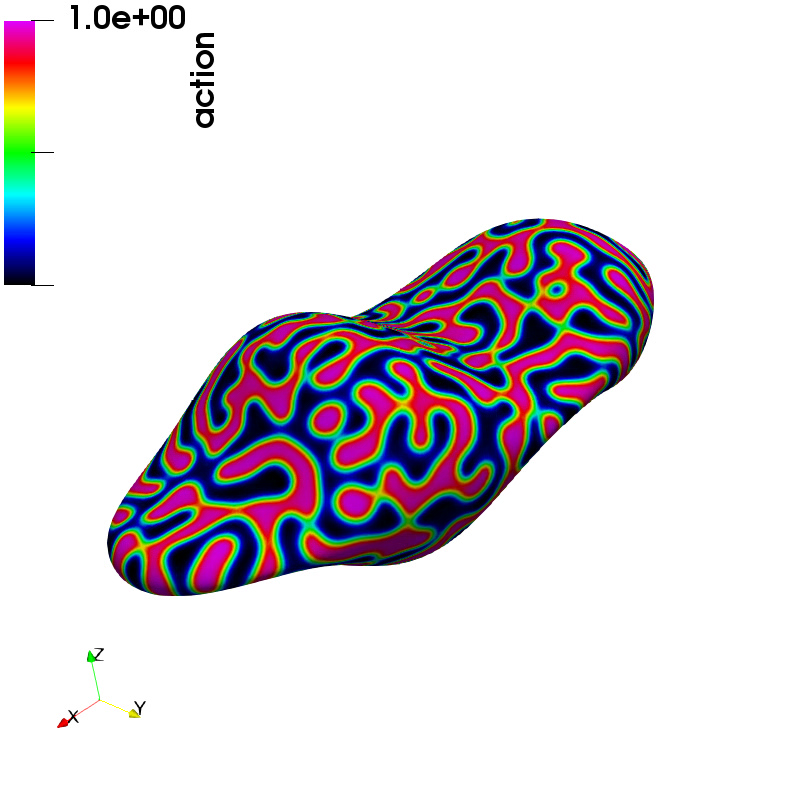}
        \put(35,80){\small{$t = 0.5$}}
\end{overpic}
\vskip .3cm
\begin{overpic}[width=.20\textwidth,viewport=110 150 700 580, clip,grid=false]{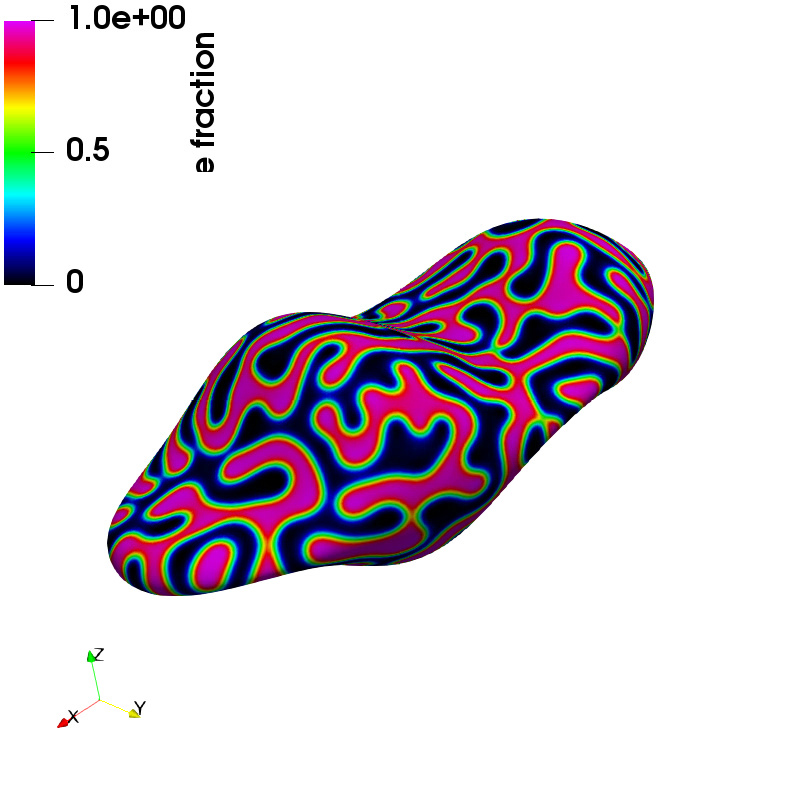}
        \put(35,80){\small{$t = 1$}}
\end{overpic}
\begin{overpic}[width=.20\textwidth,viewport=110 150 700 580, clip,grid=false]{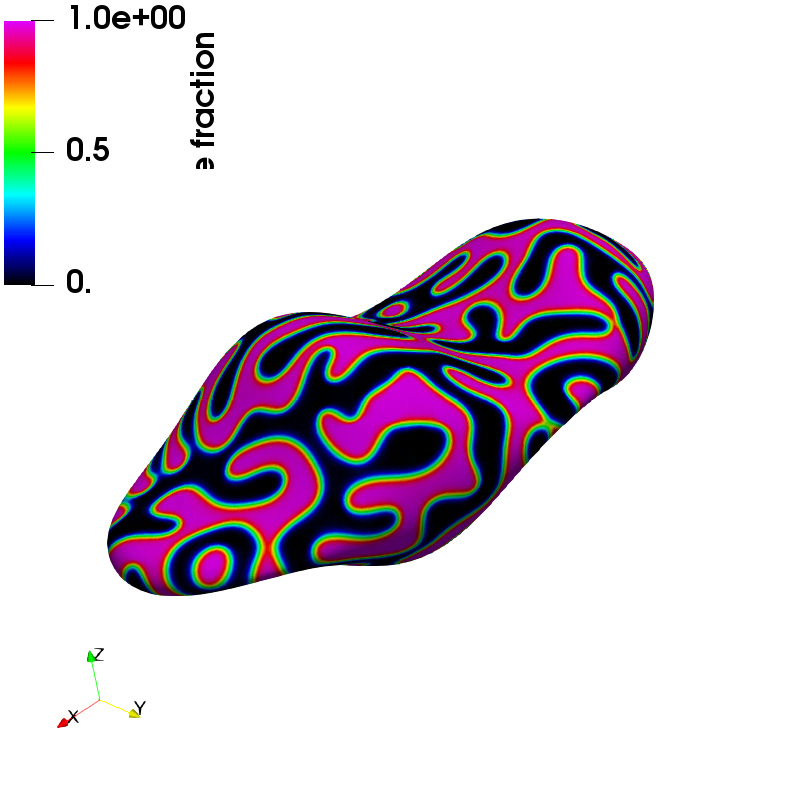}
        \put(35,80){\small{$t = 10$}}
\end{overpic}
\begin{overpic}[width=.20\textwidth,viewport=110 150 700 580, clip,grid=false]{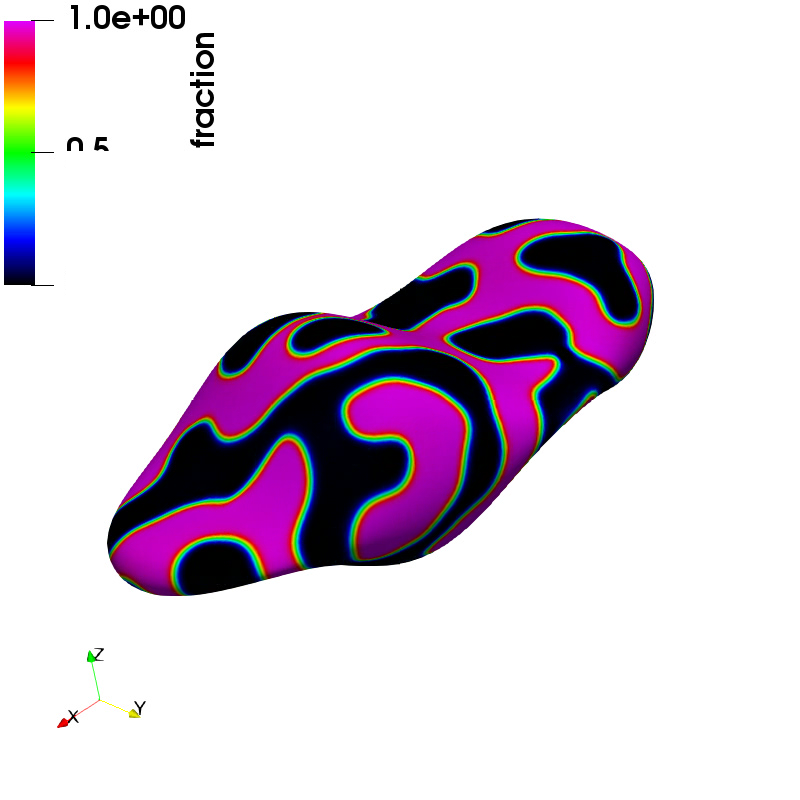}
        \put(35,80){\small{$t = 50$}}
\end{overpic}
\begin{overpic}[width=.20\textwidth,viewport=110 150 700 580, clip,grid=false]{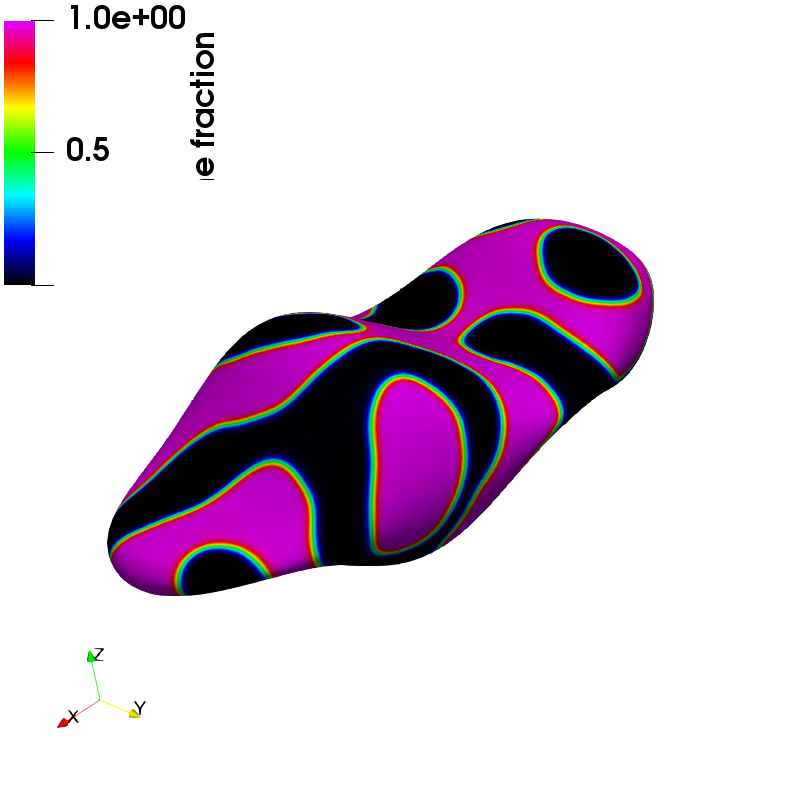}
        \put(35,80){\small{$t = 100$}}
\end{overpic}
\vskip .3cm
\begin{overpic}[width=.20\textwidth,viewport=110 150 700 580, clip,grid=false]{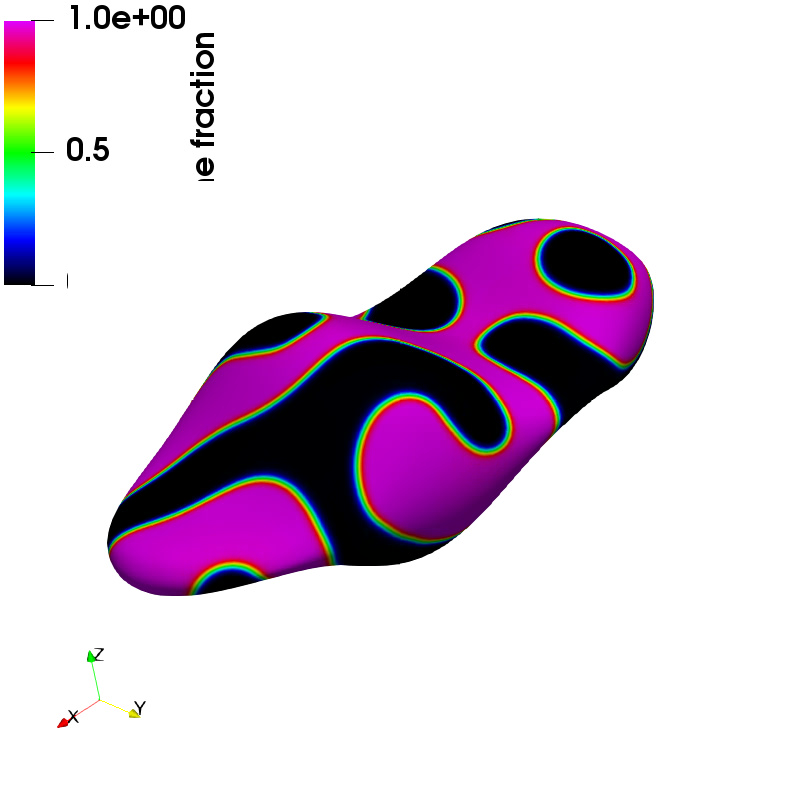}
        \put(35,80){\small{$t = 200$}}
\end{overpic}
\begin{overpic}[width=.20\textwidth,viewport=110 150 700 580, clip,grid=false]{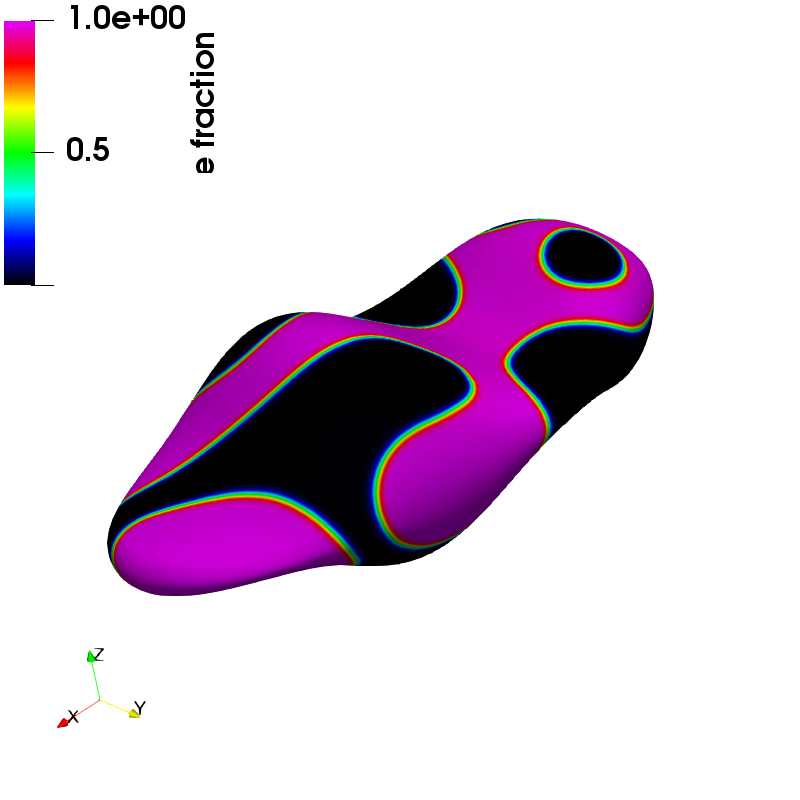}
        \put(35,80){\small{$t = 500$}}
\end{overpic}
\begin{overpic}[width=.20\textwidth,viewport=110 150 700 580, clip,grid=false]{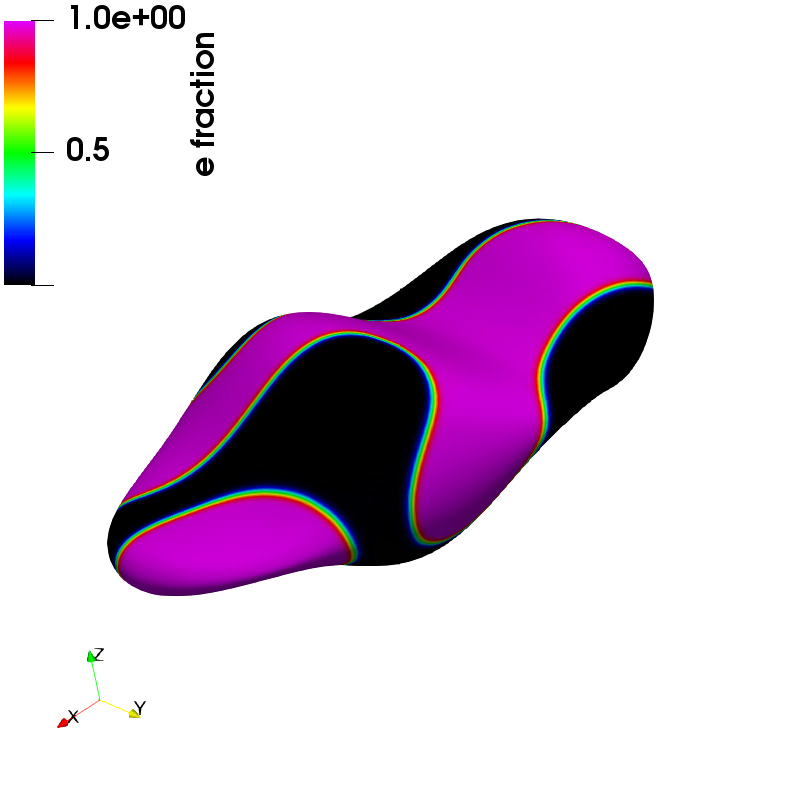}
        \put(35,80){\small{$t = 1000$}}
\end{overpic}
\hspace{1.2cm}
\begin{overpic}[width=.08\textwidth,grid=false]{figures_AC_sphere_legend.png}
\end{overpic}
\end{center}
\caption{\label{cellCH_0_10}
Evolution of the numerical solution of the Cahn--Hilliard equation on the idealized cell surface for $t \in (0, 1000]$
computed on a mesh with size $h =0.031$.
Time step is $\Delta t =0.01$ for $t \in (0, 1]$ and $\Delta t =1$ for $t \in (1, 1000]$.
The simulation was started with a random initial condition depicted in the top left panel.}
\end{figure}

In the last subsection, we apply the Cahn--Hilliard
model to simulate the conservative evolution of phase separation on the idealized cell.
The evolution of the numerical solution to the Cahn--Hilliard equation
for $t \in [0, 1000]$ is shown in Fig.~\ref{cellCH_0_10}. Just like on the sphere and on the spindle,
we see the initial fast phase that ends around $t = 1$ and is followed by the slowdown phase.
Thus, the time steps are chosen like for the tests on the sphere and on the spindle.

\begin{figure}
\begin{center}
\begin{overpic}[width=.18\textwidth,viewport=110 150 700 580, clip,grid=false]{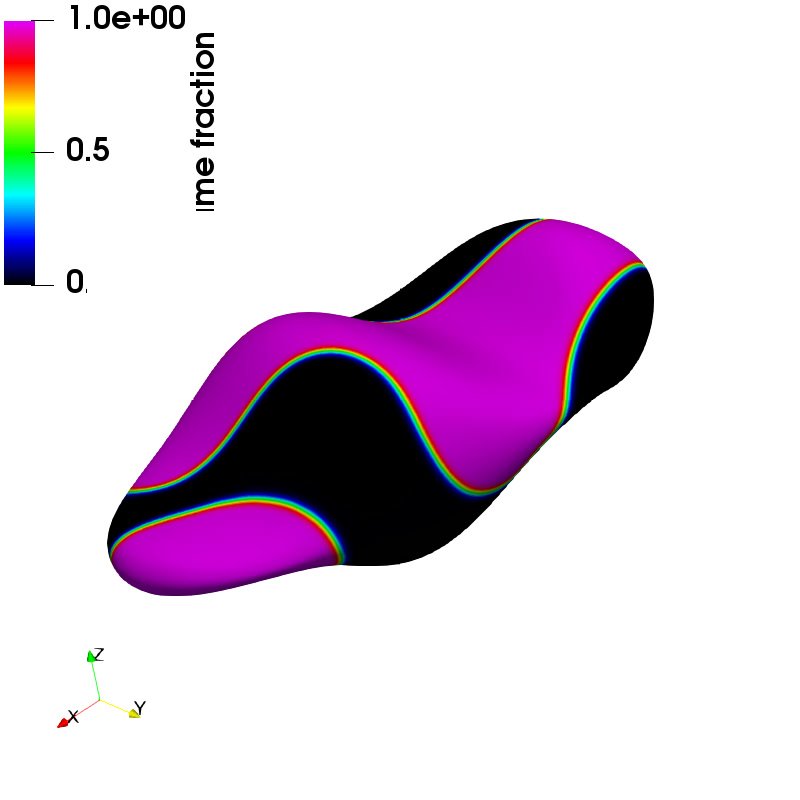}
        \put(35,80){\small{$t = 2000$}}
\end{overpic}
\begin{overpic}[width=.18\textwidth,viewport=110 150 700 580, clip,grid=false]{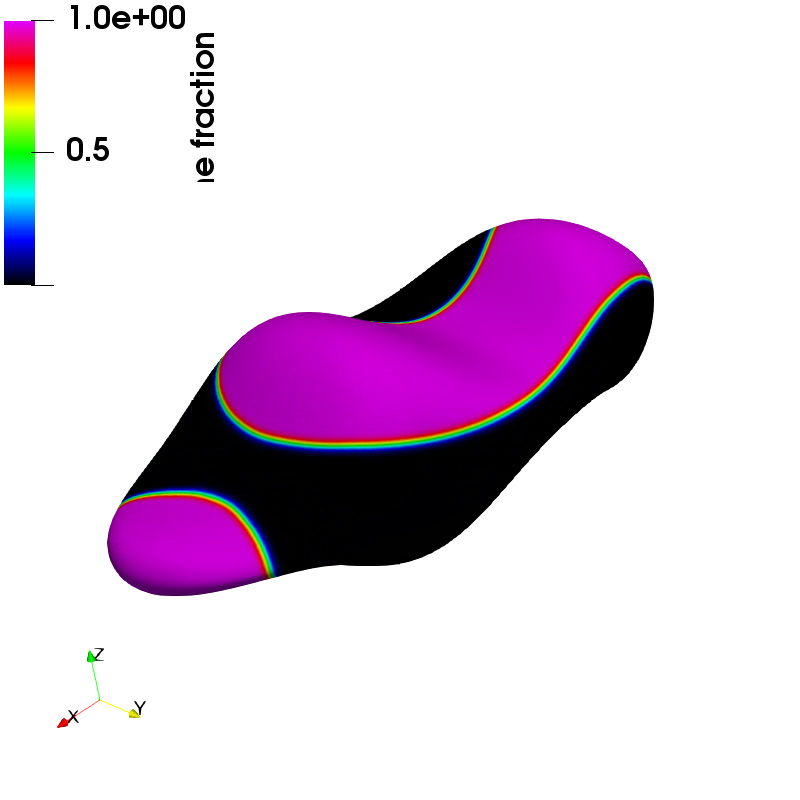}
        \put(35,80){\small{$t = 8000$}}
\end{overpic}
\begin{overpic}[width=.18\textwidth,viewport=110 150 700 580, clip,grid=false]{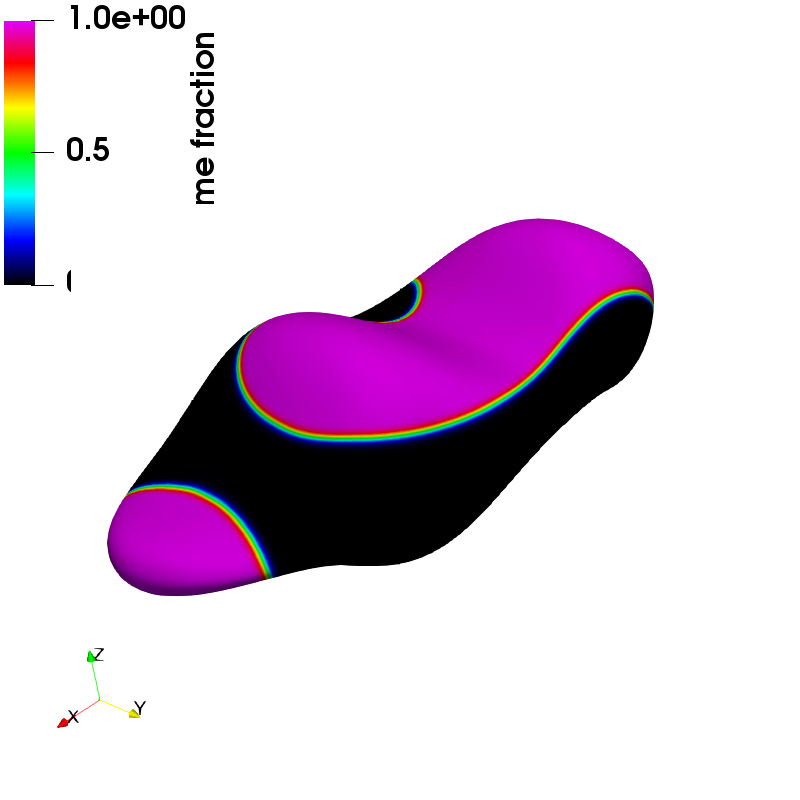}
        \put(35,80){\small{$t = 15000$}}
\end{overpic}
\begin{overpic}[width=.18\textwidth,viewport=110 150 700 580, clip,grid=false]{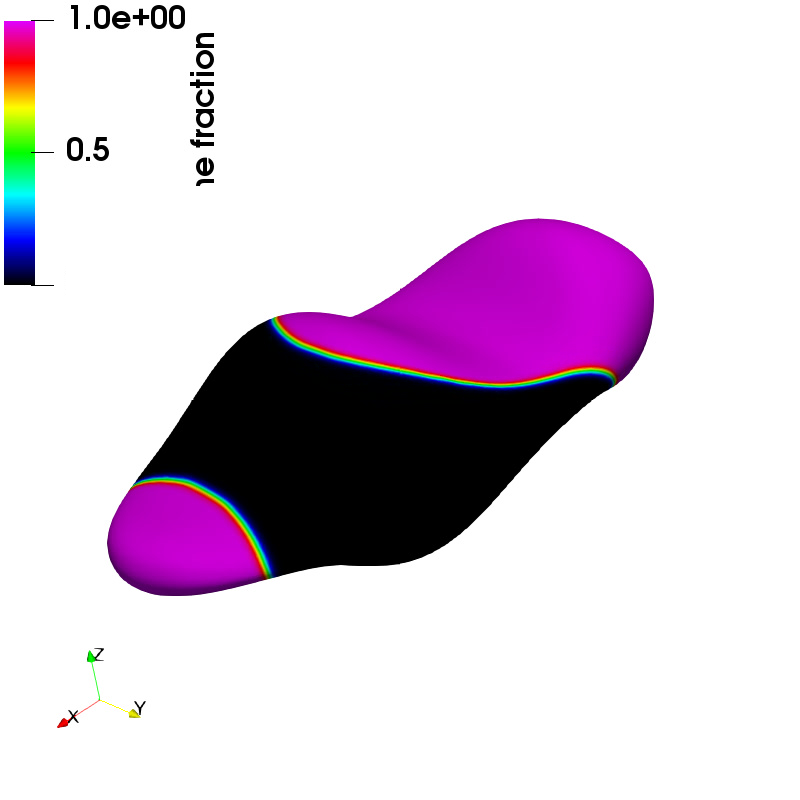}
        \put(35,80){\small{$t = 25000$}}
\end{overpic}
\begin{overpic}[width=.18\textwidth,viewport=110 150 700 580, clip,grid=false]{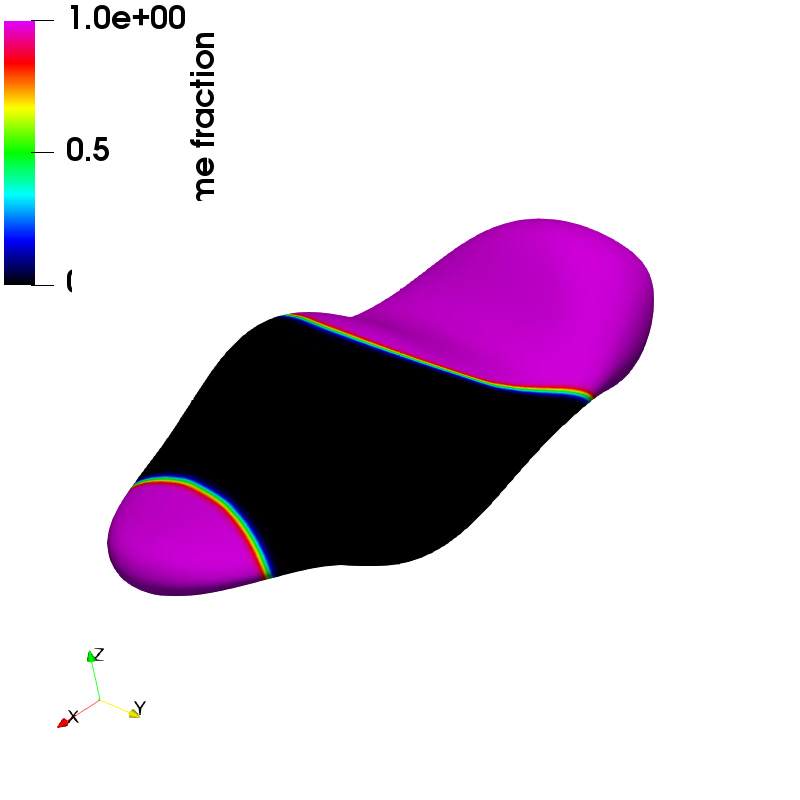}
        \put(35,80){\small{$t = 30000$}}
\end{overpic}
\end{center}
\caption{\label{cellCH_20_300}
Evolution of the numerical solution of the Cahn--Hilliard equation on the idealized cell surface for $t \in [2000, 30000]$
computed on a mesh with size $h =0.031$ and with time step is $\Delta t =1$.
The legend is as in Fig.~\ref{cellCH_0_10}.}
\end{figure}

\begin{figure}
\begin{center}
\begin{overpic}[width=.28\textwidth,viewport=60 250 810 580, clip,grid=false]{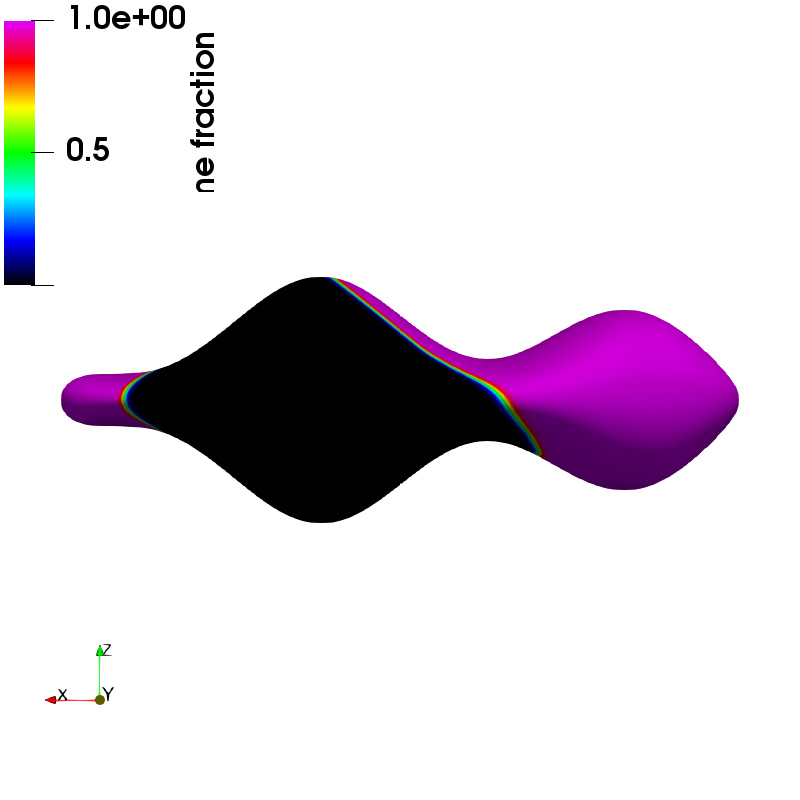}
        \put(25,45){\small{$t = 32000$}}
\end{overpic}
\begin{overpic}[width=.28\textwidth,viewport=60 250 810 580, clip,grid=false]{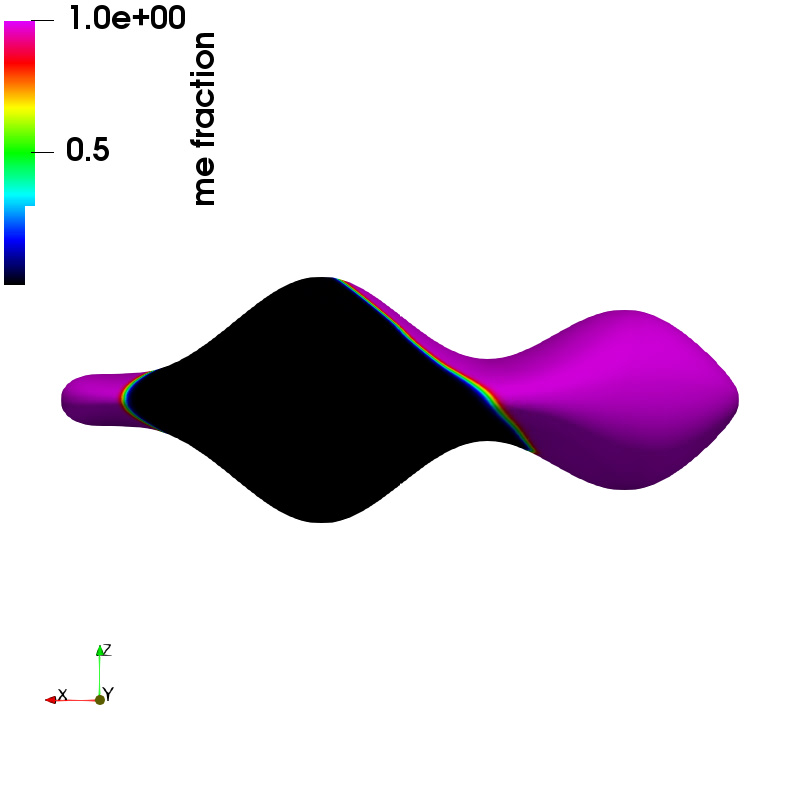}
        \put(25,45){\small{$t = 34000$}}
\end{overpic}
\begin{overpic}[width=.28\textwidth,viewport=60 250 810 580, clip,grid=false]{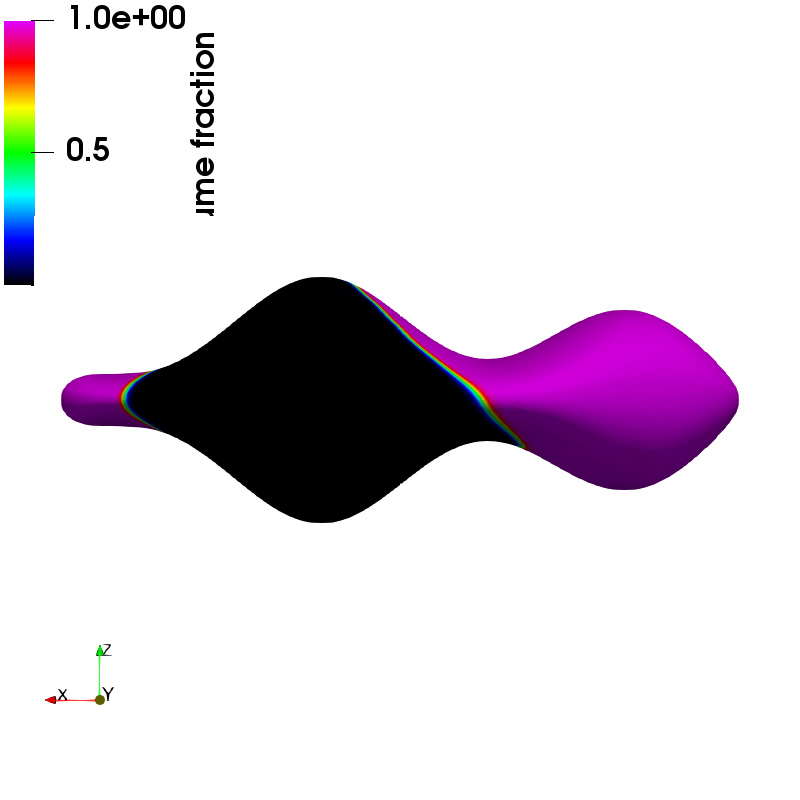}
        \put(25,45){\small{$t = 36000$}}
\end{overpic}
\hspace{.3cm}
\begin{overpic}[width=.08\textwidth,grid=false]{figures_AC_sphere_legend.png}
\end{overpic}
\end{center}
\caption{\label{cellCH_320_360}
Evolution of the numerical solution of the Cahn--Hilliard equation on the idealized cell surface for $t \in [32000, 36000]$
computed on a mesh with size $h =0.031$ and with time step is $\Delta t =1$.
Side view.}
\end{figure}

The evolution of the numerical solution
for $t \in [2000, 30000]$ is shown in Fig.~\ref{cellCH_20_300}.
After $t = 25000$, we observe a further deceleration in the process of dissipation of the interfacial energy.
This is even more evident in Fig.~\ref{cellCH_320_360}, which displays a side view
of the solution for $t \in [32000, 36000]$.
Since there is no visible change in the position of the interface
between phases after  $t = 36000$, we consider the rightmost panel in Fig.~\ref{cellCH_320_360}
close to the steady state. Just like for the sphere (see Sec.~\ref{sec:sphere_CH}), we see
one large and one small pink domain (i.e., $c = 1$): the small domain is positioned around the `beak'
and the large one is at the opposite end.

Finally, we note that although the shape of the surface clearly affects the steady
state of the phases, its influence on the spinodal decomposition, initial and intermediate
stages of phase separation is less evident from our results. This agrees with the limiting behavior
of the models: interface motion by  the geodesic curvature for the Allen--Cahn
equation~\cite{elliott2010modeling} and the minus  Laplace-Beltrami operator
of the geodesic curvature along the interface for the Cahn--Hilliard equation
~\cite{cahn1996cahn,o2016cahn}. Hence, the evolution of the interfaces between lateral phases
is largely driven by their intrinsic curvature and the interplay with the membrane shape is not explicit.

\section{Conclusions}\label{sec:concl}
We performed a computational study of lateral phase separation  and coarsening on surfaces of biological interest.
To model these processes, we considered both the Allen--Cahn (conservative model) and the Cahn--Hilliard
(non-conservative model) equations posed on surfaces. Although in this work we assume
the surfaces to be rigid, our longer term goal is to solve numerically the Allen--Cahn and Cahn-Hilliard equation
on evolving shapes. In fact, biological membranes exhibit shape transitions and shape instabilities,
which need to be accounted for in a realistic model. This need
dictated our choice for the numerical approach.

We considered the trace finite element method (TraceFEM), which
was shown in \cite{lehrenfeld2018stabilized} to handle relatively easily evolving surfaces.
TraceFEM is a geometrically unfitted method that has an important
additional advantage: surfaces can be defined implicitly and
no knowledge of the surface parametrization is required.  This allows flexible numerical treatment of complex shapes,
like the ones found in cell biology.

After validating the accuracy of our implementation of TraceFEM with benchmark problems,
we applied it to simulate phase transition on a series of surfaces of increasing
geometric complexity using both the surface Allen--Cahn and Cahn--Hilliard models.
We compared the numerical results produced by the two models
on a sphere and found that the Cahn--Hilliard model
successfully reproduces the spinodal decomposition experimentally
observed in giant vesicles in \cite{VEATCH20033074}.
Both models were also compared on the surface of a spindle
with the aim of getting a preliminary insight into the formation of microdomains
in bacteria \cite{Bramkamp_Lopez2015}. Finally, we presented the results on
a more complex surface that represents an idealized cell. For both the sphere and the idealized cell,
we let the simulations run until sufficiently close to the steady state to understand
the role of certain geometric characteristics on the final equilibrium.

\bibliographystyle{siam}
\bibliography{literatur}{}

\end{document}